%% file: 24RoughVolJump1108.tex
\documentclass[10pt]{article}
\usepackage{indentfirst, latexsym,bm}
\usepackage{amsmath}
\usepackage{pifont}
\usepackage{amsfonts}
\usepackage{mathrsfs}
\usepackage{array}
\usepackage{multirow} 
\usepackage{graphicx}
\usepackage{subfigure}
\usepackage{picinpar}
\RequirePackage[colorlinks,citecolor=blue,urlcolor=blue,linkcolor=blue]{hyperref} 

\usepackage{authblk}

\RequirePackage[numbers]{natbib}
\usepackage{enumitem} 

\usepackage{booktabs} 
\usepackage{threeparttable}
\usepackage{mathrsfs,amsfonts,amsmath,amssymb} 
\usepackage{color}

\makeatletter
\def\namedlabel#1#2{\begingroup
	#2%
	\def\@currentlabel{#2}%
	\phantomsection\label{#1}\endgroup
}
\makeatother

\newcommand{\blue}{\color{blue}}

\newcommand\email[1]{\href{mailto:#1}{ \nolinkurl{#1}}}

\newtheorem{theorem}{Theorem}[section]
\newtheorem{definition}[theorem]{Definition}
\newtheorem{lemma}[theorem]{Lemma}
\newtheorem{corollary}[theorem]{Corollary}
\newtheorem{proposition}[theorem]{Proposition}
\newtheorem{remark}[theorem]{Remark}
\newtheorem{condition}[theorem]{Condition}
\newtheorem{example}{Example}[section]

\def\blue{\color{blue}}

\def\blemma{\begin{lemma}}
	\def\elemma{\end{lemma}}
\def\bproposition{\begin{proposition}}
	\def\eproposition{\end{proposition}}
\def\ttheorem{\begin{theorem}}
	\def\etheorem{\end{theorem}}
\def\bcorollary{\begin{corollary}}
	\def\ecorollary{\end{corollary}}
\def\bremark{\begin{remark}}
	\def\eremark{\end{remark}}
\def\bcondition{\begin{condition}}
	\def\econdition{\end{condition}}

\def\benumerate{\begin{enumerate}}
	\def\eenumerate{\end{enumerate}}
\def\bitemize{\begin{itemize}}
	\def\eitemize{\end{itemize}}

\def\beqlb{\begin{eqnarray}}
	\def\eeqlb{\end{eqnarray}}
\def\beqnn{\begin{eqnarray*}}
	\def\eeqnn{\end{eqnarray*}}
\def\ar{\!\!\!&}

\def\proof{\noindent{\it Proof.~~}}
\def\qed{\hfill$\Box$\medskip}

\begin{document} 

 \title{\bf Path-dependent Fractional Volterra Equations and the Microstructure of Rough Volatility Models driven by Poisson Random Measures}

 \author{Ulrich Horst\footnote{Department of Mathematics, and School of Business and Economics,  Humboldt-Universit\"at zu Berlin, Unter den Linden 6, 10099 Berlin; email: horst@math.hu-berlin.de. Horst gratefully acknowledges financial support from the DFG CRC/TRR 388 ``Rough Analysis, Stochastic Dynamics and Related Fields", Project B02.} \quad\   Wei Xu\footnote{School of Mathematics and Statistics, Beijing Institute of Technology, No. 5, South Street, Zhongguancun, Haidian District, 100081 Beijing; email: xuwei.math@gmail.com}
 \quad\  and \quad Rouyi Zhang\footnote{Department of Mathematics, Humboldt-Universit\"at zu Berlin, Unter den Linden 6, 10099 Berlin; email: rouyi.zhang@hu-berlin.de}
    }
  \maketitle

\begin{abstract}
   We consider a microstructure foundation for rough volatility models driven by Poisson random measures. In our model the volatility is driven by self-exciting arrivals of market orders as well as self-exciting arrivals of limit orders and cancellations. The impact of market order on future order arrivals is captured by a Hawkes kernel with power law decay, and is hence persistent. The impact of limit orders on future order arrivals is temporary, yet possibly long-lived. After suitable scaling the volatility process converges to a fractional Heston model driven by an additional Poisson random measure. The random measure generates occasional spikes and clusters of spikes in the volatility process. 
   Our results are based on novel existence and uniqueness of solutions results for stochastic path-dependent Volterra equations driven by Poisson random measures.   
 	\medskip
 	\smallskip
 	
 	\noindent \textbf{\textit{MSC2020 subject classifications.}} 
 	 Primary 60G55, 60F05;  secondary 60G22.
 	 
 	  \smallskip
 	 
 \end{abstract}

 \input{1.IntroductionNoPlots}

\input{2.Benchmark}

 \input{3.Ctightness}

 \input{4.AccumulationNew}

 \input{RegularityMaximal}

 \input{NolinearVolterraRiccatiNew}

 
 

\bibliographystyle{plain}
\bibliography{Reference}

 \end{document}

%% file: 1.IntroductionNoPlots.tex
 \section{Introduction}
 \label{sec:Introduction}
 \setcounter{equation}{0}

 We establish a scaling limit for a family of order driven financial market models where the arrivals of market and limit orders are governed by a Hawkes process with heavy-tailed kernel. Both market and limit order arrivals have a self-exciting - albeit very different - impact on the arrival of future orders. In our model, the rescaled volatility processes converge in law to a rough Heston model driven by an additional Poisson random measure. 
 Our scaling limit is based on a novel existence and uniqueness of solutions results for non-linear, path-dependent fractional Volterra equations and on novel C-tightness results for families of stochastic integral equations driven by Poisson random measures.    

 \subsection{Stochastic volatility and Hawkes processes}

 Stochastic volatility models have been extensively investigated in the mathematical finance literature in the last decades. 
 The classical Heston \cite{Heston1993} model assumes that the variance (squared volatility) process follows square-root mean-reverting Cox-Ingerson-Ross \cite{CoxIngersollRoss1985} process. 
 The Heston model introduces a dynamics for the underlying asset that can take into account the asymmetry and excess kurtosis that are typically observed in financial assets returns, and provides analytically tractable option pricing formulas.  However, it is unable to capture large volatility movements. 
 To account for large volatility movements, the model has been extended to jump-diffusion settings by numerous authors. 
 Bates \cite{Bates1996} adds a jump component in the asset price process. Barndorff-Nielsen and Shephard \cite{Barndorff-NielsenShephard2001} consider volatility processes of Ornstein-Uhlenbeck type driven by L\'evy processes. 
 Affine models allowing for jumps in prices and volatilities are considered in Bakshi et al.~\cite{BakshiCaoChen1997}, Duffie et al. \cite{DuffiePanSingleton2000} and Pan \cite{Pan2002}, among many others. Empirical evidence for the presence of (negatively correlated) co-jumps in returns and volatility is given in, e.g.~Eraker  \cite{Eraker2004}, Eraker et al \cite{ErakerJohannesPolson2003} and Jacod and Todorov \cite{JacodTodorov2010}.

\subsubsection{Diffusion models with self-exciting jump dynamics}

 In a standard jump model with arrival rates calibrated to historical data, jumps are inherently rare. 
 Even more unlikely are patterns of multiple jumps in close succession over hours or days. Large moves, however, tend to appear in clusters. 
 For instance, A\"{\i}t-Sahalia et al.~\cite{AitSahaliaCacho-DiazLaeven2015} report that ``from mid-September to mid-November 2008, the US stock market jumped by more than 5\% on 16 separate days. Intraday fluctuations were even more pronounced: during the same two months, the range of intraday returns exceeded 10\% during 14 days.'' 
 
 Jump clusters have been discussed in the financial econometrics literature by Andersen et al. \cite{AndersenFusariTodorov2015} and Bates \cite{Bates2019} among many others. They consider continuous-time models of self-exciting price/volatility co-jumps in intraday stock returns and volatility and show that every small intraday jump substantially increases the probability of more intraday jumps in volatility and returns. Bates \cite{Bates2019} furthermore finds that ``multi-factor models with both exogenous and self-exciting but short-lived volatility spikes'' substantially improve model fits both in-sample and out-of-sample. He also shows that such models provide more accurate predictions of implied volatility. 

 To account for self-exiting jump dynamics one needs to leave the widely applied class of L\'evy jump processes. L\'evy processes have independent increments and hence do not allow for any type of serial dependence.  Hawkes processes are capable of displaying mutually exciting jumps.  A Hawkes process is a random point process $\{N(t):t\geq 0\}$ where events arrive at random points in time $\tau_1 < \tau_2 < \tau_3 < \cdots$ according to an  intensity process $\{V(t):t\geq 0\}$ that is usually of the form
 \beqlb\label{HawkesDensity}
     V(t):= \mu(t)+ \sum_{0<\tau_i<t} \phi(t-\tau_i) = \mu(t) + \int_{(0,t)} \phi(t-s)N(ds), \quad t \geq 0,
 \eeqlb
 where the \textit{immigration density} $\mu(\cdot)$ captures the arrivals of exogenous events, and the \textit{kernel} $\phi(\cdot)$ captures the self-exciting impact of past events on the arrivals of future events. 
 
 Originally introduced by Hawkes \cite{Hawkes1971a, Hawkes1971b} to model the occurrence seismic events, Hawkes processes have recently received considerable attention in the financial mathematics and economics literature as a powerful tool to model financial time series. Applications in finance range from intraday transaction dynamics \cite{Bowsher2007} to asset price and limit order book dynamics \cite{BacryDelattreHoffmannMuzy2013,HorstXu2019} and stochastic volatility modeling \cite{ElEuchFukasawaRosenbaum2018,HorstXu2022,JaissonRosenbaum2015,JaissonRosenbaum2016}. Hawkes processes with light-tailed kernels of the form
 \[
 	\phi(t) = \alpha e^{-\beta t}, \quad \alpha>0, \, \beta > 0
 \]
  have been used in Horst and Xu \cite{HorstXu2022} to establish scaling limits for a class of continuous-time stochastic volatility models with self-exciting jump dynamics. Many of the existing stochastic volatility models including the classical Heston model \cite{Heston1993}, the Heston model with jumps \cite{Bates1996, DuffiePanSingleton2000, Pan2002}, the OU-type volatility model \cite{Barndorff-NielsenShephard2001} and the multi-factor model with self-exciting volatility spikes  \cite{Bates2019}  were obtained as scaling limits under different scaling regimes. 
 
 \subsubsection{Rough volatility models}
 
Whereas most of  the financial economics literature focuses on semi-martingale volatility models with jumps, especially self-exciting jumps, the quantitative finance literature focuses mostly on rough volatility models. The analysis by Gatheral et al \cite{GatheralJaissonRosenbaum2018} and many others suggests that historical volatility time series are much rougher than those of Brownian martingales and that log-volatility is better modeled by a fractional Brownian motion with a Hurst parameter $H < 1/2$.  Although estimating the precise degree of roughness of volatility is subtle and challenging empirically, see \cite{ChongHoffmannLiuRosenbaumSzymanski2023a, ChongHoffmannLiuRosenbaumSzymanski2023b,ContDas2022,Fukasawa2021} and references therein, the use of rough volatility models is by now an established paradigm in the quantitative finance literature when modeling equity markets and pricing options. Rough volatility provides excellent fits to market data; it reproduces very good behavior of the implied volatility surfaces, in particular the at-the-money skew as shown in, e.g. \cite{BayerFrizGatheral2016}. We refer to the recent book by Bayer et al.~\cite{Bayer-et-al2023} for a comprehensive discussion of the many advantages of rough volatility models. 

The prototype rough volatility model is the rough Heston model. Other popular rough volatility models include the (mixed) rough Bergomi model \cite{JacquierMartiniMuguruza2018} and the rough SABR model \cite{FukasawaGatheral2022}. The rough Heston model can be described by an affine Volterra process \cite{AbiJaberLarssonPulido2019} and admits a semi-explicit representation of the characteristic function in terms of a fractional Riccati equation. Microstructure foundations for the  rough Heston models based on self-exciting market order dynamics were first considered by Rosenbaum and co-workers \cite{ElEuchFukasawaRosenbaum2018,JaissonRosenbaum2016}. They assert the weak convergence of the rescaled {\sl integrated} intensity processes  of nearly unstable Hawkes processes with heavy-tailed kernels of the form
\begin{equation} \label{kernel1}
	\phi(t) =  \big(1+  t \big)^{-\alpha-1}, \quad  \alpha \in  \Big(\frac 1 2, 1\Big)
\end{equation}
to the integral of a rough fractional diffusion. Their results were refined in our recent work \cite{HorstXuZhang2023a}, where we proved  the  weak convergence of the rescaled intensity processes to a fractional diffusion model.\footnote{We emphasize that convergence of the processes cannot be inferred from the convergence of the integrated processes.} 

Despite the great popularity of stochastic volatility models with self-exciting jump dynamics in the financial economics community and the great popularity of rough stochastic volatility models in the quantitative finance community, no effort has so far been made to unify both approaches. This paper provides a first step towards a common mathematical framework within which to study rough volatility models with self-exciting jump dynamics and their scaling limits.

\subsection{Our contributions}

Scaling self-exciting jumps into rough volatility models is mathematically extremely challenging. To overcome this challenge we follow a slightly different approach. Instead of working with ``true'' jumps, our model features occasional spikes and clusters of spikes in the volatility process in the spirit of the ``self-exciting [...] volatility spikes” observed in \cite{Bates2019}. 

Specifically, we consider a family of financial market models where orders to buy or sell an asset arrive according to Hawkes process where the impact of market orders on future order arrivals is captured by a power-law kernel of the form \eqref{kernel1} as in \cite{ElEuchFukasawaRosenbaum2018,HorstXuZhang2023a,JaissonRosenbaum2016}. Extending previous work on the microstructure of rough volatility the intensity is additionally driven by exogenous events that may or may not change asset prices. 

For the sake of exposition we shall think of the exogenous events as limit order placements and cancellations that change the state of the order book, say the volume imbalance at the top of the book.\footnote{There is a substantial economic literature that shows that volume imbalances at the top of the order book have an important impact of market order arrivals; see e.g.~\cite{Cebiroglu-Horst} and references therein.}
Many other interpretations are possible, though. We may also think of different types of markets orders (that would change prices) sent by different traders with different forms of impact on future order arrivals, of self-exciting signals to buy or sell an asset, or simply of exogenous news events that increase price volatility.    

We assume that the exogenous events arrive at a much higher frequency than market orders but that their impact on the volatility is much smaller. If we think of the exogenous events as limit orders, then this assumption is empirically well justified.\footnote{Similar assumptions have been make in many limit order book models;  see \cite{HorstKreher2018,HorstPaulsen2017} and references therein.} Our key assumption is that the impact of the exogenous events lasts for a random amount of time with a heavy-tailed life-time distribution before it dies.


\subsubsection{Our model}

The arrival of market orders is described by by an increasing sequence of adapted random times $\{\tau^\mathtt{m}_k\}_{k \in \mathbb{N}}$ whereas the arrival of limit orders is described by an increasing sequence of adapted random times $\{\tau^\mathtt{l}_k\}_{k \in \mathbb{N}}$; their respective life-lengths are described by a sequence of i.i.d.~positive random variables $\{\ell_k\}_{k \in \mathbb{N}}$ with distribution $\nu(dy)$.
Our intensity process then takes the form     
\begin{equation} \label{def:V2}
\begin{split}
 V(t)
 &:= \mu+ \Lambda(t)+ \sum_{k=1}^{N^{\mathtt{m}}(t)}\zeta^{\mathtt{m}} \cdot \phi(t-\tau^{\mathtt{m}}_k) + \sum_{k=1}^{N^{\mathtt{l}}(t)} \zeta^{\mathtt{l}} \cdot \mathbf{1}_{\{\ell_k>t-\tau^{\mathtt{l}}_k\}}\\
 &=
 \mu+ \Lambda(t)+ \int_0^t\int_0 \zeta^{\mathtt{m}} \cdot \phi(t-s)N^{\mathtt{m}}(ds)+ \int_0^t \int_0^\infty\zeta^{\mathtt{l}} \cdot \mathbf{1}_{\{y>t-s\}}N^{\mathtt{l}}(ds,dy),\quad t\geq 0,
 \end{split}
 \end{equation}
where $N^{\mathtt{m}}(ds)$ is a Hawkes process, and $N^{\mathtt{l}}(ds,dy)$ is a marked Hawkes process with respective intensities 
 \beqnn
	\lambda^{\mathtt{m}}\cdot V(s-)\, ds \quad \mbox{and} \quad 
	\lambda^{\mathtt{l}} \cdot V(s-)\, ds\,\nu(dy). 
 \eeqnn

Our main contribution is to prove that a sequence of suitably rescaled volatility processes $\{V^{(n)}\}_{n \geq 1}$ converges in distribution to the unique weak solution of a stochastic Volterra equation of the form
 \begin{equation} \label{I1}  
 \begin{split}
 V_*(t)=\,
 & V(0) \big(1-F(t)\big)+
 \int_0^t  c_1 \cdot f(s)ds+
 \int_0^t c_2 \cdot f(t-s)\sqrt{V_*(s)}dB(s)\\  
 &+\int_0^t\int_0^{\infty}\int_0^{V_*(s)}\Big(c_3 \cdot \int_{(t-s-y)^+}^{t-s} f(r)dr\Big) \widetilde{N}(ds,dy,dz),
 \end{split}
 \end{equation}
where $B$ is a Brownian motion, $\widetilde{N}(ds,dy,dz)$ is a compensated Poisson random measure with intensity 
 \beqnn
	c_4\cdot ds \,\nu_*(dy)\, dz\quad \mbox{and} \quad 
	\nu_*(dy)=\alpha(1+\alpha)y^{-\alpha-2}dy 
 \eeqnn
that captures the impact of limit orders on the volatility, $F$ and $f$ denote the Mittag-Leffler distribution and density function, and $c_1, c_2, c_3,c_4$ are scaling constants of which three can be chosen independently. The case $c_3 = 0$ or $c_4=0$ corresponds to the model studied in \cite{HorstXuZhang2023a}.

 \subsubsection{Mathematical contributions}

 Several mathematical challenges are to be overcome to establish the convergence of the sequence of rescaled volatility processes $\{V^{(n)}\}_{n \geq 1}$ to a rough Heston model with spikes. 
 First,  as already argued in \cite{JaissonRosenbaum2015} one of the main challenges when analyzing scaling limits of Hawkes processes with heavy-tailed kernels is to prove the $C$-tightness of the sequence of the intensity processes. 
 This problem has been overcome in the benchmark model \cite{HorstXuZhang2023a} by introducing a novel technique to verify the $C$-tightness of a sequence {\sl c\`adl\`ag} processes based on the classical Kolmogorov-Chentsov tightness criterion for {\sl continuous} processes. 
 We extend this method to account for the presence of stochastic integrals driven by Poisson random measures. 
Our key observation is that the integral processes corresponding to the limit order arrivals decomposes into a sequence of continuous processes that turns out to be $C$-tight plus a sequence of discontinuous ``remainder'' terms that converge to zero in $L^2$. 
 It is not difficult to show that the former sequence satisfies the classical Kolmogorov-Centsov criterion. 
 The key is to prove that the remainder terms satisfy the assumptions in \cite{HorstXuZhang2023a} and are hence $C$-tight. 
 With the $C$-tightness in hand, we conclude the convergence to zero in probability from the previously established $L^2$-convergence.  

 The second challenge is to identify the weak accumulation points. 
 We generalize the weak convergence result established in \cite[Section 4.2]{Xu2024b} for stochastic Volterra integrals by applying the general theory of weak convergence of It\^o's stochastic integrals with respect to infinite-dimensional semimartingales, due to Kurz and Protter \cite{KurtzProtter1996} to characterize the accumulation points. 
 The challenge is the time-dependence of the integrator of the integral w.r.t.~the marked Hawkes point process $N^{\mathtt{l}}(ds,dy)$ in \eqref{def:V2}, which prevents us from rewriting the rescaled stochastic integral as an integral with respect to some $(L^2)^\#$-martingale in the sense of \cite{KurtzProtter1996}. The problem can be overcome by utilizing the previously established $C$-tightness of the rescaled volatility processes. 
 This allows us to focus on the finite-dimensional distributions of the integral process, which in turns allows us to ``drop'' the time-dependence of the integrand. As in our accompanying work \cite{HorstXuZhang2023a} the $C$-tightness of the sequence of rescaled volatility processes is key to carry out this step. 
 We are unaware of any method to identify a candidate scaling limit without a priori knowing that the rescaled processes $\{V^{(n)}\}_{n \geq 1}$ are $C$-tight. 

 The third - and main - challenge is to prove weak uniqueness of accumulation points. 
It is well known that the rough Heston model can be described in terms of an affine Volterra process that admits a semi-explicit representation of the characteristic function in terms of a fractional Riccati equation; see  \cite{AbiJaberLarssonPulido2019} for details.  Our model is much more involved. 
 We prove that the characteristic function of the accumulation points can be represented in terms of the unique solution to a non-linear fractional Volterra equation that features a form of path-dependence, which originates from the random life-lengths of limit order impacts. Specifically, we prove that for any $\lambda\geq 0$ and function $g\in L^\infty(\mathbb{R}_+;\mathbb{R}_+)$ it holds that
 \begin{align}\label{charF1}
    \mathbf{E}\Big[\exp\{-\lambda\cdot V_*(T)-g*V_*(T)\}\Big]=
    \exp\Big\{-V_*(0)\cdot L_K*\psi^\lambda_{g} (T)-c_1*\psi^\lambda_{g} (T)\Big\},\quad T\geq 0,
 \end{align}
 where $\psi^\lambda_{g} $ is the unique non-negative solution to a non-linear fractional Volterra equation of the form  
  \begin{align*}
 D^{\alpha} \psi^\lambda_{g} +\psi^\lambda_{g}  =
  g-
 \frac{|c_2|^2}{2} \cdot  | \psi^\lambda_{g} |^2- \mathcal{V}\circ\psi^\lambda_{g} \quad \mbox{with}\quad
 I^{1-\alpha} \psi_g^\lambda(0+)=   \lambda .
 \end{align*}
%
 Here, $D^\alpha$ denotes the Riemann-Liouville fractional derivative and $I^{1-\alpha}$ the fractional integral operator,  and $\mathcal{V}$ is a nonlinear operator acting on a locally integrable function $f$ according to
 \beqnn
 \mathcal{V}\circ f(x)\:=
 \int_0^\infty \Big( \exp\Big\{-\int_{(x-y)^+}^x c_3 \cdot f(r)dr\Big\}-1+ \int_{(x-y)^+}^x c_3 \cdot f(r)dr \Big)\cdot c_4\cdot\nu_*(dy),\quad x\geq 0.
 \eeqnn

 Our uniqueness of solutions result for the above Volterra equation allows us to establish the characteristic function formula \eqref{charF1} for the accumulation points and hence their weak uniqueness and the weak convergence of our sequence of rescaled volatility processes. 
 Establishing the characteristic functional formula requires a series a priori estimates for the accumulation points, including the H\"older continuity of their sample paths.  

 Fractional Riccati and Volterra equations naturally arise rough volatility models and have been studied by many authors. Affine Volterra processes were first studied in \cite{AbiJaberLarssonPulido2019} where explicit exponential-affine representations of the Fourier–Laplace functional in terms of the unique solution of an associated system of deterministic integral equation of convolution type were provided. 
 They have been extended to affine Volterra processes with jumps in \cite{BondiLivieriPulido2024}. A family of fractional Riccati equations whose solutions take the form of power series is analyzed in \cite{CallegaroGrasselliPages2021}. Our equation is very different, due to the non-linear operator ${\cal V}$. To prove the existence of a unique solution to our Volterra equation we first prove that any solution is non-negative, continuous and then apply results established in aforementioned works, especially the power-series expansion establish in \cite{CallegaroGrasselliPages2021} to prove a priori estimate for the solutions near the origin. 
 Due to the continuity of the solutions the priori estimate extends to the entire time interval. 
 This provides us with a candidate function space within which to search for solutions. 
 Existence and uniqueness of solutions is then established using a fixed-point argument. 

 The remainder of this paper is organized as follows. In Section \ref{sec:model}, we introduce our benchmark model and state the scaling result. Section \ref{sec:tightness} establishes the $C$-tightness of the family of rescaled volatility processes. Section \ref{sec:accummulation} characterizes its weak accumulation points. Section \ref{Sec.Regular} establishes regularity conditions on the weak accumulation points that are key to the weak uniqueness of accumulation points and the weak convergence result that are established in Section \ref{Sec.VolterraRiccati}. 
 
 \medskip
 
 {\sl Notation.} We frequently use the following notation. For any $x\in\mathbb R$, we put $x^+:= x\vee 0$ and denote by $[x]$ be the integer part of $x$. 
By $f*g$ we denote the convolution of functions $f,g$ on $\mathbb R_+$. By $\Delta_h$ and $\nabla_h$ we denote the forward and backward difference operator with step size $h>0$, i.e.,
\[
\Delta_h f(x):= f(x+h)-f(x)\quad \mbox{and}\quad \nabla_h f(x):= f(x)-f(x-h).
\]

%% file: 2.Benchmark.tex
 \section{The model}
 \label{sec:model}
 \setcounter{equation}{0}
 
 In this section, we introduce a stochastic volatility model that generalizes the rough volatility models studied in \cite{ElEuchFukasawaRosenbaum2018,HorstXuZhang2023a,JaissonRosenbaum2016} to a jump-type regime, and study its scaling limit. 
 We assume throughout that all random variables and stochastic processes are defined on a common probability space $(\Omega, \mathscr{F} , \mathbf{P})$ endowed with a filtration $\{\mathscr{F}_t : t \geq 0\}$ that satisfies the usual hypotheses. 
 The convergence concept for stochastic processes we use will be weak convergence in the space $\mathbf{C}(\mathbb{R}_+;\mathbb{R}^d)$ of all $\mathbb{R}^d$-valued continuous functions on $\mathbb{R}_+$ endowed with the uniform topology or in the space $\mathbf{D}(\mathbb{R}_+;\mathbb{R}^d)$ of all $\mathbb{R}^d$-valued c\`{a}dl\`{a}g functions on $\mathbb{R}_+$ endowed with the Skorokhod topology; see \cite{Billingsley1999,JacodShiryaev2003} for details.

 \subsection{The benchmark model}

 We consider an order-driven market where asset prices are driven by incoming orders to buy or sell an asset. 
 The market order arrival times are described by an increasing sequence of adapted random times $\{\tau^\mathtt{m}_k\}_{k\geq1}$. 
 The impact of each market order on the price is described by an independent sequence of i.i.d.~$\mathbb R$-valued random variables $\{\xi^\mathtt{m}_k \}_{k\geq 1}$ with distribution $\nu^\mathtt{m}(du)$.  
 In terms of these sequences we define the random point measure
 \beqnn
 N^\mathtt{m}(ds, du) := \sum_{k=1}^\infty {\bf 1}_{\{\tau^\mathtt{m}_k \in ds, \xi^\mathtt{m}_k \in du\}}
 \eeqnn
 on $(0,\infty) \times \mathbb R$ and assume that the logarithmic price process $\{P(t) : t \geq 0\}$ satisfies the dynamics 
 \beqlb\label{eqn.Price}
 P(t) \ar=\ar  P(0)+   \sum_{\tau^\mathtt{m}_k\leq t}  \xi^\mathtt{m}_k = P(0)+  \int_0^t \int_{\mathbb R}   u \, N^\mathtt{m}(ds,du).
 \eeqlb 
 
 We assume that the measure $N^\mathtt{m}(ds,du)$ is a \textsl{marked Hawkes point measure}, that is, the embedded point process 
 \beqnn \label{Nm}
 N^\mathtt{m}(t) := N^\mathtt{m}\big((0,t],\mathbb{R}\big),\quad t\geq 0,
 \eeqnn
 is a Hawkes process with an intensity process $\{V(t) : t \geq 0\}$ that will be specified below. 
 In \cite{ElEuchFukasawaRosenbaum2018,HorstXuZhang2023a,JaissonRosenbaum2016}  it is assumed that the intensity process is of the form 
 \beqlb\label{eqn.Vol1}
 	V(t):= \mu + \Lambda(t) + \sum_{k=1}^{N^\mathtt{m}(t)}\zeta^\mathtt{m}  \cdot \phi(t-\tau_k) , 
 \eeqlb
 where the positive constant $\mu$ describes the arrival rate of exogenous orders, the non-negative function $\Lambda$ represents the combined impact of all the events that arrived prior to the time zero on future arrivals, the kernel 
 \beqlb\label{def:kernel}
 \phi(t):=\alpha \cdot \big(1+  t \big)^{-\alpha-1}, \quad t\geq 0, \  \alpha \in  \Big(\frac 1 2, 1\Big)
 \eeqlb
 specifies the self-exciting impact of past order arrivals on future arrivals, and the positive constant $\zeta^\mathtt{m}$ measures the impact of each child order on the overall order arrival dynamics. 

 In our recent work \cite{HorstXuZhang2023a} we have shown that the intensity process converges to a fractional diffusion and that the logarithmic price process converges to a rough Heston-type model after suitable scaling. In this paper we allow for a more general intensity that accounts for occasional spikes in the volatility process, akin to a series of self-exciting jumps in the volatility process. The occurrence of self-exciting jumps in volatility and/or price processes is well documented in the financial economics literature (see \cite{Bates2019} and references therein) and cannot be captured by standard rough volatility models. 
 

 \subsubsection{The intensity process}

 We assume that the intensity process is driven by incoming market orders as in \eqref{eqn.Vol1} with the power-law kernel \eqref{def:kernel} whose impact on future order arrivals decays slowly to zero but never completely vanishes, and additional exogenous events that do not change prices but still have an impact on asset price volatility. We may think of the exogenous events as self-exciting buying or selling signals that increase the intensity of market order arrivals or of self-exciting limit order placements or cancellations that change the buy-sell side imbalance of order book volumes.\footnote{There is a substantial economic literature that shows that volume imbalances at the top of the order book have an important impact of market order arrivals; see e.g.~\cite{Cebiroglu-Horst} and references therein.} 

 For the sake of exposition we will indeed think of exogenous events as limit order placements or cancellations (negative placements); many other interpretations are possible, though. We assume that limit orders/cancellations arrive at a much higher frequency than market orders but that their impact on the volatility and hence the price dynamics is much smaller. 
Both assumptions are empirically well justified. 

 We also assume that the impact of limit orders/cancellations lasts for a random amount of time before its death. This is a modeling assumption that awaits empirical justification. 
 We shall see that the different modeling assumptions - decaying but persistent impact vs.~constant but temporary impact - on market and limit order arrivals translate into very different dynamics of the intensity process. 

 The arrivals of limit orders are described by an increasing sequence of adapted random times $\{\tau^\mathtt{l}_k\}_{k \geq 1}$ and their respective life-lengths are described by a sequence of i.i.d. positive random variables $\{\ell_k\}_{k \geq 1}$ with distribution 
 \beqlb\label{def:nu}
 \nu(dy) =  (1+\alpha)\cdot (1+  y )^{-\alpha-2}\, dy,\quad \alpha \in  \Big(\frac 1 2, 1\Big). 
 \eeqlb
 In particular, the impact of certain orders will last very long. These orders will trigger a large number of child events, which, in the scaling limit will result in occasional spikes in the volatility process. Overall, our intensity process takes the form     
 \begin{equation} \label{def:V}
 \begin{split}
 V(t)
 &:= \mu + \Lambda(t) + \sum_{k=1}^{N^{\mathtt{m}}(t)}\zeta^{\mathtt{m}} \cdot \phi(t-\tau^{\mathtt{m}}_k) + \sum_{k=1}^{N^{\mathtt{l}}(t)} \zeta^{\mathtt{l}} \cdot \mathbf{1}_{\{\ell_k>t-\tau^{\mathtt{l}}_k\}}\\
 &=
 \mu + \Lambda(t) + \int_0^t  \zeta^{\mathtt{m}} \cdot \phi(t-s)N^{\mathtt{m}}(ds)+ \int_0^t \int_0^\infty\zeta^{\mathtt{l}} \cdot \mathbf{1}_{\{y>t-s\}}N^{\mathtt{l}}(ds,dy),\quad t\geq 0,
 \end{split}
 \end{equation}
 where $N^{\mathtt{m}}(ds)$ is a Hawkes process and $N^{\mathtt{l}}(ds,dy)$ is a marked Hawkes process on $(0,\infty)^2$ with respective intensities 
 \beqnn
	\lambda^{\mathtt{m}}\cdot V(s-)\, ds \quad \mbox{and} \quad 
	\lambda^{\mathtt{l}} \cdot \alpha \cdot V(s-)\, ds\, \nu(dy). 
 \eeqnn
 The positive scaling constants $\lambda^\mathtt{m}$ and $\lambda^\mathtt{l}$ will be further specified below. We call the following vector \textit{the characteristics} of our volatility process: 
 \beqlb \label{def:char}
	\big(\mu,\Lambda,\zeta^\mathtt m,\lambda^\mathtt m,\zeta^\mathtt l,\lambda^\mathtt l, \phi, \nu \big) .
 \eeqlb
 

\subsubsection{Representation via Poisson random measures}

 In what follows we represent the intensity process as an integral process driven by two martingale measures. 
Using the argument given in \cite{HorstXuZhang2023a}, on an extension of the original probability space, we can define two time-homogeneous Poisson random measures $N_\mathtt{m}(ds,dz)$ and $N_\mathtt{l}(ds,dy,dz)$ on $(0,\infty)^2$ and $(0,\infty)^3$ with respective intensities 
 \beqnn
	\lambda^\mathtt{m}\cdot ds\, dz \quad \mbox{and} \quad \lambda^\mathtt l\cdot \alpha \cdot ds\, \nu(dy)\, dz
 \eeqnn
such that the intensity process can be represented as
 \beqnn
 V(t)\ar=\ar \mu + \Lambda(t) +  \int_0^t  \int_0^{V(s-)} \zeta^\mathtt{m}\cdot \phi(t-s) N_\mathtt{m}(ds, dz)
  +  \int_0^t \int_0^\infty  \int_0^{V(s-)} \zeta^\mathtt{l} \cdot \mathbf{1}_{\{y> t-s\}}N_\mathtt{l}(ds,dy,dz).
 \eeqnn
 The respective compensated Poisson random measures are denoted as 
 \beqnn
 \widetilde{N}_\mathtt{m}(ds,dz):= N_\mathtt{m}(ds,dz)- \lambda^\mathtt{m}\cdot ds \, dz \quad \text{and}\quad \widetilde{N}_\mathtt{l}(ds,dy,dz):= N_\mathtt{l}(ds,dy,dz)- \lambda^\mathtt l\cdot \alpha  \cdot \, ds\, \nu(dy)\, dz.
 \eeqnn

 To represent the intensity process as an integral process driven by the $(\mathscr{F}_t)$-martingale measures $\widetilde{N}_\mathtt{m}(ds,dz)$ and $\widetilde{N}_\mathtt{l}(ds,dy,dz)$ we set 
 \beqnn
 \beta:= \zeta^\mathtt{m}\cdot \lambda^\mathtt{m}+ \zeta^\mathtt{l}\cdot \lambda^\mathtt{l}  
 \eeqnn
 and recall that the resolvent $R$ associated with the kernel $\phi$ is given as the unique solution of the resolvent equation
 \beqlb \label{RH}
 R(t)=\beta\cdot\phi(t)+ \beta\cdot\phi*R(t),\quad t\geq 0. 
 \eeqlb
 Furthermore, we introduce the two-parameter function 
 \beqlb\label{RRH}
 \mathcal{R}(t,y):= \mathbf{1}_{\{y>t\}}+ \int_0^t R(t-s) \mathbf{1}_{\{y>s\}}ds,\quad t\geq 0, \, y> 0,
 \eeqlb
 that represents the mean impact of an event with given life-length on the intensity process. In terms of these quantities the intensity process $V$ can be represented as the unique solution to the integral equation
 \begin{align}    \label{SVE1}  
 \begin{split}
 V(t)=& ~\mu + \mu\int_0^t R (s)ds +  \Lambda(t) +\int_0^t R(t-s)\Lambda(s)ds \cr
 & ~ + \int_0^t \int_0^{V(s-)} \frac{\zeta^\mathtt{m}}{\beta}\cdot R(t-s) \widetilde{N}_\mathtt{m}(ds,dz)
 \\
  & ~ +  \int_0^t \int_0^\infty  \int_0^{V(s-)}  \zeta^\mathtt{l}\cdot  \mathcal{R} (t-s,y)\widetilde{N}_\mathtt{l}(ds,dy,dz),\quad t\geq 0.
  \end{split}
 \end{align}
 

 \subsection{Scaling limit of the volatility process}

 We are now going to introduce a sequence of rescaled intensity processes where the intensity of order arrivals tends to infinity, the impact of an individual order on asset prices tends to zero and the average number of (market and limit) child events tends to one. 

 The dynamics of the $n$-th volatility process $\{V_n(t) : t \geq 0 \}$ is defined in terms of an underlying Hawkes process and an additional Hawkes point process akin to \eqref{def:V} with corresponding characteristics
 \beqnn
	\big(\mu_n,\Lambda_n,\zeta^\mathtt{m}_n,\lambda^\mathtt{m}_n,\zeta^\mathtt{l}_n,\lambda^\mathtt{l}_n, \phi, \nu \big).
 \eeqnn
 The corresponding Poisson random measures are denoted ${N}_{\mathtt{m},n}(ds,dz)$ and ${N}_{\mathtt{l},n}(ds,dy,dz)$, respectively. Our goal is to establish the convergence in law of the sequence of rescaled volatility processes defined by 
 \beqnn
 V^{(n)}(t) :=\frac{V_n(nt)}{n^{2 \alpha - 1}}, \quad t \geq 0,\, n\geq 1. 
 \eeqnn
 
 \begin{remark}
 Once the convergence of the above processes has been established the same arguments as in \cite{HorstXuZhang2023a} can be used to establish the joint convergence of the price-volatility process (after suitable scaling of the price process).  
 We hence solely focus on the volatility processes in what follows. 
 \end{remark}
 
 Applying a change of variables to \eqref{SVE1} shows that the rescaled volatility process satisfies
 \begin{align}\label{rescaledV}
 \begin{split}
 V^{(n)}(t)=& \frac{\mu_n}{n^{2\alpha-1}} + 
 \frac{\mu_n}{n^{2\alpha-1}} \int_0^{nt} R_n(s)ds +  \frac{\Lambda_n(nt)}{n^{2\alpha-1}} +\frac{1}{n^{2\alpha-2}}\int_0^t R_n\big(n(t-s)\big)\Lambda_n(ns)ds \\
 &+\frac{\zeta^\mathtt{m}_n}{\beta_n \cdot n^{2\alpha-1}}\int_0^t  \int_0^{V^{(n)}(s-)}  R_n\big(n(t-s)\big) \widetilde{N}^{(n)}_\mathtt{m}( ds , dz)\\   &+\frac{\zeta^\mathtt{l}_n}{ n^{2\alpha-1}}\int_0^t \int_0^\infty \int_0^{V^{(n)}(s-)} \mathcal{R}_n\big(n(t-s),ny\big)\widetilde{N}^{(n)}_\mathtt{l}(ds,dy,dz), \quad t\geq 0,
 \end{split}
 \end{align}
 where $R_n$ is the unique solution to the resolvent equation \eqref{RH} but with $\beta$ replaced by
 \beqnn
 \beta_n := \zeta^\mathtt{m}_n\cdot \lambda^\mathtt{m}_n+\zeta^\mathtt{l}_n  \cdot  \lambda^\mathtt{l}_n ,
 \eeqnn
 and the two-parameter function $\mathcal{R}_n$ is defined in terms of the resolvent $R_n$ as in \eqref{RRH}. Furthermore, 
 \begin{gather*}
   \widetilde{N}^{(n)}_\mathtt{m}(ds, dz)
   := 
   N_{\mathtt{m},n}(n\cdot ds ,n^{2\alpha-1}\cdot dz)- n^{2\alpha}\cdot \lambda^\mathtt{m}_n\cdot ds\, dz,\\ 
   \widetilde{N}^{(n)}_\mathtt{l}(ds,dy,dz)
   := 
   N_{\mathtt{l},n}(n\cdot ds, n\cdot dy,n^{2\alpha-1}\cdot dz)-n^{2\alpha}\cdot \lambda^\mathtt{l}_n\cdot \alpha \cdot ds\, \nu(n\cdot dy)\, dz . 
 \end{gather*}

 We assume throughout that our model parameters satisfy the following condition. 
 In particular, we assume that limit orders arrive at a much higher frequency than market orders but that they have a much smaller impact on the volatility of asset prices.  

 \begin{condition}\label{main.Condition}
 The kernel $\phi$ is given by \eqref{def:kernel} and the life-time distribution $\nu$ is given by \eqref{def:nu} for some $\frac 1 2 < \alpha < 1$. Moreover, the model parameters satisfy the following conditions.
  \begin{enumerate}
   \item[(1)] For each $n\geq 1$, we have that $\beta_n<1$ and the function $\Lambda_n$ is of the form 
   \beqnn
   \Lambda_n(t) =V_{0,n}\cdot (1+t)^{-\alpha} ,\quad t\geq 0. 
   \eeqnn 
   
  \item[(2)] There exist four non-negative constants  $\zeta_*^{\mathtt{m}},\lambda_*^{\mathtt{m}},\zeta_*^{\mathtt{l}},\lambda_*^{\mathtt{l}}$  such that
 \beqnn
 \zeta_*^{\mathtt{m}}\cdot\lambda_*^{\mathtt{m}}+\zeta_*^{\mathtt{l}}\cdot\lambda_*^{\mathtt{l}}=1
 \quad \mbox{and}\quad 
 \zeta^{\mathtt{m}}_n\to \zeta_*^\mathtt{m}, \quad 
 \lambda_n^{\mathtt{m}}\to\lambda_*^{\mathtt{m}},\quad \frac{\zeta_n^{\mathtt{l}}}{n^{\alpha-1}} \to \zeta_*^{\mathtt{l}},\quad 
 \frac{\lambda^{\mathtt{l}}_n}{n^{1-\alpha}}  \to \lambda_*^{\mathtt{l}} \quad \mbox{as }n\to\infty.
 \eeqnn  
 
 \item[(3)] There exist three constants $V_*(0) \in \mathbb{R}_+$, $a\geq 0$ and $b>0$ such that as $n\to\infty$, 
 \beqlb
  \frac{V_{0,n}}{n^{2\alpha-1}}\to V_*(0),\quad \frac{\mu_n}{n^{\alpha-1}}	\to a , \quad 
 n^\alpha(1-\beta_n)\to  b .
 \eeqlb
  \end{enumerate}
\end{condition}

 To formulate our main results, we first recall some notation and terminologies from the theory of fractional equations. 
 In what follows we denote by $\Gamma(\cdot)$ the Gamma function and set 
 \beqnn
 \gamma := \frac{b}{ \Gamma(1-\alpha)}. 
 \eeqnn
 The Mittag-Leffler probability density function $f^{\alpha, \gamma}$ and its distribution function $F^{\alpha, \gamma}$ with parameter $(\alpha, \gamma)$ are given, respectively, by 
 \begin{gather*}
 	f^{\alpha,\gamma}(t)=\gamma \cdot t^{\alpha-1}\cdot E_{\alpha,\alpha}(-\gamma\cdot t^\alpha) \quad \mbox{and} \quad
 	F^{\alpha,\gamma}(t)=\int_0^t f^{\alpha,\gamma}(s)ds, \quad t> 0,
 \end{gather*}
 where $E_{\alpha,\alpha}$ is the Mittag-Leffler function with parameter $(\alpha,\alpha)$ given by
 \begin{equation*}
 E_{\alpha,\alpha}(x):= \sum_{k=0}^{\infty}\frac{x^n}{\Gamma(\alpha(n+1))},\quad x\in \mathbb{R}.
 \end{equation*}  
 The function $E_{\alpha,\alpha}$ is locally H{\"o}lder continuous with index $\alpha$.
 Additionally, there exists a constant $C>0$ such that uniformly in $t>0$, 
 \begin{align}\label{upper bound of mittag leffler density functions}
 F^{\alpha,\gamma}(t)\leq  C\cdot t^{\alpha},\quad f^{\alpha,\gamma}(t)\leq C\cdot t^{\alpha-1}
 \quad\mbox{and}\quad 
 \big(f^{\alpha,\gamma}\big)^{'}(t)\leq C\cdot t^{\alpha-2}.   
 \end{align}  
 Moreover, we shall repeatedly use the following important functions 
\begin{equation}\label{func:K}
    K(t):= \gamma\cdot \frac{t^{\alpha-1}}{\Gamma(\alpha)} \quad \mbox{and} \quad 
    L_K(t) := \frac 1 \gamma \cdot \frac{t^{-\alpha}}{\Gamma(1-\alpha)},
    \quad t> 0. 
 \end{equation}
 The function $L_K$ is the \textsl{resolvent of the first kind}, and $f^{\alpha,\gamma}$ is the \textsl{resolvent of the second kind} for $K$, i.e.,
 \begin{equation} \label{LK}
 L_K* K\equiv 1 \quad \mbox{and} \quad K=f^{\alpha,\gamma}+f^{\alpha,\gamma}*K.
 \end{equation}

 We are now ready to state the main results of this paper. 
 We start with the following result on the existence of accumulation points of the sequence of rescaled volatility processes and their probabilistic representation. 
 
 \begin{theorem}\label{thm1}
 Under Condition \ref{main.Condition}, the following holds.
 \begin{enumerate}
 \item[(1)] The sequence of rescaled processes $\{V^{(n)}\}_{n\geq 1}$ is $C$-tight in  $\mathbf{D}(\mathbb{R}_+;\mathbb{R}_+)$. 
 	
 \item[(2)] Any accumulation point $V_* \in \mathbf{C}(\mathbb{R}_+;\mathbb{R}_+)$ is a weak solution to the stochastic Volterra equation
 \begin{align}    \label{eq1}
 \begin{split}
 V_*(t)=  & ~ V_*(0) \cdot \big(1-F^{\alpha,\gamma}(t)\big) + \frac{a}{b}\cdot F^{\alpha,\gamma}(t)
 + \int_0^t f^{\alpha,\gamma}(t-s)\cdot\frac{\zeta_*^\mathtt m\sqrt{\lambda^\mathtt{m}_*}}{b}\sqrt{V_*(s)}\,dB(s)\\  
 & \ +\int_0^t\int_0^{\infty}\int_0^{V_*(s)}\Big(\int_{(t-s-y)^+}^{t-s} \frac{\zeta_*^\mathtt{l}}{b}\cdot f^{\alpha,\gamma}(r)dr \Big) \widetilde{N}(ds,dy,dz),
 \end{split}
 \end{align}
 where $B$ is a Brownian motion and $\widetilde{N}(ds,dy,dz)$ is a compensated Poisson random measure with intensity 
 \beqnn
  \lambda^\mathtt{l}_* \cdot ds\,\nu_*(dy)\,dz \quad \mbox{and} \quad 
  \nu_*(dy)=\alpha(1+\alpha)y^{-\alpha-2}\, dy. 
 \eeqnn
 
 \item[(3)] The stochastic Volterra equation (\ref{eq1}) is equivalent to 
 \begin{align}\label{eq2}
 \begin{split}
 V_*(t)=  & ~ V_*(0)+ \int_0^tK(t-s)
 	\Big(\frac{a}{b}-V_*(s)\Big)ds+\int_0^t K(t-s)\,
 	\frac{\zeta_*^\mathtt m\sqrt{\lambda^\mathtt{m}_*}}{b}\sqrt{V_*(s)}dB(s) \\
 &+\int_0^t\int_0^\infty\int_0^{V_*(s)}\Big(\int_{(t-s-y)^+}^{t-s} \frac{\zeta_*^\mathtt{l}}{b}\cdot K(r)dr \Big) \widetilde{N}(ds,dy,dz). 
 \end{split}
 \end{align} 
 \end{enumerate} 
 \end{theorem}
 
 \begin{remark}
 	By Remark~2.2 in \cite{Xu2024b}, there exists a constant $C>0$ such that for any $t\geq 0$,
 	\beqlb\label{eqn.210}
 	\int_0^t\int_0^{\infty} \Big(\int_{(t-s-y)^+}^{t-s}  f^{\alpha,\gamma}(r)dr \Big)^2 \nu_*(dy)ds + \int_0^t\int_0^{\infty} \Big(\int_{(t-s-y)^+}^{t-s}   K(r)dr \Big)^2 \nu_*(dy)ds \leq C\cdot t^\alpha. 
 	\eeqlb 
 	Hence the two stochastic Volterra integrals with respect to the compensated Poisson random measure $\widetilde{N}(ds,dy,dz)$ in (\ref{eq1}) and (\ref{eq2}) are well-defined as It\^o's integrals. 
 \end{remark}
 
 The next result establishes important a-priori estimates for weak solutions of our stochastic Volterra equations (\ref{eq1}) and (\ref{eq2}). To state the result, we recall that for $\kappa\in(0,1]$, the $\kappa$-H\"older coefficient of a real-valued function $f$ on $[0,T]$ is defined by 
 \beqnn
 \big\|f\big\|_{C_T^{\kappa}}:= \sup_{0\leq x_1<x_2\leq T}\frac{|f(x_1)-f(x_2)|}{|x_1-x_2|^\kappa}. 
 \eeqnn
 
  \begin{theorem}\label{MainThm.02}
 Let $V_*$ be any weak solution to the stochastic Volterra equation (\ref{eq1}) or (\ref{eq2}). Then  the following hold.
 
 \begin{enumerate}
 	\item[(1)] The process $V_*$ is almost surely  H\"older continuous of any order strictly less that  $ \alpha-\frac{1}{2}$. 
 	
 	\item[(2)] For any $\kappa\in (0,\alpha-\frac{1}{2})$ and $p\geq 0$, there exists a constant $C>0$ such that for any $T\geq 0$,
 	\beqnn
 	 \mathbf{E}\Big[  \big\|V_{*}\big\|_{C^\kappa_T}^p \Big] \leq  C\cdot (1+ T)^{p(\alpha-\kappa)}.
 	\eeqnn
 	
 	\item[(3)] For each $p\geq 0$, there exists a constant $C>0$ such that for any $T>0$,
 	\beqnn 
 	\mathbf{E}\bigg[  \sup_{t\in[0,T]}\big| V_*(t) \big|^p \bigg] \leq  C (1+T)^{p\alpha}. 
 	\eeqnn 
 \end{enumerate}
 	
 \end{theorem}
 
The preceding theorem is key to prove the weak uniqueness of solutions to our stochastic Volterra equation. Along with the previously established results, weak uniqueness of solutions implies the convergence in law of the sequence of rescaled volatility processes to a unique limit. Specifically, we have the following result. 
 
 \begin{theorem} \label{thmcf}
 For any $\lambda\geq 0$ and $g\in L^\infty(\mathbb{R}_+;\mathbb{R}_+)$ the following hold. 
 \begin{enumerate}
 \item[(1)] There exists a unique (locally integrable) global solution  $\psi^\lambda_{g} \in \mathbf{C}((0,\infty);\mathbb{R}_+)$ to the nonlinear Volterra integral equation
 \begin{align} \label{PSI}   
 	\psi^\lambda_{g}  (t)=\lambda\cdot f^{\alpha,\gamma}(t)+g*f^{\alpha,\gamma}(t)-\frac 1 2 \cdot \Big(\frac{\zeta_*^\mathtt{m}\sqrt{\lambda^\mathtt{m}_*}}{b}\Big)^2 \cdot |\psi^\lambda_{g} |^2* f^{\alpha,\gamma}(t) -
 	\big(\mathcal{V}\circ\psi^\lambda_{g} \big)* f^{\alpha,\gamma}(t) ,
 \end{align}
 where $\mathcal{V}$ is a nonlinear operator that acts on a locally integrable function $f$ according to
 \beqnn
 \mathcal{V}\circ f(x) :=
 \int_0^\infty \Big(
 \exp\Big\{-\int_{(x-y)^+}^x \frac{\zeta_*^\mathtt{l}}{b}\cdot f(r)dr\Big\}-1+
 \int_{(x-y)^+}^x \frac{\zeta_*^\mathtt{l}}{b}\cdot f(r)dr
 \Big)\cdot  \lambda^\mathtt{l}_* \cdot \nu_*(dy),\quad x\geq 0.
 \eeqnn
 
 \item[(2)] The Laplace functional of any weak solution $V_*$ of our stochastic Volterra equation admits the representation
 	\begin{align} 	\label{CF}
 		\mathbf{E}\Big[\exp \big\{-\lambda\cdot V_*(T)-g*V_*(T) \big\}\Big]=
 		\exp\Big\{-V_*(0)\cdot L_K*\psi^\lambda_{g} (T)-\frac{a}{b}*\psi^\lambda_{g}(T)  \Big\},\quad T\geq 0.  
 	\end{align}
 In particular, the stochastic Volterra equations (\ref{eq1}) and (\ref{eq2}) admit a unique weak solution. 	
 
  \end{enumerate}
 	 
 \end{theorem}
   
 As an immediate corollary of the above theorems we obtain the convergence in law of the sequence of rescaled volatility processes. 

 \begin{corollary} 
 The sequence of rescaled volatility processes $\{V^n : n \in \mathbb N \}$ converges in law to the unique weak solution to \eqref{eq1}.
 \end{corollary}

 Three main challenges arise when analyzing the limiting dynamics of our volatility processes. The first is to prove the $C$-tightness of rescaled processes, which will be achieved in Section \ref{sec:tightness}. 
 The second is the  characterization of the weak accumulation points, especially the limiting contributions of the limit order flow to the volatility process. 
 The second challenge is addressed in Section \ref{sec:accummulation}. The third challenge is to establish the uniqueneess of weak accumulation points, that is, the uniqueness of solutions to the non-linear Volterra equation \eqref{PSI}. 
 This is achieved in Section \ref{Sec.VolterraRiccati}. 
 The uniqueness result strongly hinges on the regularity conditions established in Section  \ref{Sec.Regular}.

%% file: 3.Ctightness.tex
 \section{$C$-tightness}
 \label{sec:tightness}
 \setcounter{equation}{0}
 
 In this section we establish the $C$-tightness of the sequence of rescaled volatility processes and hence the existence of a weak continuous accumulation point.
 To simplify notation we express the rescaled processes as
 \beqlb
 	V^{(n)}(t) =I^{(n)}(t) + \frac{\zeta^\mathtt{m}_n}{\beta_n}\cdot J_1^{(n)}(t)+\frac{\zeta^\mathtt{l}_n}{ n^{\alpha-1}}\cdot J_2^{(n)}(t), \quad t \geq 0, 
 \eeqlb
 where the processes on the right-hand side of the above equation are defined by
 \begin{align}\label{decomV}
 \begin{split}
  I^{(n)}(t)
  &:=
  \frac{\mu_n}{n^{2\alpha-1}}+\frac{\mu_n}{n^{2\alpha-1}} \cdot \int_0^{nt}R_n(s)ds 
   +  \frac{\Lambda_n(nt)}{n^{2\alpha-1}} + \int_0^{nt} R_n(nt-s) \frac{\Lambda_n(s)}{n^{2\alpha-1}} ds,\\
 J_1^{(n)}(t)
 &:= 
 \int_0^t \int_0^{V^{(n)}(s-)} \frac{1}{n^{2\alpha-1}} \cdot R_n\big(n(t-s)\big) \widetilde{N}^{(n)}_\mathtt{m}(ds , dz),\\
 J_2^{(n)}(t)
 &:= 
 \int_0^t \int_0^\infty \int_0^{V^{(n)}(s-)} \frac{1}{n^{\alpha}} \cdot\mathcal{R}_n\big(n(t-s),ny\big)\widetilde{N}^{(n)}_\mathtt{l}(ds, dy, dz).
 \end{split}
 \end{align} 
 The boundedness and convergence of the two sequences  $\{\zeta^\mathtt{m}_n/\beta_n\}_{n\geq 1} $ and $\{\zeta^\mathtt{l}_n/n^{\alpha-1}\}_{n\geq 1}$ follow from Condition~\ref{main.Condition}(2). Moreover, the $C$-tightness of the sequence $\{I^{(n)}\}_{n\geq 1}$ has already been established in \cite{HorstXuZhang2023a}. 

 \begin{proposition}[Corollary 3.2 in \cite{HorstXuZhang2023a}] \label{Prop.ConvergenceRn}
 As $n\to\infty$, we have that
 \beqnn
 \sup_{t\geq 0} \Big|I^{(n)}(t) -V_*(0) \cdot \big(1-F^{\alpha,\gamma}(t)\big)-  \frac{a}{b}\cdot F^{\alpha,\gamma}(t) \Big| \to 0. 
 \eeqnn  
 \end{proposition}
  
 Establishing the $C$-tightness of the sequences  $\{J^{(n)}_1\}_{n\geq 1}$ and  $\{J^{(n)}_2\}_{n\geq 1}$ is more challenging. 
 We start with a series of a priori estimates that will allow us to prove the desired tightness combining the classical Kolmogorov tightness criterion and the tightness criterion for c\`adl\`ag processes established in \cite{HorstXuZhang2023a}.


\subsection{A priori estimates}

 We first recall the following two results that were established in \cite{HorstXuZhang2023a,Xu2024b}.
 They provide growth estimates for the resolvents $\{R_n\}_{n\geq 1}$ defined in \eqref{RH} as well as their derivative, denoted by $\{R'_n\}_{n\geq 1}$ and the two-parameter functions $\{\mathcal{R}_n\}_{n \geq 1}$ defined in \eqref{RRH}. 

\begin{proposition}[Proposition 3.6 and 3.7 in \cite{HorstXuZhang2023a}] \label{Prop.RH01}
	There exists a constant $C>0$ such that for any $n\geq 1$ and $t\geq 0$,
	\begin{align}\label{ubRH}
	R_n(t) \leq C(1+t)^{\alpha-1} 
	\quad\mbox{and}\quad 
	\big| R'_n(t) \big| \leq C\cdot (1+t)^{\alpha-2},\quad t\geq 0. 
	\end{align}
 
\end{proposition} 

 \begin{proposition}[Proposition 4.9 in \cite{Xu2024b}]\label{Prop.RH02}
For any $p>1+\alpha$, there exists a constant $C>0$ such that for any $t\geq 0$ and $n\geq 1$,
\begin{align}\label{ubR}
 \int_0^t \int_0^\infty \Big|n^{-\alpha}\cdot\mathcal{R}_n(ns,ny) \Big|^p\cdot n^{\alpha+1} \cdot \nu(n\cdot dy)\, ds \leq C (1+t)^{\alpha(p-1)}.
\end{align}
 \end{proposition}
 
The next proposition provides key moment estimates for stochastic integrals with respect to Poisson random measures that involve stochastic integral boundaries.

 \begin{proposition}[Theorem D.1 in \cite{Xu2024b}]
 \label{2pine}
 Let $p\geq 1$ and $T\geq 0$, and let $X:=\{X(t):t\geq 0\}$ be a non-negative c\`adl\`ag $(\mathscr{F}_t)$-progressive process such that
 \beqnn
 L:= \sup_{t\in[0,T]}\textbf{E}\Big[\big|X(t)\big|^p\Big] <\infty. 
 \eeqnn
 Furthermore, let $\widetilde{N}_0(ds,dy,dz)$ be a compensated $(\mathscr{F}_t)$-Poisson random measure on $\mathbb{R}_+^3$ with intensity 
 \beqnn
 \eta \cdot ds\,\nu_0(dy)\,dz
 \eeqnn
 for some constant $\eta>0$ and $\sigma$-finite measure $\nu_0(dy)$ on $\mathbb{R}_+$. 
 Then, for any interval $A\subset\mathbb{R}_+$ and any real-valued function $f$ on $[0,T]\times \mathbb R_+$ that satisfies 
 \beqnn
 \int_0^T \int_A \big|f(s,y)\big|^{2p} \, ds\,\nu_0(dy) <\infty,
 \eeqnn 
 there exists a constant $C>0$ that depends only on $p$ and $L$ such that 
 \beqnn
 \begin{split}	 \mathbf{E}\Big[\Big|\int_0^T\int_A\int_0^{X(s-)}f(s,y)\widetilde{N}_0(ds,dy,dz)\Big|^{2p}\Big]
 	&\leq C\Big(\eta\cdot\int_0^T\int_A |f(s,y)|^2\nu_0(dy)ds\Big)^p \\
 	&\quad +C\cdot \eta\cdot\int_0^T\int_A \big|f(s,y) \big|^{2p}\nu_0(dy)ds.
  \end{split}
 \eeqnn
\end{proposition}

 The above proposition will be repeatedly applied to the stochastic integrals in \eqref{decomV}. In this case the stochastic process $X$ is given by the rescaled volatility process $V^{(n)}$. 
 To apply this proposition in our setting, we thus need to verify that the processes $\{V^{(n)}\}_{n\geq 1}$ are uniformly $L^p$-bounded. This is shown by the following lemma. 

\begin{lemma} \label{V2p}
 For any $p\geq 0$, there exists a constant $C>0$ such that for any $i=1,2$, $T\geq 0$ and $n\geq 1$, 
 \beqnn
 \sup_{t\in[0,T]} \mathbf{E}\Big[\big|J_i^{(n)}(t)\big|^{2p}\Big]\leq C \cdot(1+T)^{2p\alpha}. 
 \eeqnn
 In particular,
 \beqnn
 \sup_{t\in[0,T]} \mathbf{E}\Big[\big|V^{(n)}(t)\big|^{2p}\Big] \leq C\cdot (1+T)^{2p\alpha}. 
 \eeqnn
 \end{lemma}

 \proof 
Let us fix $p \geq 0$ and $T>0$. In view of Proposition \ref{Prop.ConvergenceRn}, the uniform bound on the rescaled volatility processes follows from the bounds on the integral processes  $\{J_i^{(n)}\}_{n\geq 1}$ for $i=1,2$. 

To overcome the problem that these processes are not (local) martingales, due to the time-dependence of their integrands, for any $t\geq 0$ and $n\geq 1$ we introduce the following auxiliary processes on $[0,t]$: 
 \begin{align*}
 J^{(n)}_{1,t}(r) &:= \int_0^r \int_0^{V^{(n)}(s-)}  \frac{1}{ n^{2\alpha-1}}\cdot R_n\big(n(t-s)\big)  \widetilde{N}^{(n)}_{\mathtt m}(ds,dz), \quad r \in [0,t],\\
 J^{(n)}_{2,t}(r) &:= 
 \int_0^r\int_0^\infty \int_0^{V^{(n)}(s-)}   \frac{1}{ n^{\alpha}}\cdot \mathcal{R}_n \big(n(t-s),ny\big)\widetilde{N}^{(n)}_\mathtt{l}(ds, dy, dz), \quad r \in [0,t]. 
 \end{align*}

 The integrands of the processes $ J^{(n)}_{i,t}$ are uniformly bounded, due to Proposition \ref{Prop.RH01}, and $\widetilde N^{(n)}_\mathtt{m}(ds,dz)$ and $\widetilde N^{(n)}_\mathtt{l}(ds,dy,dz)$ are compensated Poisson random measures on $[0,t] \times \mathbb{R}_+^2 $.
 Hence, the auxiliary processes are martingales on $[0,t]$. By construction, they satisfy 
 \beqnn
	J^{(n)}_{i,t}(0) \overset{\rm a.s.}=0 \quad \mbox{and} \quad J^{(n)}_{i,t}(t)\overset{\rm a.s.}=J_i^{(n)} (t).
 \eeqnn

 In view of Jensen's inequality it suffices to establish the desired bounds for $2p=2^k$ with $k\in  \mathbb{N}$. 
 Since
 \beqlb  \label{expectation}
 \mathbf{E}\Big[ \big|J_i^{(n)} (t) \big|^{2^k}\Big]= \mathbf{E}\Big[\big|J^{(n)}_{i,t}(t)\big|^{2^k}\Big],\quad i=1,2 ,
 \eeqlb
 it suffices to consider the auxiliary processes for which we proceed by induction, starting with the case $k=0$. Taking expectations on both sides of \eqref{rescaledV} we get that 
 \beqlb \label{eqn.Upperbound.02}
 \mathbf{E}\Big[\big|V^{(n)}(t)\big|\Big]= \mathbf{E}\Big[ V^{(n)}(t) \Big] =  I^{(n)}(t)+ 
 \frac{\zeta^\mathtt{m}_n}{\beta_n}\cdot\mathbf{E}\Big[ J_{1,t}^{(n)}(t)\Big]
 +\frac{\zeta^\mathtt{l}_n}{ n^{\alpha-1}}\cdot \mathbf{E}\Big[J_{2,t}^{(n)}(t)\Big]=I^{(n)}(t),\quad t\geq 0.
 \eeqlb
 Hence Proposition \ref{Prop.ConvergenceRn} yields a constant $C>0$ that is independent of $T$ such that 
 \beqnn
 \sup_{n\geq 1}\sup_{t\in[0,T]} \mathbf{E}\Big[\big|V^{(n)}(t)\big| \Big]\leq C  . 
 \eeqnn

Using the same arguments as in the proof of \cite[Lemma 3.9]{HorstXuZhang2023a} we can also show that there exists a constant $C>0$ that does not depend on $T$ such that
\beqnn
  \sup_{n\geq 1}\sup_{t\in[0,T]}\mathbf{E}\Big[ \big|J_1^{(n)} (t) \big| \Big]\leq C. 
 \eeqnn 
 
 It remains to establish the moment estimates on $\{J^{(n)}_2\}_{n\geq 1}$, accordingly, $\{J^{(n)}_{2,t}\}_{n\geq 1}$. 
 Using the Burkholder-Davis-Gundy inequality and then applying Jensen's inequality, there exists a constant $C>0$ that is again independent of $n$ and $t$ such that 
 \beqnn
 \mathbf{E}\Big[ \big|J_2^{(n)} (t) \big| \Big] 
  \ar=\ar \mathbf{E}\Big[\big|J^{(n)}_{2,t}(t)\big| \Big]  \cr
 \ar\leq\ar C \cdot \mathbf{E}\Big[ \Big|\int_0^t \int_0^\infty\int_0^{V^{(n)}(s-)}  \frac{\big|\mathcal{R}_n\big(n(t-s),ny\big)\big|^2}{n^{2\alpha}}  N^{(n)}_\mathtt l(ds,dy,dz)\Big|^{1/2}\Big]\cr
 \ar\leq\ar C \cdot \Big(\mathbf{E}\Big[\int_0^t \int_0^\infty \int_0^{V^{(n)}(s-)} \frac{\big|\mathcal{R}_n\big(n(t-s),ny\big)\big|^2}{n^{2\alpha}}  N^{(n)}_\mathtt l(ds,dy,dz)\Big]\Big)^{1/2} \cr
 \ar\leq\ar C \cdot \Big(\int_0^t \int_0^\infty\mathbf{E}\big[V^{(n)}(s)\big]\cdot \lambda_n^\mathtt l\cdot n^{2\alpha}\cdot  \frac{\big|\mathcal{R}_n\big(n(t-s),ny\big)\big|^2}{n^{2\alpha}}  \nu(n\cdot dy)ds\Big)^{1/2}. 
 \eeqnn 
 In view of \eqref{eqn.Upperbound.02} and using that $\lambda^{(n)}_\mathtt l\sim n^{1-\alpha}$ and the inequality \eqref{ubR} we can further conclude that 
 \beqnn
 \sup_{t\in[0,T]}\mathbf{E}\Big[ \big|J_2^{(n)} (t) \big| \Big] 
 \ar\leq \ar  C \Big(\int_0^T \int_0^\infty n^{\alpha+1}\cdot
 \Big|\frac{\mathcal{R}_n(ns,ny)}{n^{\alpha}}\Big|^2 \nu(n\cdot dy)ds\Big)^{1/2} 
 \leq  C \cdot (1+T)^{\alpha/2}  
 \eeqnn
 for some constant $C>0$ that is independent of $n$ and $T$. Putting the preceding estimates together proves the desired result for $k=0$. 
 
 We now proceed under the assumption that the inequality holds for $2p=2^k$ for some $k \geq 0$ and prove that it holds for $2p=2^{k+1}$.  By Proposition \ref{Prop.ConvergenceRn} and Condition~\ref{main.Condition}(2), there exists a constant $C>0$ such that for all $t\in[0,T]$ and $ n\geq 1$, 
 \beqnn
 \mathbf{E}\Big[\big|V^{(n)}(t)\big|^{2p}\Big]\leq C\cdot \Big(1
+\mathbf{E}\Big[\big|J_1^{(n)}(t)\big|^{2p}\Big]
 +
\mathbf{E}\Big[\big|J_2^{(n)}(t)\big|^{2p}\Big]
 \Big).
 \eeqnn
 Following the arguments given in the proof of  \cite[Lemma 3.9]{HorstXuZhang2023a} one can show that uniformly in $T\geq 0$,
 \beqnn
 \sup_{n\geq 1}\sup_{t\in[0,T]} \mathbf{E}\Big[\big|J_1^{(n)}(t)\big|^{2p}\Big]\leq C\cdot (1+T)^{2p\alpha} .
 \eeqnn 
 To obtain a similar estimate for the processes $\{J_2^{(n)}\}_{n\geq 1}$ we first notice that the induction hypothesis for $2p=2^k$ allows us to apply Proposition \ref{2pine}. Hence, we have that uniformly in $T\geq 0$, 
 \begin{align*}
 \sup_{t\in[0,T]}\mathbf{E}\Big[\big|J_2^{(n)}(t)\big|^{2p}\Big] 
 =\sup_{t\in[0,T]}\mathbf{E}\Big[\big|J_{2,t}^{(n)}(t)\big|^{2p}\Big] 
 \leq & ~  C \cdot\Big| \int_0^T \int_0^\infty \Big| \frac{\mathcal{R}_n(ns,ny)}{n^{\alpha}} \Big|^2 n^{\alpha+1} \nu(n\cdot dy) ds\Big|^p \\
    &+ C \cdot\int_0^T \int_0^\infty \Big|\frac{\mathcal{R}_n(ns,ny)}{n^{\alpha}} \Big|^{2p} n^{\alpha+1} \nu(n\cdot dy) ds.  
 \end{align*}
 Applying Proposition \ref{2pine} again allows us to conclude that for some constant $C>0$ independent of $T$, 
 \beqnn
 \sup_{n\geq 1}\sup_{t\in[0,T]}\mathbf{E}\Big[\big|J_2^{(n)}(t)\big|^{2p}\Big] \leq C\cdot(1+T)^{\alpha p}+C\cdot(1+T)^{\alpha(2p-1)}  
 \leq C \cdot (1+T)^{\alpha(2p-1)} 
 \eeqnn    
 and hence that
 \beqnn
  \sup_{n\geq 1}\sup_{t\in[0,T]} 
  \Big(\mathbf{E}\Big[\big|J_1^{(n)}(t)\big|^{2p}\Big] +
  \mathbf{E}\Big[\big|J_2^{(n)}(t)\big|^{2p}\Big] +  \mathbf{E}\Big[\big|V^{(n)}(t)\big|^{2p} \Big]\Big)\leq
 C\cdot (1+T)^{2p\alpha} . 
 \eeqnn
  \qed

 Armed with the preceding moment estimates, the $C$-tightness of the sequence $\{J_1^{(n)}(t):t\geq 0\}_{n\geq 1}$ can be proved in the same way as in \cite{HorstXuZhang2023a}. The detailed proof is omitted. 
 
 \begin{lemma}
 The sequence $\{J_1^{(n)} \}_{n\geq 1}$ is $C$-tight.  
 \end{lemma}

\subsection{$C$-tightness of $\{J_2^{(n)}\}_{n\geq 1}$}
\label{ctightVn}
 
 It remains to show that the sequence $\{J_2^{(n)}\}_{n\geq 1}$ that captures the impact of limit orders on the market dynamics is $C$-tight. In terms of the two-parameter function
 \beqnn
 \mathcal{R}_n(nt,ny)= \mathbf{1}_{\{y>t\}}+n\cdot \int_0^tR_n\big(n(t-s)\big)\mathbf{1}_{\{y>s\}}ds,\quad t\geq 0,\ y>0, 
 \eeqnn
 introduced in \eqref{RRH}, the process $\{J_2^{(n)}(t):t\geq 0\}$ can be decomposed as
 \beqnn
 J^{(n)}_2(t) = J_{2,c}^{(n)}(t)  + \varepsilon^{(n)}_1 (t),
 \eeqnn
 where the {\sl continuous} processes $ J_{2,c}^{(n)} $ and the ``error term'' $ \varepsilon^{(n)}_1 $ are given by, respectively,
 \beqnn
 J_{2,c}^{(n)}(t)\ar :=\ar    \int_0^t\int_0^{\infty}\int_0^{V^{(n)}(s-)}
 \Big(\int_{(t-s-y)^+}^{t-s} \frac{R_n(nr)}{n^{\alpha-1}}  dr\Big)
 \widetilde{N}_{\mathtt{l}}^{(n)}(ds,dy, dz), \cr
 \varepsilon^{(n)}_1 (t) \ar:=\ar  
 \int_0^t\int_0^{\infty}\int_0^{V^{(n)}(s-)}  \frac{\mathbf{1}_{\{y>t-s\}}}{n^\alpha}\,
 \widetilde{N}_{\mathtt{l}}^{(n)}( ds, dy, dz) .
 \eeqnn

 The $C$-tightness of $\{J_{2}^{(n)}\}_{n\geq 1}$ is obtained by separately proving that the sequences $\{ J_{2,c}^{(n)} \}_{n\geq 1}$ and $ \{\varepsilon^{(n)}_1\}_{n\geq 1}$ are $C$-tight. Moreover, we will also see that  $\{\varepsilon^{(n)}_1\}_{n\geq 1}$ converges weakly to zero.

\begin{remark}
 The continuity of the process $J_{2,c}^{(n)}$ is due to the inner integral that depends on the time variable and is  continuous on the whole time interval. 
 Indeed, since $ \widetilde{N}_{\mathtt{l}}^{(n)}(ds,dy, dz)=  N_{\mathtt{l}}^{(n)}(ds,dy, dz) -  n^{2\alpha}\cdot \lambda^{\mathtt{l}}_n\cdot \alpha \cdot ds\,\nu(n\cdot dy)\,dz$ we have that 
 \beqlb \label{eqn.301}
 \begin{split}
  J_{2,c}^{(n)}(t) = &   \int_0^t\int_0^{\infty}\int_0^{V^{(n)}(s-)}\Big(\int_{(t-s-y)^+}^{t-s}\frac{R_n(nr)}{n^{\alpha-1}} dr\Big)
  N_{\mathtt{l}}^{(n)}(ds,dy, dz)\cr
  & -   \int_0^tV^{(n)}(s-)  ds \int_0^{\infty} \Big(\int_{(t-s-y)^+}^{t-s} \frac{R_n(nr)}{n^{\alpha-1}} dr\Big)
  n^{2\alpha}\cdot \lambda^{\mathtt{l}}_n\cdot \alpha \cdot \nu(n\cdot dy)dz.
  \end{split}
 \eeqlb
 The continuity of the second term on the right-hand side of (\ref{eqn.301}) follows from the elementary properties of convolution.  
 Since 
 \beqnn
 \int_0^t ds \int_0^\infty V^{(n)}(s-) \cdot n^{2\alpha}\cdot \lambda^{\mathtt{l}}_n\cdot \alpha \cdot  ds \nu(n\cdot dy) =n^{2\alpha}\cdot \lambda^{\mathtt{l}}_n\cdot \alpha \cdot \int_0^t V^{(n)}(s-) ds<\infty ,\quad a.s., \quad 
 \eeqnn
 for any $t\geq 0$, 
 the point process  $\{(\tau_k, y_k,z_k)\}_{k\in\mathbb{N}_+}$ corresponding to $N^{(n)}_\mathtt{l}$ has only finitely many jumps in any bounded time interval and hence the first term on the right-hand side of (\ref{eqn.301}) can be written as a finite sum with the $k$-th summand  being of the form
 \beqnn
 \int_{(t-\tau_k -y_k)^+}^{t-\tau_k} \frac{R_n(nr)}{n^{\alpha-1}} dr \in \mathbf{C}(\mathbb{R};\mathbb{R}_+).
 \eeqnn
This yields the continuity of $J_{2,c}^{(n)}$.  
\end{remark}


 \subsubsection{$C$-tightness of $\{J_{2,c}^{(n)}\}_{n\geq 1}$}

 In this section, we prove the $C$-tightness of $\{J_{2,c}^{(n)}\}_{n\geq 1}$ and first 
 recall the forward difference operator $\Delta_h$ and backward difference operator $\nabla_h$ introduced at the end of the introduction.
  To simplify notation, we introduce the rescaled resolvent and its integral
 \beqlb\label{Rn}
 R^{(n)} (t):= \frac{R_n(nr)}{n^{\alpha-1}}
 \quad \mbox{and}\quad 
 \mathcal{I}_{R^{(n)}}(t):= \int_0^t R^{(n)}(s)ds ,\quad n\geq 1,\  t\geq 0, 
 \eeqlb
 which allow us to write the process $J_{2,c}^{(n)} $ as 
 \beqlb\label{I2}
 J_{2,c}^{(n)}(t)\ar =\ar    \int_0^t\int_0^{\infty}\int_0^{V^{(n)}(s-)}
 \nabla_y \mathcal{I}_{R^{(n)}}(t-s)
 \widetilde{N}_{\mathtt{l}}^{(n)}(ds,dy, dz),\quad t\geq 0. 
 \eeqlb  
 The sequence of continuous processes $\{ J_{2,c}^{(n)} \}_{n\geq 1}$ is $C$-tight if it meets the Kolmogorov-Chentsov condition, that is, if for each $T\geq 0$, there exist constants $C ,p, \kappa>0 $   such that for any $h\in(0,1)$, 
 \beqlb\label{eqn.KolTight} 
 \sup_{n\geq 1} \mathbf{E}\Big[ \big|\Delta_h J_{2,c}^{(n)} (t)\big|^{p}\Big]\leq C\cdot  h^{1+\kappa} . 
 \eeqlb
 The proof uses the following uniform estimates on the functions $\{\mathcal{I}_{R^{(n)}}\}_{n\geq 1}$ and their increments.  
 A direct consequence of \eqref{ubRH} tells that for some constant $C>0$ that is independent of $n$ and $t$,
 \beqlb\label{UpperBIRn}
 \sup_{n\geq 1} \mathcal{I}_{R^{(n)}} (t) \leq  C \cdot t^\alpha ,\quad t\geq 0.
 \eeqlb

 \begin{proposition}\label{Prop:Appendix1}
 There exists a constant $C>0$ such that for any $n\geq 1$ and $t,y,h\geq 0$,  
 \beqlb 
  \big|\nabla_y \mathcal{I}_{R^{(n)}}(t)\big|
  = \big|\Delta_y \mathcal{I}_{R^{(n)}}(t-y) \big|
  \ar\leq\ar C \cdot \Big ( \Big(   \big|(t-y)^+\big|^{\alpha-1}\cdot y\Big)\wedge t^\alpha\wedge y^\alpha \Big) ,\label{1stDiff} \\
  \big|\Delta_h \nabla_y\mathcal{I}_{R^{(n)}}(t)\big| 
  =  \big| \nabla_y\Delta_h\mathcal{I}_{R^{(n)}}(t)\big|
  \ar\leq\ar C\cdot \Big(  h^\alpha \wedge \Big( h\cdot\big|(t-y)^+\big|^{\alpha-1} \Big)\wedge \Big( h\cdot \big|(t-y)^+\big|^{\alpha-2} \cdot y\Big) \Big) . \qquad \label{2rdDiff}
 \eeqlb 
 \end{proposition}
 \proof
 The inequality (\ref{1stDiff}) follows from \eqref{ubRH} and the inequality $|x^\alpha-z^\alpha|\leq x^\alpha\wedge |x-z|^\alpha$ for $x\geq z\geq 0$, i.e., 
 \beqnn
 \big|\nabla_y \mathcal{I}_{R^{(n)}}(t)\big| 
 = \int_{t-y}^t \frac{R_n(nr)}{n^{\alpha-1}}dr 
 \ar\leq\ar C  \cdot \int_{(t-y)^+}^t r^{\alpha-1}dr \cr
 \ar\leq\ar  C \cdot \Big ( \Big(   \big|(t-y)^+\big|^{\alpha-1} \cdot y\Big)\wedge \big|t^\alpha-(t-y)^\alpha\big|  \Big)\cr
 \ar\leq\ar C \cdot \Big ( \Big(   \big|(t-y)^+\big|^{\alpha-1}\cdot y\Big)\wedge t^\alpha\wedge y^\alpha \Big) ,
 \eeqnn
 for some constant $C>0$ that is independent of $n$ and $t$. 
 To prove the inequality (\ref{2rdDiff}), we first use (\ref{1stDiff}) with $y=h$ to get
 \beqnn
 \big|\Delta_h \nabla_y \mathcal{I}_{R^{(n)}}(t)\big|
 = \left |
 \Delta_h \mathcal{I}_{R^{(n)}}(t)-\Delta_h \mathcal{I}_{R^{(n)}}(t-y)\right | 
 \leq 
 \left |\Delta_h \mathcal{I}_{R^{(n)}}(t)\right|+
 \left |\Delta_h \mathcal{I}_{R^{(n)}}(t-y)\right| \leq 
 C\cdot h^\alpha,
 \eeqnn
 uniformly in $n\geq 1$ and $t,y,h\geq 0$. 
 On the other hand, by using the second inequality in (\ref{ubRH}) and then the mean value theorem we see that $(R^{(n)})'(t)\leq C\cdot t^{\alpha-2}$ uniformly in $t\geq 0$ and  
 \beqnn
 \big|\Delta_h \nabla_y \mathcal{I}_{R^{(n)}}(t)\big|
 \ar\leq\ar C\cdot \int_0^y \int_0^h \big| R^{(n)'}(t+\tilde{h}-\tilde{y}) \big|\, d\tilde{h}\,d\tilde{y}\cr
 \ar\leq\ar 
 C\int_0^y \int_0^h\big|(t+\tilde{h}-\tilde{y})^+\big|^{\alpha-2}\,d\tilde{h}\,d\tilde{y} \cr
 \ar\leq\ar C\cdot  \Big( h\cdot\big|(t-y)^+\big|^{\alpha-1} \Big)\wedge \Big( h\cdot \big|(t-y)^+\big|^{\alpha-2} \cdot y\Big). 
 \eeqnn 
 \qed

 \begin{lemma} \label{ctightness of I2}
 For any $p\geq 3$ and $T\geq 0$, there exists constants $C,\kappa>0$ such that for any  $h\in (0,1)$,  
 \beqlb \label{Iltight}
   \sup_{n\geq 1}\sup_{t\in[0,T]}  \textbf{E}\Big[\big|\Delta_h J_{2,c}^{(n)}(t)\big|^{2p}\Big]\leq C\cdot h^{1+\kappa}.
 \eeqlb
 Moreover, the sequence $\{J_{2,c}^{(n)}\}_{n\geq 1}$ is $C$-tight.
 \end{lemma} 
 \proof 
 For any $t,h\geq 0$, by \eqref{I2} we have   $\Delta_h J_{2,c}^{(n)}(t)=  J^{(n)}_{2,1}(t,h)+  J^{(n)}_{2,2}(t,h)$ with 
 \begin{align*}
  J^{(n)}_{2,1}(t,h)&:=   \int_t^{t+h}\int_0^{\infty}\int_0^{V^{(n)}(s-)} \nabla_y \mathcal{I}_{R^{(n)}_H}(t+h-s) \widetilde{N}_\mathtt{l}^{(n)}(ds, dy, dz), \\  
  J^{(n)}_{2,2}(t,h)
  &:= \int_0^t\int_0^\infty\int_0^{V^{(n)}(s-)} \Delta_h\nabla_y \mathcal{I}_{R^{(n)}_H}(t-s)  \widetilde{N}_\mathtt{l}^{(n)}( ds, dy,  dz).
 \end{align*}
 Depending on the value of the integration variable $y$, we can further decompose the preceding functions as 
 \beqnn
 J^{(n)}_{2,1}(t,h)= I^{(n)}_{1}(t,h)+ I^{(n)}_{2}(t,h)
 \quad \mbox{and}\quad 
 J^{(n)}_{2,2}(t,h)= I^{(n)}_{3}(t,h)+ I^{(n)}_{4}(t,h)+I^{(n)}_{5}(t,h),
 \eeqnn
 where 
 \beqnn
 I^{(n)}_{1}(t,h)
 \ar:=\ar\int_t^{t+h}\int_0^{t+h-s}\int_0^{V^{(n)}(s-)}
 \nabla_y \mathcal{I}_{R^{(n)}}(t+h-s)\widetilde{N}_{\mathtt{l}}^{(n)}(ds, dy, dz),\cr
 I^{(n)}_{2}(t,h)
 \ar:=\ar \int_t^{t+h}\int_{t+h-s}^{\infty}\int_0^{V^{(n)}(s-)}\mathcal{I}_{R^{(n)}}(t+h-s)\widetilde{N}_{\mathtt{l}}^{(n)}(  ds, dy, dz), \cr
 I^{(n)}_{3}(t,h)
 \ar:=\ar \int_0^t\int_0^{t-s}\int_0^{V^{(n)}(s-)}\Delta_h\nabla_y\mathcal{I}_{R^{(n)}}(t-s)\widetilde{N}_{\mathtt{l}}^{(n)}(ds,dy, dz) ,\cr
 I^{(n)}_{4}(t,h)
 \ar:=\ar \int_0^t\int_{t-s}^{t+h-s}\int_0^{V^{(n)}(s-)}\Big(\nabla_y\mathcal{I}_{R^{(n)}}(t+h-s)-\mathcal{I}_{R^{(n)}}(t-s)\Big)\widetilde{N}_{\mathtt{l}}^{(n)}(ds,dy,dz), \cr
 I^{(n)}_{5}(t,h)
 \ar:=\ar   \int_0^t\int_{t+h-s}^{\infty}\int_0^{V^{(n)}(s-)}\Delta_h \mathcal{I}_{R^{(n)}}(t-s)\widetilde{N}_{\mathtt{l}}^{(n)}(ds,dy,dz). 
 \eeqnn 
 It suffices to verify the desired moment condition separately for the processes $\tilde I^{(n)}_i$ with $i=1, ..., 5$, i.e.~to show that for each $p>0$ there exist constants $C >0$ such that
 \beqnn
 \sup_{n\geq 1}\sup_{t\in[0,T]}  \mathbf{E}\Big[\big|I^{(n)}_i (t,h)\big|^{2p}\Big]\leq C\cdot h^{\alpha(p-1)},\quad i=1,\cdots,5.  
 \eeqnn
 Since $p\geq 3> \frac{1}{\alpha}+1$, we have $\alpha(p-1)>1$ and then the desired inequality (\ref{Iltight}) holds. 

 In what follows we establish the estimate for the processes $\tilde I^{(n)}_1$.  
 An application of Proposition \ref{2pine} along with the change of variables $s=t+h-s$ shows that
 \beqlb\label{eqn301}
 \mathbf{E}\Big[\big| I^{(n)}_1(t)\big|^{2p}\Big]\leq 
 C \cdot \bigg|\int_0^h \int_0^s \frac{\big|\nabla_y\mathcal{I}_{R^{(n)}}(s)\big|^2}{y^{\alpha+2}} dy ds \bigg|^p+C\int_0^h  \int_0^s \frac{\big|\nabla_y\mathcal{I}_{R^{(n)}}(s)\big|^{2p}}{y^{\alpha+2}} dy ds,
 \eeqlb
 uniformly in $t\in[0,T]$, $h\in(0,1)$ and $n\geq 1$. 
 To further bound the above integrals, for $q\geq 1$ we apply (\ref{1stDiff}) to get 
 \beqnn
 \int_0^h \int_0^s \frac{\big|\nabla_y\mathcal{I}_{R^{(n)}}(s)\big|^q}{y^{\alpha+2}} dy ds
 \ar\leq\ar C\cdot \int_0^h \int_0^s \frac{\big|s^\alpha \wedge \big( y\cdot |(s-y)^+|^{\alpha-1} \big) \big|^q}{y^{\alpha+2}} dy ds\cr
 \ar\leq\ar C\int_0^h \int_0^{s/2}\frac{(s-y)^{q(\alpha-1)}}{y^{\alpha+2-q}}dy ds +C\int_0^h \int_{s/2}^s\frac{s^{q\alpha}}{y^{\alpha+2}}dy ds\cr
 \ar\leq\ar  C\int_0^h s^{q(\alpha-1)} \int_0^{s/2}y^{q-\alpha-2}dy ds+
 C\int_0^h s^{q\alpha}\cdot s^{-\alpha-1}ds\\
 \ar\leq\ar C \int_0^h s^{\alpha q-q+q-\alpha -1}ds+ C h^{\alpha(q-1)}\\
 \ar\leq\ar C\cdot h^{\alpha(q-1)}.
 \eeqnn
 Plugging this with $q\in\{2,2p\}$ back into the right side of (\ref{eqn301}), we have uniformly in $h\in[0,T]$ that
 \beqnn
 \sup_{n\geq 1}\sup_{t\in[0,T]}\mathbf{E}\Big[\big|  I^{(n)}_1(t,h)\big|^{2p}\Big]
 \leq C\cdot \big( h^{p\alpha} + h^{\alpha(2p-1)} \big)
 \leq C\cdot h^{\alpha(p-1)}. 
 \eeqnn
 Similar estimates can be established for the remaining processes using similar arguments. 
 \qed


 \subsubsection{$C$-tightness of $\{\varepsilon_1^{(n)} \}_{n\geq 1}$}

 In this section we prove the $C$-tightness of the sequence $\{\varepsilon_1^{(n)} \}_{n\geq 1}$.
 The key to our tightness proof is the following $C$-tightness criterion for c\`adl\`ag stochastic processes introduced in \cite{HorstXuZhang2023a}. 

 \begin{lemma}[Lemma 3.5 in \cite{HorstXuZhang2023a}] \label{tightness condition}
 Let $\{X^{(n)} \}_{n\geq 1}$ be a sequence of c\`adl\`ag stochastic processes defined on a common probability space such that $$\sup_{n \geq 1}  \mathbf{E} \Big[\big|X^{(n)}(0)\big|^q\Big] < \infty$$ for some $q>0$. 
 Then, the sequence $\{X^{(n)} \}_{n\geq 1}$ is $C$-tight if for any $T\geq 0$ and some constant $\theta>2$, the following two  conditions hold.
 \begin{itemize}
 	\item[(1)] There exist constants $C >0 $, $p\geq 1$, $m\in\{1,2, ...\}$ and pairs $\{ (a_i,b_i) \}_{i=1,\cdots ,m}$ satisfying 
 	\beqnn
 	a_i\geq 0,\quad b_i>0,\quad   b_i + \frac{a_i}{\theta}  >1 , 
 	\eeqnn
    such that for all $n\geq 1$ and $h \in (0,1)$,
 	\beqlb \label{condition2}
 	 \sup_{t\in [0,T]}	\mathbf{E}\Big[\big|\Delta_h X^{(n)}(t)\big|^{p}\Big]\leq 
 		C\cdot \sum_{i=1}^m \frac{h^{b_i}}{n^{a_i}}. 
 	\eeqlb 
   \item[(2)] $\displaystyle\sup_{k=0,1,\cdots,[Tn^\theta]} \sup_{h\in[0,1/n^{\theta}]} \big|\Delta_h X^{(n)}(k/n^\theta) \big|  \overset{\rm p}\to 0$  as $n\to\infty$.
\end{itemize} 
 \end{lemma}


In \cite{HorstXuZhang2023a} we proved that any sequence of c\`adl\`ag  stochastic processes that satisfies the conditions of the above lemma can be approximated in probability by a sequence of piecewise linear interpolations on the time grids $\left\{ \frac{k}{n^\theta} : k=0,1,2, ...\right\}$ that satisfy the Kolmogorov tightness criterion for continuous processes, due to Condition (1). The second condition guarantees that the original sequence and the sequence of linear interpolations converge to the same limit.   

We apply the above lemma to the sequence $\{\varepsilon_1^{(n)}\}_{n\geq 1}$. Since $\varepsilon^{(n)}_1(0)\overset{\rm a.s.}=0$ for all $n\geq 1$, it suffices to check Conditions (1) and (2) to establish the desired tightness.
\begin{proposition}
 The sequence $\{\varepsilon_1^{(n)} \}_{n\geq 1} $ is $C$-tight. 
\end{proposition}
\proof
 To verify that the sequence $\{\varepsilon_1^{(n)}\}_{n\geq 1}$ satisfies Condition (1) and (2) in Lemma \ref{tightness condition}, we first decompose the increment process   $\Delta_h\varepsilon^{(n)}_1(t)=\varepsilon^{(n)}_1(t+h)-\varepsilon^{(n)}_1(t) $ into the following two terms
 \beqnn
 \varepsilon^{(n)}_{11}(t,h)\ar:=\ar
   \int_t^{t+h}\int_{t+h-s}^{\infty}\int_0^{V^{(n)}(s-)}\frac{1}{n^\alpha} \widetilde{N}_{\mathtt{l}}^{(n)}(ds,dy,dz),\cr
 \varepsilon^{(n)}_{12}(t,h) \ar:=\ar - \int_0^{t}\int_{t-s}^{t+h-s}\int_0^{V^{(n)}(s-)}\frac{1}{n^{\alpha}}\widetilde{N}_{\mathtt{l}}^{(n)}(ds,dy,dz). 
 \eeqnn

\medskip

{\bf Step 1. The sequence $\{\varepsilon_1^{(n)}\}_{n\geq 1}$ satisfies 
Condition (1).} We consider the sequences $\{\varepsilon_{11}^{(n)}\}_{n\geq 1}$ and $\{\varepsilon_{12}^{(n)}\}_{n\geq 1}$ separately, starting with the first one.
\begin{itemize}
\item For $\varepsilon^{(n)}_{11}$, by Proposition \ref{2pine} and a change of variables we have 
\beqnn
  \sup_{t\in[0,T]} \mathbf{E}\Big[|\varepsilon^{(n)}_{11}(t,h)|^{2p}\Big]
  \leq  
   C\cdot \bigg|  \frac{1}{n^{\alpha-1}} \int_0^{h} \bar{\nu} (ns)  
  ds\bigg|^{p}+
   \frac{C}{n^{2p\alpha-(\alpha+1)}}\cdot \int_0^{h} \bar{\nu}(ns)ds,
\eeqnn
 where $\bar{\nu}(x):=\nu(x,\infty)= (1+x)^{-\alpha-1}$. Clearly, 
 \beqnn
 \int_0^{h} \bar{\nu} (ns) ds  \leq h.
 \eeqnn
 Moreover, using the inequality $|x^\alpha-y^\alpha|\leq |x-y|^\alpha$ for $x,y\geq 0$,  
 \beqnn
 \int_0^{h} \bar{\nu} (ns) ds  
 \ar=\ar \frac{1}{n\alpha}\Big(1 -\Big(\frac{1}{1+nh}\Big)^\alpha \Big)
 \leq  \frac C n\cdot \Big(1-\frac{1}{1+nh}\Big)^\alpha 
 \leq   C\cdot \frac{h^{\alpha}}{n^{1-\alpha}} .
 \eeqnn
 Therefore, for some constant $C>0$ that is independent of $h\in(0,1)$,
 \beqnn
 \sup_{t\in[0,T]}\mathbf{E} \Big[\big|\varepsilon^{(n)}_{11}(t,h)\big|^{2p}\Big]\leq C\cdot \Big( h^{p\alpha} + \frac{h}{n^{2p\alpha-(\alpha+1)}}  \Big).
 \eeqnn

 \item For $\varepsilon^{(n)}_{12}$, again by Proposition \ref{2pine} and a change of variables,
 \beqnn
 \sup_{t\in[0,T]}\mathbf{E} \Big[|\varepsilon^{(n)}_{12}(t,h)|^{2p}\Big]
 \leq 
 C\cdot \bigg| \frac{1}{n^{\alpha-1}}  \int_0^t\int_{ns}^{n(s+h)} \nu(dr) 
  ds\bigg|^{p}+ \frac{C}{n^{2p\alpha-(\alpha+1)}}\int_0^t\int_{ns}^{n(s+h)} \nu(dr) 
  ds.
 \eeqnn
 For the double integrals we compute: 
 \beqnn
  \int_0^t\int_{ns}^{n(s+h)} \nu(dr) 
 ds
 \ar=\ar   \int_0^t \int_{ns}^{n(s+h)}(\alpha+1)(1+r)^{-\alpha-2}drds\cr 
 \ar\leq\ar nh\cdot \int_0^t  (\alpha+1)(1+ns)^{-\alpha-2} ds\cr
 \ar\leq\ar  h \int_0^\infty  (\alpha+1)(1+s)^{-\alpha-2} ds \cr
 \ar=\ar h.
 \eeqnn
 Moreover, using the inequality $|x^\alpha-y^\alpha|\leq |x-y|^\alpha$ for $x,y\geq 0$ again,  
 \beqnn
  \int_0^t\int_{ns}^{n(s+h)} \nu(dr) 
  ds
 \ar=\ar 
 \int_0^t
 \Big[\big(1+ns\big)^{-\alpha-1}
 -\big(1+n(s+h)\big)^{-\alpha-1} \Big] ds\\
 \ar=\ar \frac{1}{\alpha n}
 \big(1-(1+nt)^{-\alpha}-(1+nh)^{-\alpha}+
 (1+nt+nh)^{-\alpha} \big)\\
 \ar\leq \ar  \frac{1}{\alpha n}\Big( 1-\Big(\frac{1}{1+nh}\Big)^{\alpha}\Big)\\
 \ar\leq\ar\frac{1}{\alpha n} \cdot \Big(\frac{n h }{1+nh }\Big)^\alpha\\
 \ar\leq\ar C\cdot \frac{h^{\alpha}}{n^{1-\alpha}}.
 \eeqnn
 As a result, there exists a constant $C>0$ that is independent of $h\in(0,1)$ such that
 \beqnn
 \sup_{t\in[0,T]}\mathbf{E} \Big[ \big|\varepsilon^{(n)}_{12}(t,h)\big|^{2p}\Big]\leq  C\cdot \Big( h^{p\alpha} + \frac{h}{n^{2p\alpha-(\alpha+1)}}  \Big) . 
 \eeqnn
\end{itemize}
Altogether, we arrive at the moment estimate
\beqnn
  \sup_{t\in[0,T]}  \mathbf{E}\Big[\big|\Delta_h\varepsilon^{(n)}_{\mathtt{l}}(t)\big|^{2p}\Big]\leq  C\cdot \Big( h^{p\alpha} + \frac{h}{n^{2p\alpha-(\alpha+1)}}  \Big), 
\eeqnn
 which is of the form \eqref{condition2} with
 \beqnn
  p\geq 2 ,\quad	m=2,\quad a_1=0,\quad b_1=p\alpha ,\quad a_2= 2p\alpha-(\alpha+1),\quad b_2= 1, \quad b_2+a_2/\theta > 1.
 \eeqnn

\medskip
\noindent {\bf Step 2. The sequence $\{\varepsilon^{(n)}_{1}\}_{n\geq 1}$ satisfies Condition (2).}  It holds that
 \beqnn
 \sup\limits_{k=0,1,\cdots,[Tn^\theta]}\sup\limits_{h\in[0,1/n^\theta]}
 \Big|\Delta_h \varepsilon^{(n)}_1(k/n^\theta) \Big| \leq
 \sup\limits_{k=0,1,\cdots,[Tn^\theta]}\sup\limits_{h\in[0,1/n^{\theta}]}
 \Big(\big|{\varepsilon}_{11}^{(n)}(k/n^\theta,h)\big|+\big|{\varepsilon}_{12}^{(n)}(k/n^\theta,h)\big|\Big).
 \eeqnn
  Let us first consider the sequence $\{\varepsilon_{11}^{(n)}\}_{n \geq 1}$.
 There exists a constant $C>0$ independent of $k$, $n$ and $h$ such that 
 \beqnn
  \sup_{h\in[0,1/n^{\theta}]}
 \Big|\varepsilon_{11}^{(n)}(k/n^\theta,h)\Big|
 \ar \leq\ar C\cdot  \sup_{h\in[0,1/n^{\theta}]}
  \int_{k/n^\theta}^{k/n^\theta+h}\int_{k/n^\theta +h-s}^\infty\int_0^{V^{(n)}(s-)}\frac{1}{n^\alpha} N_{\mathtt{l}}^{(n)}( ds, dy, dz) \cr
  \ar \ar + C\cdot  \sup_{h\in[0,1/n^{\theta}]} \int_{k/n^\theta}^{k/n^\theta+h}\int_{k/n^\theta +h-s}^\infty
  n\cdot V^{(n)}(s) \nu(n\cdot dy)ds \cr
  \ar\leq\ar   
  C  \int_{(k+1)/n^\theta}^{k/n^\theta+h}\int_0^\infty\int_0^{V^{(n)}(s-)}\frac{1}{n^\alpha} N_{\mathtt{l}}^{(n)}( ds, dy, dz)  
  +  C  \int_{k/n^\theta}^{(k+1)/n^\theta} 
  n\cdot V^{(n)}(s)  ds
 \eeqnn
 By the definition of the compensated Poisson random measure, the stochastic integral on the right-hand side of the second inequality can be  bounded by
 \beqnn
  \bigg|\int_{k/n^\theta}^{(k+1)/n^\theta}\int_0^\infty\int_0^{V^{(n)}(s-)}\frac{1}{n^\alpha} \widetilde{N}_{\mathtt{l}}^{(n)}( ds,  dy, dz)\bigg| +
 C\cdot \int_{k/n^\theta}^{(k+1)/n^\theta}  n\cdot V^{(n)}(s)ds. 
 \eeqnn
 Therefore, we have uniformly in $n\geq 1$ that
 \beqnn
 \sup_{k=0,1,\cdots,[Tn^\theta]}\sup_{h\in[0,1/n^{\theta}]}
 \Big|\varepsilon_{11}^{(n)}(k/n^\theta,h)\Big|
 \ar\leq\ar C\cdot \sup_{k=0,1,\cdots,[Tn^\theta]}   \bigg|\int_{k/n^\theta}^{(k+1)/n^\theta}\int_0^\infty\int_0^{V^{(n)}(s-)}\frac{1}{n^\alpha} \widetilde{N}_{\mathtt{l}}^{(n)}( ds,  dy,  \cdot dz)\bigg|\cr 
 \ar\ar  + 2 C\cdot \sup_{k=0,1,\cdots,[Tn^\theta]} \int_{k/n^\theta}^{(k+1)/n^\theta} 
 n\cdot V^{(n)}(s)  ds. 
 \eeqnn 
It remains to prove the convergence in probability to zero of the above terms. 

\begin{itemize}
    \item For any $\eta>0$, by Chebyshev's inequality,  
    \beqnn
    \lefteqn{\mathbf{P}\bigg( 
    	\sup_{k=0,1,\cdots,[Tn^\theta]}   \Big|\int_{k/n^\theta}^{(k+1)/n^\theta}\int_0^\infty\int_0^{V^{(n)}(s-)}\frac{1}{n^\alpha} \widetilde{N}_{\mathtt{l}}^{(n)}( ds,  dy,  \cdot dz)\Big|>\eta \bigg)} \ar\ar\cr
    	\ar\leq\ar \sum_{k=0}^{[Tn^\theta]}\mathbf{P}\bigg(
    	\Big|\int_{k/n^\theta}^{(k+1)/n^\theta}\int_0^\infty\int_0^{V^{(n)}(s-)}\frac{1}{n^\alpha} \widetilde{N}_{\mathtt{l}}^{(n)}( ds,  dy,  \cdot dz)\Big|>\eta\bigg)\cr
    	\ar\leq\ar \frac{1}{\eta^{2p}}\sum_{k=0}^{[Tn^\theta]}\mathbf{E}
    	\bigg[\Big|\int_{k/n^\theta}^{(k+1)/n^\theta}\int_0^\infty\int_0^{V^{(n)}(s-)}\frac{1}{n^\alpha} \widetilde{N}_{\mathtt{l}}^{(n)}( ds,  dy,  \cdot dz)\Big|^{2p}\bigg]\cr
    	\ar\leq\ar C\cdot n^\theta \cdot \sup_{k=0,1,\cdots,[Tn^\theta]}\mathbf{E}
    	\bigg[\Big|\int_{k/n^\theta}^{(k+1)/n^\theta}\int_0^\infty\int_0^{V^{(n)}(s-)}\frac{1}{n^\alpha} \widetilde{N}_{\mathtt{l}}^{(n)}( ds,  dy,  \cdot dz)\Big|^{2p}\bigg]\cr
    	\ar\leq\ar  Cn^\theta \sup_{k=0,1,\cdots,[Tn^\theta]}\bigg( \Big| \int_{k/n^\theta}^{(k+1)/n^\theta} n^{1-\alpha} ds\Big|^p
    	+ \int_{k/n^\theta}^{(k+1)/n^\theta} n^{1+\alpha-2\alpha p} ds \bigg)\cr
    	\ar\ar\cr
    	\ar\leq\ar C\cdot \Big(n^{p(1-\alpha-\theta)+\theta}+n^{1+\alpha-2\alpha p}\Big),
    \eeqnn 
     which converges to zero as $n \to \infty$ for all $p>2$. Here the second last inequality follows from Proposition \ref{2pine} and the fact that $\nu(0,\infty)= 1$.

 \item  For any $\eta>0$, we use Chebyshev's inequality again to obtain that 
 \beqnn
 \mathbf{P}\bigg(
 \sup_{k=0,1,\cdots,[Tn^\theta]} \int_{k/n^\theta}^{(k+1)/n^\theta} 
 n\cdot V^{(n)}(s)  ds  > \eta  \bigg)
 \ar\leq\ar \frac{n^2}{\eta^2}\cdot  \sum_{k=0}^{[Tn^\theta]} 
 \mathbf{E}  \bigg[
 \Big|\int_{k/n^\theta}^{(k+1)/n^\theta} 
  V^{(n)}(s)  ds\Big|^2\bigg].
 \eeqnn 
 By H{\"o}lder's inequality and Lemma \ref{V2p}, we see that 
 \beqnn
 \mathbf{E}  \bigg[
 \Big|\int_{k/n^\theta}^{(k+1)/n^\theta} 
 V^{(n)}(s)  ds\Big|^2\bigg] 
 \leq \frac{1}{n^\theta} \cdot 
 \int_{k/n^\theta}^{(k+1)/n^\theta} 
 \mathbf{E}  \big[|V^{(n)}(s)|^2\big]   ds  \leq \frac{C}{n^{2\theta}},
 \eeqnn
 for some constant $C>0$ independent of  $k$ and $n$. 
 Hence as $n\to\infty$,
  \beqnn
   \mathbf{P}\bigg(
  \sup_{k=0,1,\cdots,[Tn^\theta]} \int_{k/n^\theta}^{(k+1)/n^\theta} 
  n\cdot V^{(n)}(s)  ds  > \eta  \bigg) 
  \leq \frac{C}{n^{\theta-2}} \to 0.
  \eeqnn  
\end{itemize}

 Let us now consider the sequence $\{\varepsilon_{12}^{(n)}\}_{n \geq 1}$. Similar to the above, we can derive the following upper bound: 
 \beqnn
 \sup_{k=0,1,\cdots,[Tn^\theta]}\sup_{h\in[0,1/n^{\theta}]}
 \Big|\varepsilon_{12}^{(n)}(k/n^\theta,h)\Big|
 \ar\leq\ar 
 \sup_{k=0,1,\cdots,[Tn^\theta]}
 \bigg|\int_0^{k/n^\theta} \int_{T}^{(k+1)/n^\theta }\int_0^{V^{(n)}(s-)}
 \frac{1}{n^\alpha}\widetilde{N}_{\mathtt{1}}^{(n)}( ds, dy, dz)\bigg| \\
\ar\ar +
 2\cdot 
 \int_0^{T}n\cdot V^{(n)}(s)ds \cdot \sup_{k=0,1,\cdots,[Tn^\theta]} \int_{k/n^\theta}^{(k+1)/n^\theta }
 \nu(n\cdot dy)  . 
 \eeqnn
 \begin{itemize}
 \item For any $\eta>0$, by Chebyshev's inequality and  Proposition \ref{2pine} we have 
 \beqlb\label{est1}
 \lefteqn{\mathbf{P} \bigg(\sup_{k=0,1,\cdots,[Tn^\theta]}
 \Big|\int_0^{k/n^\theta} \int_{k/n^\theta}^{(k+1)/n^\theta }\int_0^{V^{(n)}(s-)}
 \frac{1}{n^\alpha}\widetilde{N}_{\mathtt{1}}^{(n)}( ds, dy, dz)\Big|>\eta \bigg)}\ar\ar\cr
 \ar\leq\ar C\cdot n^\theta \cdot  \sup_{k=0,1,\cdots,[Tn^\theta]} \mathbf{E} \bigg[ 
 \Big|\int_0^{k/n^\theta} \int_{k/n^\theta}^{(k+1)/n^\theta }\int_0^{V^{(n)}(s-)}
 \frac{1}{n^\alpha}\widetilde{N}_{\mathtt{1}}^{(n)}( ds, dy, dz)\Big|^{2p} \bigg]\cr
 \ar\leq\ar C\cdot n^\theta \cdot  
 \Big|\frac{T}{n^{\alpha-1}}\cdot  \sup_{k=0,1,\cdots,[Tn^\theta]}   \int_{k/n^\theta}^{(k+1)/n^\theta } \nu(n\cdot dy) \Big|^{p} \cr
 \ar\ar + C\cdot n^\theta \cdot   
 \frac{T}{n^{2p\alpha-(\alpha+1)}}  \sup_{k=0,1,\cdots,[Tn^\theta]}  \int_{k/n^\theta}^{(k+1)/n^\theta } \nu(n\cdot dy) . 
 \eeqlb  
 Recalling that $\nu(dy)=(\alpha+1)(1+y)^{-\alpha-2}dy$, we see that
 \beqlb\label{eqn.201}
 \sup_{k=0,1,\cdots,[Tn^\theta]} \int_{k/n^\theta}^{(k+1)/n^\theta}
 \nu(n\cdot dy) \leq  \frac{\alpha+1}{n^{\theta-1}}.
 \eeqlb 
 Plugging this back into the right-hand side of (\ref{est1}), the probability can be bounded by
 \beqnn  
 \frac{C}{n^{p(\theta+\alpha-2)-\theta}}  + \frac{C}{n^{2p\alpha-\alpha-2}}  ,
 \eeqnn
 which goes to $0$ as $n\to \infty$ for all $p>2\theta$. 
 
\item
  For any $\eta>0$, by using Chebyshev's inequality, (\ref{eqn.201}) and then Lemma \ref{V2p} we have 
  \beqnn
 \lefteqn{ \mathbf P \bigg ( 
 \int_0^{T}n\cdot V^{(n)}(s)ds \cdot \sup_{k=0,1,\cdots,[Tn^\theta]} \int_{k/n^\theta}^{(k+1)/n^\theta }
 \nu(n\cdot dy) > \eta \bigg) }\ar\ar\cr
 \ar\leq\ar \mathbf P \bigg ( 
 \int_0^{T}n\cdot V^{(n)}(s)ds \cdot \frac{\alpha+1}{n^{\theta-1}} > \eta \bigg) 
 \leq \frac{1}{\eta}\cdot \frac{\alpha+1}{n^{\theta -2}}   \int_0^{T} \mathbf{E}\big[V^{(n)}(s)\big]ds \leq  \frac{C}{n^{\theta-2}},
  \eeqnn
  which goes to $0$ as $n\to\infty$, since $\theta>2$.
  \qed 
  \end{itemize}

%% file: 4.AccumulationNew.tex
 \section{Accumulation points} 
 \label{sec:accummulation}
 \setcounter{equation}{0}

 In this section, we characterize the accumulation points of the sequence of rescaled volatility processes. By identifying the accumulation points of each term in \eqref{decomV} we obtain a stochastic Volterra equation that each accumulation point satisfies. 

 Following the same arguments as in the proof of \cite[Proposition 3.13]{HorstXuZhang2023a}, we can identify the accumulation points of the sequence $\{J^{(n)}_1 \}_{n\geq 1}$ in terms of the Mittag-Leffler density function $f^{\alpha, \gamma}$.

 \begin{proposition}[Proposition 3.13 in \cite{HorstXuZhang2023a}] \label{Prop.3.13}
 Any accumulation point $(V_*,J_{1,*}) \in \mathbf{C}(\mathbb{R}_+;\mathbb{R}_+\times \mathbb R)$ of the sequence $\{ (V^{(n)},J_1^{(n)} ) \}_{n\geq 1}$ is a weak solution to the stochastic equation 
 \begin{equation}\label{V*} 
 J_{1,*}(t)=
 \int_0^t  f^{\alpha,\gamma}(t-s) \cdot \frac{ \sqrt{\lambda^\mathtt m_*}}{b}\cdot\sqrt{V_*(s)} dB_s,\quad t\geq 0.
 \end{equation}
 \end{proposition}

 Identifying the accumulation points of the sequence $\{J^{(n)}_2 \}_{n\geq 1}$ is more challenging. 
 The process $J^{(n)}_2$ has been decomposed in Section \ref{ctightVn} as 
 \beqnn
 J^{(n)}_2=J_{2,c}^{(n)} + \varepsilon^{(n)}_1 .
 \eeqnn
 We refine this decomposition by further splitting the integral process as 
 \beqnn
 J^{(n)}_2=\hat{J}^{(n)}_2+\varepsilon^{(n)}_1+\varepsilon^{(n)}_2,
 \eeqnn
 where the two processes $\{\hat{J}^{(n)}_2 \}_{n\geq 1}$ and $\{\varepsilon^{(n)}_2 \}_{n\geq 1}$ are defined as follows: for $t\geq 0$,
 \beqnn
 \hat{J}_2^{(n)}(t)
 \ar:=\ar
 \int_0^t\int_0^{\infty}\int_0^{V^{(n)}(s-)}
 \Big(\int_{(t-s-y)^+}^{t-s}\frac{ f^{\alpha,\gamma}(r)}{b} dr\Big)
 \widetilde{N}_\mathtt l^{(n)}(ds, dy,dz), \cr 
 \varepsilon^{(n)}_2(t)\ar:=\ar 
 \int_0^t\int_0^{\infty}\int_0^{V^{(n)}(s-)}
 \Big(\int_{(t-s-y)^+}^{t-s} R_H^{(n)}(r)dr-\int_{(t-s-y)^+}^{t-s}\frac{f^{\alpha,\gamma}(r)}{b} dr
 \Big)\widetilde{N}_{\mathtt{l}}^{(n)}( ds, dy, dz).
 \eeqnn
 
 The $C$-tightness of the sequences $\{J_{2,c}^{(n)} \}_{n\geq 1}$ and $\{\varepsilon^{(n)}_1\}_{n\geq 1}$ has already been established. 
 By Lemma 5.4 in \cite{Xu2024b}, the sequence $\{\hat{J}_2^{(n)}\}_{n\geq 1}$ is also $C$-tight. Hence, the sequence 
 \beqnn
 \big\{\varepsilon^{(n)}_2\big\}_{n\geq 1}=\big\{J_{2,c}^{(n)}-\hat{J}^{(n)}_2 \big\}_{n\geq 1}
 \eeqnn
 is $C$-tight as well. 
 Furthermore, using the same argument as in the proof of \cite[Proposition 4.12]{Xu2024b} it is not difficult to show that for any $T>0$,
 \beqnn
    \lim_{n\to\infty}\sup_{t\in[0,T]}\mathbf{E}\Big[\big|\varepsilon^{(n)}_i(t)\big|^2\Big]= 0,\quad i=1,2.
 \eeqnn
 This shows that the two sequences $\{\varepsilon^{(n)}_i\}_{n\geq 1}$, $i=1,2$, converge weakly to the zero process in $\mathbf D(\mathbb{R}_+;\mathbb{R})$. 
 It hence remains to identify the weak accumulation points of the sequence $\{\hat{J}^{(n)}_2\}_{n\geq 1}$.


 \subsection{Weak convergence of stochastic Volterra integrals}

 To establish our weak convergence result, we generalize the weak convergence result established in \cite[Section 4.2]{Xu2024b} for stochastic Volterra integrals by applying the general theory of weak convergence of It\^o's stochastic integrals with respect to infinite-dimensional semimartingales, due to Kurz and Protter \cite{KurtzProtter1996} that we briefly recall for the readers' convenience. 
 
 \subsubsection{Infinite-dimensional stochastic integration}
 
 Let $\mathbb{H}$ be a separable Banach space endowed with a norm $\|\cdot\|_{\mathbb{H}}$.  We first recall the definition of $\mathbb{H}^\#$-semimartingales and the corresponding stochastic integrals. To this end, we denote 
 by $\mathcal{S}_0$ be the collection of processes of the form
 \beqnn
 X(t):=\sum_{k=1}^m \xi_{k}(t)\varphi_k ,
 \eeqnn
 where $m\in\mathbb{Z}_+$, $\{ \varphi_k \}_{k\geq 1}\subset \mathbb{H}$ and $\xi_{k}$ is a process of the form
 \beqnn
 \xi_{k}(t):= \sum_{i=0}^\infty \eta_{k,i}\cdot \mathbf{1}_{[\tau_{k,i},\tau_{k,i+1})}(t),
 \eeqnn
 in which $0=\tau_{k,0}\leq\tau_{k,1}\leq\cdots$ are $(\mathscr{F}_t)$-stopping times and $\eta_{k,i}$ is a  $\mathscr{F}_{\tau_{k,i}}$-measurable $\mathbb{R}$-valued random variable. 

 \begin{definition}
 We say	$Y$ is a $(\mathscr{F}_t)$-adapted \textit{$\mathbb{H}^\#$-semimartingale}, if it is a $\mathbb{R}$-valued stochastic process indexed by $\mathbb{H} \times \mathbb{R}_+$ such that:
 \begin{enumerate}
 		\item[$\bullet$] for each $\varphi\in \mathbb{H}$,  $Y(\varphi):=\{ Y(\varphi,t):t\geq 0 \}$ is a c\`adl\`ag $(\mathscr{F}_t)$-semimartingale with $Y(\varphi,0)\overset{\rm a.s.}=0$;
 		
 		\item[$\bullet$] for each $t\geq 0$, $a_1,\cdots,a_m\in\mathbb{R}$ and $\varphi_1,\cdots,\varphi_m\in\mathbb{H}$, 
 		\beqnn
 		Y\bigg(\sum_{k=1}^m a_k \varphi_k,t\bigg)\overset{\rm a.s.}=\sum_{k=1}^m a_k Y(\varphi_k,t).
 		\eeqnn
 	\end{enumerate}

 \end{definition}
 
 For $X\in\mathcal{S}_0$, its stochastic integral with respect to $Y$ is defined by
 \beqnn
 X_-\cdot Y(t) =\sum_{k=1}^m \int_0^t \xi_k(s-)d Y(\varphi_k,s), \quad t\geq 0.
 \eeqnn 
The definition of the stochastic integral can be extended to all c\`adl\`ag, $\mathbb{H}$-valued stochastic processes $X$ by approximation; see  \cite[Page 226]{KurtzProtter1996} for details. To introduce the notion of a standard $\mathbb{H}^\#$-semimartingale we introduce for $t\geq 0$ the set  
 \beqlb\label{eqn.Ht}
 \mathcal{H}_t(Y):=  \Big\{ \sup_{s\leq t} \big|X_-\cdot Y(s)\big| : X \in\mathcal{S}_0,\, \sup_{s\leq t} \big\|X(s) \big\|_{\mathbb{H}}\leq 1  \Big\}.
 \eeqlb

 \begin{definition}
  The $\mathbb{H}^\#$-semimartingale $Y$ is  {\rm standard} if $\mathcal{H}_t(Y)$
  is stochastically bounded for each $t\geq 0$. 
  Moreover,  a sequence of $\mathbb{H}^\#$-semimartingles $\{ Y^{(n)} \}_{n\geq 1}$ is said to be {\rm uniformly tight} if $\cup_{n\geq 1} \mathcal{H}_t(Y^{(n)})$ is stochastically bounded for each $t\geq 0$. 
 \end{definition}
 
 For $k,d\in\mathbb{Z}_+$, a process $\boldsymbol{X} := (X_{ij})_{i\leq k,j\leq d} \in \mathbf{D}(\mathbb{R}_+; \mathbb{H}^{k\times d})$ and a $d$-dimensional $(\mathscr{F}_t)$-adapted  $\mathbb{H}^\#$-semimartingale $\boldsymbol{Y} := (Y_{j})_{j\leq d}$, the stochastic integral $\boldsymbol{X}_-\cdot\boldsymbol{Y}$ is defined by
 \beqnn
 \boldsymbol{X}_-\cdot \boldsymbol{Y}(t):= \bigg( \sum_{j=1}^d X_{ij,-}\cdot Y_{j}(t)  \bigg)_{i\leq k},\quad t\geq 0. 
 \eeqnn 

 The following convergence result for infinite-dimensional stochastic integrals has been established in \cite{KurtzProtter1996}.

 \begin{proposition}\label{kurztheorem} 
 For $n\geq 1$, let $\boldsymbol{X}^{(n)} \in  \mathbf{D}(\mathbb{R}_+; \mathbb{H}^{d\times k})$ and  $\boldsymbol{Y}^{(n)}$ be a standard  $k$-dimensional $\mathbb{H}^\#$-semimartingale.   
 If the sequence $\{\boldsymbol{Y}^{(n)}\}_{n \geq 1}$ is uniformly tight and 
 \beqnn
 	\big(\boldsymbol{X}^{(n)},\boldsymbol{Y}^{(n)}\big)\Rightarrow \big(\boldsymbol{X},\boldsymbol{Y}\big)
 \eeqnn
 in the sense that for each choice of $\varphi_1,\cdots,\varphi_m\in \mathbb{H}$, 
 \beqnn
 \big(\boldsymbol{X}^{(n)}, \boldsymbol{Y}^{(n)}(\varphi_1,\cdot),\cdots,\boldsymbol{Y}^{(n)}(\varphi_m,\cdot)\big) \to \big(\boldsymbol{X},\boldsymbol{Y}(\varphi_1,\cdot),\cdots,\boldsymbol{Y}(\varphi_m,\cdot)\big),
 \eeqnn
 weakly in $\mathbf{D}(\mathbb{R}_+; \mathbb{H}^{d\times k}\times \mathbb{R}^{k\times m})$, 
 then there exists a filtration $\{\mathscr{F}_t:t\geq 0\}$ such that $\boldsymbol{X}$ is $(\mathscr{F}_t)$-adapted, $\boldsymbol{Y}$ is an $(\mathscr{F}_t)$-adapted, standard, $k$-dimensional $\mathbb{H}^\#$-semimartingale and 
 \beqnn
 (\boldsymbol{X}^{(n)}, \boldsymbol{Y}^{(n)},\boldsymbol{X}^{(n)}_-\cdot \boldsymbol{Y}^{(n)})\Rightarrow (\boldsymbol{X}, \boldsymbol{Y},\boldsymbol{X}_-\cdot \boldsymbol{Y}). 
 \eeqnn 
 \end{proposition}

\subsubsection{Weak limit of  $\{\hat{J}^{(n)}_2\}_{n\geq 1}$}

 To identify the weak limit of the sequence $\{\hat{J}^{(n)}_2\}_{n\geq 1}$, we first generalize Proposition~\ref{kurztheorem} to a sequence of stochastic Volterra integrals $\{ Z^{(n)} \}_{n\geq 1}$ with  
 \beqlb\label{int-gen}
 Z^{(n)}(t):= \int_0^t \int_0^\infty \int_0^{\xi^{(n)}(s-)}G(t,s,y)\widetilde{N}^{(n)}(ds,dy,dz), \quad t \geq 0
 \eeqlb 
 for some process $\xi^{(n)} \in \mathbf{D}(\mathbb{R}_+;\mathbb{R}_+)$, a function $G$ on $(0,\infty)^3$ and a compensated Poisson random measure $\widetilde N^{(n)}(ds,dy,dz)$ on $[0, \infty)^3$  with intensity $ds \, \mu_0^{(n)}(dy)\,dz$ where $\mu_0^{(n)}(dy)$ is a $\sigma$-finite measure on $(0,\infty)$. We assume that all processes are defined on a common probability space that supports an additional compensated Poisson process $\widetilde N_*(ds,dy,dz)$ with intensity measure $ds\,\mu^*_0(dy)\,dz$ with $\mu_0^*(dy)$ being a $\sigma$-finite measure on $(0,\infty)$.

 \begin{condition} \label{ass:L2}
We assume that the following two cconditions hold.
  
  \begin{enumerate}
  	\item[(1)] There exists a constant $C>0$ such that for any non-negative measurable function $f$ on $\mathbb{R}_+$,
  	\beqnn
  	\sup_{n\geq 1} \int_0^\infty f(y)  \mu^{(n)}_0(dy) \leq C\cdot \int_0^\infty f(y)  \mu^*_0(dy). 
  	\eeqnn
  	
  	\item[(2)]  The function $G$ is continuous in the second argument and  
  	\beqlb \label{L2 bounded condition}
    \int_0^T \int_0^\infty \big| G(t,s,y) \big|^2 \mu^*_0(dy)ds<\infty, \quad t,T\geq 0,\, n\geq 1. 
  	\eeqlb
  \end{enumerate}  
 \end{condition}
 
 Let $\mu^*(dy,dz) :=\mu^*_0(dy)\,dz$ and $L^2(\mu^*)$ be the Hilbert space of square integrable functions on $(0,\infty)^2$ with respect to $\mu^*$ endowed with the norm $\|\cdot\|_{L^2(\mu^*)}$.
 In terms of the compensated random measures we can introduce the following measure valued processes: 
 \beqnn
  {\bf \widetilde N}^{(n)}(t) := \widetilde N^{(n)}\big((0,t],dy,dz\big) 
 \quad \mbox{and} \quad 
 {\bf \widetilde N}^*(t) := \widetilde N^*\big((0,t],dy,dz\big),
 \quad t \geq 0.
 \eeqnn  
 Under Condition \ref{ass:L2} the above processes are standard $L^2(\mu^*)^\#$-semimartigales. 
 In particular, the notion of weak convergence of $\{{\bf \widetilde N}^{(n)}\}_{n \geq 1}$ is well defined, and we can formulate our weak convergence result for stochastic Volterra integrals.
 
 \begin{lemma}\label{Lemma.401}
 Suppose that Condition~\ref{ass:L2} holds, that $\big\{ {\bf \widetilde N}^{(n)} \big\}_{n \geq 1}$ is uniformly tight and that for any $T\geq 0$,
 \beqlb\label{eqn.402}
 \sup_{t\in[0,T]}\sup_{n\geq 1}\mathbf{E}\big[\xi^{(n)}(t)\big]<\infty.
 \eeqlb
 If  $ \big(Z^{(n)},\xi^{(n)},{\bf \widetilde N}^{(n)}\big) \Rightarrow \big(Z^*, \xi^*,{\bf \widetilde N^*} \big)$, then the limit  satisfies that
 \beqnn
 Z^*(t)= \int_0^t \int_0^\infty \int_0^{\xi^*(s-)}G(t,s,y)\widetilde N^*(ds,dy,dz), \quad t \geq 0. 
 \eeqnn
 \end{lemma}
 \proof Due to the dependence of the integrand $F$ on the time variable we cannot directly utilize Theorem \ref{kurztheorem} to prove the desired weak convergence.
 To overcome this problem we utilize the tightness of the sequence $\{Z^{(n)}\}_{n \geq 1}$. 
 Due to the tightness it is enough to establish the convergence of the finite dimensional distributions, i.e.~to show that
 \beqlb\label{finite dimensional distribution Zd}
 \Big (Z^{(n)}(t_1),\cdots, Z^{(n)}(t_d)\Big)
 \to
 \Big (Z^*(t_1),\cdots, Z^*(t_d)\Big)
 \eeqlb 
 in distribution for all $d\in\mathbb N_+$, and for all $t_1<\cdots <t_d$. 
 This allows us to ``drop'' the dependence of  integrand on the time variable and to consider instead the integral processes:
 \beqnn
 Z_i^{(n)}(t)
 \ar:=\ar \int_0^t\int_0^\infty \int_0^{\xi^{(n)}(s-)}G(t_i,s,y)\widetilde N^{(n)}(ds,dy,dz),\cr
 Z^*_i(t)
 \ar:=\ar\int_0^t\int_0^\infty\int_0^{\xi_*(s-)}G(t_i,s,y)\widetilde N^* (ds,dy,dz),
 \eeqnn
  for $i=1,\cdots, d$ and $t\geq 0$. 
 In fact, since $Z^{(n)}(t_i)\overset{\rm a.s.}=Z_i^{(n)}(t_i)$ and $Z^*(t_i)\overset{\rm a.s.}=Z^*_i(t_i)$ for all $i=1,\cdots,d$, it suffices to prove that 
 \beqlb\label{eqn.401}
 \big( Z^{(n)}_1, ..., Z^{(n)}_d \big) 
 \to
 \big( Z^*_{1}, ..., Z^*_d \big),
 \eeqlb
 weakly in $\mathbf{D}(\mathbb{R}_+;\mathbb{R}^d)$ as $n\to\infty$. 
 
 Let $T>0$.  
 Under Condition~\ref{ass:L2}(2), for each $i=1,\cdots, d$ and $\varepsilon\in(0,1)$ we can find a measurable function $G_i^\varepsilon(s,y)$ on $(0,\infty)^2$ that is continuous in $s$ and satisfies 
 \beqlb\label{eqn.404}
 \int_0^T \int_0^\infty\big|G(t_i,s ,y)-G_i^\varepsilon(s,y)\big|^2\mu^*_0(dy)ds\leq \varepsilon
 \quad \mbox{and}\quad 
\int_0^\infty  \sup_{s\in[0,T]}\big|G_i^\varepsilon(s,y)\big|^2\mu^*_0(dy)<\infty. 
 \eeqlb
 Associated to these functions, for $i=1,\cdots,d$ and $n\geq 1$ we introduce the processes
 \beqnn
 Z_{i,\varepsilon}^{(n)}(t)
 \ar:=\ar \int_0^t\int_0^\infty \int_0^{\xi^{(n)}(s-)}G_i^\varepsilon(s,y)\widetilde N^{(n)}(ds,dy,dz),\cr
 Z^*_{i,\varepsilon}(t)
 \ar:=\ar\int_0^t\int_0^\infty\int_0^{\xi_*(s-)}G_i^\varepsilon(s,y)\widetilde N^* (ds,dy,dz). 
 \eeqnn
 By the Burkholder-Davis-Gundy inequality, (\ref{eqn.402}) and Condition~\ref{ass:L2}(1), we have uniformly in $n\geq 1$,
 \beqnn
 \mathbf{E}\Big[ \sup_{t\in[0,T]} \big|Z_{i,\varepsilon}^{(n)}(t)-Z_{i}^{(n)}(t)\big|^2 \Big]
 \ar\leq\ar C\cdot  \int_0^T \int_0^\infty\big|G(t_i,s ,y)-G_i^\varepsilon(s,y)\big|^2\mu^*_0(dy)ds,
 \eeqnn
 which goes to $0$ as $\varepsilon\to0$. Similarly,  $\mathbf{E}\big[ \sup_{t\in[0,T]} \big|Z_{i,\varepsilon}^*(t)-Z_{i}^*(t)\big|^2 \big]\to0$ as $\varepsilon\to0$. 
 By Theorem 3.2 in \cite[p.28]{Billingsley1999}, the weak convergence (\ref{eqn.401}) holds if for each $\varepsilon\in(0,1)$, 
 \beqlb\label{eqn.40101}
 \big( Z^{(n)}_{1,\varepsilon}, ..., Z^{(n)}_{d,\varepsilon} \big) 
 \to
 \big( Z^*_{1,\varepsilon}, ..., Z^*_{d,\varepsilon} \big),
 \eeqlb
 weakly in $\mathbf{D}(\mathbb{R}_+;\mathbb{R}^d)$ as $n\to\infty$. 
 
 We now apply Proposition~\ref{kurztheorem} with $k=1$ and $\mathbb{H}=L^2(\mu^*)$ to prove (\ref{eqn.401}).
 To this end,  we fix $\varepsilon\in(0,1)$ and introduce, for each $i=1,\cdots, d$, the $L^2(\mu^*)$-valued function 
 \beqnn
 F_i: \mathbf{D}(\mathbb{R}_+;\mathbb{R}_+)  \times [0,\infty)\to L^2(\mu^*) 
 \quad
 (x,s) \mapsto F_i(x,s)(\cdot,\cdot) 
 \eeqnn    
 where the mapping $F_i(x,s)$ is defined by
 $ F_i(x,s)(y,z):= G_i^\varepsilon(s,y)\cdot \mathbf{1}_{\{0<z\leq x(s)\}}$. 
We notice that 
\[
	F_i (\xi^{(n)},\cdot ) \in \mathbf{D}(\mathbb{R}_+;L^2(\mu^*)).
\]
 In terms of these mapping we can represent the integral process $Z^{(n)}_{i,\varepsilon}$ as 
 \beqnn
 Z_{i,\varepsilon}^{(n)}(t) = F_i\big(\xi^{(n)},-\big)\cdot {\bf \widetilde N}^{(n)}(t),\quad t\geq 0.
 \eeqnn
 To obtain the weak convergence (\ref{eqn.40101}), by Proposition~\ref{kurztheorem} it suffices to prove that 
 \beqnn
 \Big(F_1\big(\xi^{(n)},\cdot\big),\cdots, F_d\big(\xi^{(n)},\cdot \big), {\bf \widetilde N}^{(n)} \Big) \Rightarrow 
 \Big(F_1\big(\xi^*,\cdot\big),\cdots, F_d\big(\xi^*,\cdot \big), {\bf \widetilde N}^* \Big).
 \eeqnn
 
 By the continuous mapping theorem, it suffices to show that $F_i$ is a continuous mapping from $\mathbf{D}(\mathbb{R}_+;\mathbb{R}_+)$ to $\mathbf{D}(\mathbb{R}_+;L^2(\mu^*))$.
 Indeed, if $x_n\to x_*$ in $\mathbf{D}(\mathbb{R}_+;\mathbb{R}_+)$ as $n\to\infty$, by the definition of Skorokhod topology, there exists a sequence of strictly increasing, continuous bijective functions $\{\gamma_n\}_{n\geq 1}$ such that as $n\to\infty$
 \beqlb\label{eqn.403}
 \sup_{0\leq t\leq T}\big|\gamma_n(t)-t\big|\to 0
 \quad \mbox{and}\quad 
 \sup_{0\leq t\leq T} \big|x_n(\gamma_n(t)) - x_*(t)\big| \to 0,
 \eeqlb  
 which also induces that for some constant $C>0$ that is independent of $n\geq 1$ and $t\in[0,T]$,
 \beqlb\label{eqn.430}
 \big|x_n\big(\gamma_n(t)\big)\big| + \big|x_*(t)\big| \leq C.
 \eeqlb 
 Using the $L^2(\mu^*)$-boundedness of $G$ and its continuity in the second argument we conclude that
 \beqnn
  \lefteqn{\sup_{t\in[0,T]}\Big\|F_i\big(x_n,\gamma_n(t)\big)  - F_i\big(x_*,t\big) \Big\|_{L^2(\mu^*)}^2}\ar\ar\cr
 \ar=\ar \sup_{t\in[0,T]}\int_0^\infty  \int_0^\infty \big| G_i^\varepsilon\big( \gamma_n(t),y\big)\cdot \mathbf{1}_{\{0<z\leq x_n(\gamma_n(t))\}} -G_i^\varepsilon(t,y)\cdot \mathbf{1}_{\{0<z\leq x_*(t)\}} \big|^2 \mu^*_0(dy)dz\cr
 \ar\leq\ar \sup_{t\in[0,T]}\int_0^\infty  \int_0^\infty \big| G_i^\varepsilon\big(\gamma_n(t),y\big)  -G_i^\varepsilon(t,y) \big|^2 \cdot \mathbf{1}_{\{0<z\leq x_n(\gamma_n(t))\}} \mu^*_0(dy)dz\cr 
 \ar\ar +\sup_{t\in[0,T]} \int_0^\infty  \int_0^\infty \big| G_i^\varepsilon(t,y)\big|^2\cdot\big|   \mathbf{1}_{\{0<z\leq x_n(\gamma_n(t))\}} - \mathbf{1}_{\{0<z\leq x_*(t)\}} \big| \mu^*_0(dy)dz\cr
 \ar\leq\ar \sup_{s\in[0,T]}\big| x_n\big(\gamma_n(s)\big) \big| \cdot  \int_0^\infty \sup_{t\in[0,T]} \big|G_i^\varepsilon\big(\gamma_n(t),y\big)  -G_i^\varepsilon(t,y) \big|^2    \mu^*_0(dy) \cr 
 \ar\ar + \sup_{s\in[0,T]}\big| x_n\big(\gamma_n(s)\big)-x_*(s) \big| \cdot  \sup_{t\in[0,T]} \int_0^\infty \big| G_i^\varepsilon(t,y)\big|^2 \mu^*_0(dy) .
 \eeqnn 
 
 By \eqref{eqn.404} and \eqref{eqn.430}, the second term on the right-hand side of the last inequality vanishes as $n\to\infty$. 
 For the first term, by the dominated convergence theorem along with the local continuity of $G_i^\varepsilon(s,y)$ in $s$ as well as \eqref{eqn.404} and \eqref{eqn.403} it also vanishes as $n\to\infty$. As a result, $F_i$ is a continuous mapping. 
 \qed

 \begin{remark}
 We emphasize that the weak convergence, respectively the tightness of the sequence $\{\xi^{(n)}\}_{n\geq 1}$ is key to our argument. Without knowing a priori that the sequence is tight/converges, the argument does not apply. 
 \end{remark}
 
 We now use Lemma~\ref{Lemma.401} to characterize the accumulation points  of the relativly compact sequence $\{ \hat{J}_2^{(n)} \}_{n\geq 1}$ by setting 
 \beqnn
 \xi^{(n)}= V^{(n)}, \quad 
 \mu^{(n)}_0(dy)=  \frac{\lambda^{\tt l}_n}{n^{1-\alpha}}\cdot n^{\alpha+1}\cdot \alpha\cdot \nu(n\cdot dy),\quad  
 \mu^*_0(dy)=\lambda^{\tt l}_*\cdot \nu_*(dy)
 \eeqnn
 and 
 \beqnn
 G(t,s,y)=\int_{(t-s-y)^+}^{t-s} \frac{f^{\alpha,\gamma}(r) }{b} dr.
 \eeqnn
 By (\ref{eqn.210}), it is easy to see that Condition~\ref{ass:L2} is satisfied and the following result is a direct consequence of Lemma~\ref{Lemma.401}. 
 
 \begin{corollary} [Accumulation points of $\{\hat{J}_2^{(n)}\}_{n\geq 1}$]
 Any accumulation point $(V_*,J_{2,*}) $ of the sequence 
 \[
 	\big\{(V^{(n)},\tilde{J}^{(n)}_2)\big\}_{n\geq 1} 
\]	
is a weak solution to the stochastic integral equation
 \beqlb	\label{wc4}
  J_{2,*}(t)=\int_0^t\int_0^{\infty}\int_0^{V_*(s)}
  \Big(\int_{(t-s-y)^+}^{t-s} \frac{f^{\alpha,\gamma}(r)}{b}  dr\Big)\widetilde{N}(ds,dy,dz), \quad t\geq 0.
 \eeqlb 
  \end{corollary} 
 

\subsection{Proof of Theorem \ref{thm1}}

 Our preceding results show that the sequence
 $ \{ (V^{(n)}, I^{(n)}, J_1^{(n)},J_2^{(n)}, \zeta^\mathtt{m}_n/\beta_n,\zeta^\mathtt{l}_n/ n^{\alpha-1} ) \}_{n \geq 1} $
 is $C$-tight. 
 Moreover, the previous lemmas and Condition~\ref{main.Condition} show that each accumulation point is of the form 
 \beqnn
 \ar\ar\Big( V_*,V_*(0)\cdot (1-F^{\alpha,\gamma})+\frac{a}{b}F^{\alpha,\gamma}, \int_0^\cdot 
  f^{\alpha,\gamma}(\cdot-s)\cdot \frac{ \sqrt{\lambda^\mathtt{m}_*}}{b} \cdot \sqrt{V_*(s)}dB_s,\cr
 \ar\ar\qquad   \int_0^\cdot 
 \int_0^{\infty}\int_0^{V_*(s)} \Big(\int_{(\cdot-s-y)^+}^{\cdot-s}
 \frac{ f^{\alpha,\gamma}(r)}{b} dr\Big) \widetilde{N}(ds,dy,dz),\zeta^\mathtt{m}_* ,\zeta^\mathtt{l}_* \Big).
 \eeqnn
In particular, the accumulation point satisfies (\ref{eq1}). To check that the two representations \eqref{eq1} and \eqref{eq2} are equivalent, we assume that \eqref{eq2} holds and rewrite this equation as
 \beqlb\label{V'1}
  V_*(t)=\tilde{V}(t)-K*V_*(t),\quad t\geq 0,
 \eeqlb
 with 
 \beqlb 
 \tilde{V}(t)\ar:=\ar V_*(0)+ \int_0^t\frac a b\cdot K(t-s)ds 
 +\int_0^t  K(t-s)\cdot\frac{\zeta^*_\mathtt{m}\sqrt{\lambda^*_\mathtt m}}{b}\cdot\sqrt{V_*(s)}dB_s\\
 \ar\ar  +\int_0^t\int_0^{\infty}\int_0^{V_*(s)}
 \Big( \int_{(t-s-y)^+}^{t-s}
 \frac{\zeta^*_\mathtt{l}}{b}\cdot K(r)dr\Big) \widetilde{N}(ds,dy,dz),\quad t\geq 0.
 \eeqlb
 
 We now recall that the Mittag-Leffler density function $f^{\alpha,\gamma}$ is the resolvent of the second kind of the function $K$; see \eqref{LK}. 
 By applying Lemma~3 in \cite{BacryDelattreHoffmannMuzy2013} to the equation (\ref{V'1}), we see that it can be rewritten as
 \beqlb
  V_*(t)=\tilde{V}(t)-f^{\alpha,\gamma}*\tilde{V}(t), \quad t\geq 0.
 \eeqlb
 Similarly as in the argument in Section~7 in \cite{Xu2024b} and the proof of Proposition 4.10 in \cite{ElEuchFukasawaRosenbaum2018}, we also have that
 \beqnn
 f^{\alpha,\gamma}*\tilde{V}(t)
 \ar=\ar V_*(0)*f^{\alpha,\gamma}(t)+\frac a b\cdot \int_0^t f^{\alpha,\gamma}*K(s)ds \cr
 \ar\ar +\int_0^t f^{\alpha,\gamma}*K(t-s)\cdot \frac{\zeta^\mathtt{m}_*\sqrt{\lambda^\mathtt{m}_*}}{b}\cdot \sqrt{V_*(s)}dB_s \cr
 \ar\ar+\int_0^t\int_0^{\infty}\int_0^{V_*(s)}
 \Big(\int_{(t-s-y)^+}^{t-s}
 \frac{\zeta^\mathtt{l}_*}{b}\cdot f^{\alpha,\gamma}*K(r)dr \Big) \widetilde{N}(ds,dy,dz),\quad t\geq 0,
 \eeqnn
 from which we conclude that
 \beqnn
 \tilde{V}(t)-f^{\alpha,\gamma}*\tilde{V}(t)
 \ar=\ar V_*(0) \cdot \big(1-F^{\alpha,\gamma}(t)\big)+\frac a b\cdot \int_0^t \big(K-f^{\alpha,\gamma}*K\big)(s)ds\cr
 \ar\ar  +\int_0^t  \big(K-f^{\alpha,\gamma}*K\big)(t-s)\cdot \frac{\zeta^\mathtt{m}_*\sqrt{\lambda^\mathtt{m}_*}}{b}\cdot\sqrt{V_*(s)}dB_s \cr
 \ar\ar+\int_0^t\int_0^{\infty}\int_0^{V_*(s)}\int_{(t-s-y)^+}^{t-s}
 \frac{\zeta^\mathtt{l}_*}{b}\cdot \big(K-f^{\alpha,\gamma}*K\big)(r)dr \widetilde{N}(ds,dy,dz),\quad t\geq 0,
 \eeqnn
 which is equal to $V_*$ because $f^{\alpha,\gamma}=K-K*f^{\alpha,\gamma}$. Thus, \eqref{eq1} holds. 
 For the converse, it can be proved similarly that every solution of \eqref{eq1} also solves \eqref{eq2}.

%% file: RegularityMaximal.tex
 \section{Regularity and maximal inequality}
 \label{Sec.Regular}
 \setcounter{equation}{0}

 In this section we provide the proof of the H\"older continuity and the maximal inequality of the weak solutions $V^*$ to the two equivalent stochastic Volterra equations (\ref{eq1}) and (\ref{eq2}).  
 
 To this end, we recall that for a real-valued stochastic process $X$ on $[0,T]$, the Kolmogorov continuity theorem states that if for some constants $p,\kappa, C>0$, 
 \beqlb \label{kol con}
 \mathbf{E} \Big[ \big|X(t)-X(s)\big|^p \Big]\leq C\cdot \big|t-s\big|^{1+\kappa},  
 \eeqlb
 uniformly in $0 \leq s,t \leq T$, then the process has a $\theta$-H\"older continuous modification for all $0<\theta<\kappa/p$. 
 Hence we first need to establish the moment estimates for the increments of $V^*$.
 For convention, we rewrite (\ref{eq1}) as 
 \[
	 V_*(t)= I_*(t) +J_{1,*}(t) +J_{2,*}(t) 
\]
where 
 \beqnn
    I_*(t)\ar:=\ar V_*(0)\cdot \big(1-F^{\alpha,\gamma}(t)\big)+\frac{a}{b}\cdot F^{\alpha,\gamma}(t),\cr
 	J_{1,*}(t)\ar:=\ar\int_0^t f^{\alpha,\gamma}(t-s)\cdot \frac{\zeta^\mathtt{m}_*\sqrt{\lambda^\mathtt{m}_*}}{b}\cdot \sqrt{V_*(s)}dB_s, \\
 	J_{2,*}(t)\ar:=\ar \int_0^t\int_0^{\infty}\int_0^{V_*(s)}\int_{(t-s-y)^+}^{t-s}
 	\frac{\zeta^\mathtt{l}_*}{b}\cdot f^{\alpha,\gamma}(r)dr \widetilde{N}(ds,dy,dz). 
 \eeqnn
 
 The $\alpha$-H\"older continuity of Mittag-Leffler distribution $F^{\alpha,\gamma}$ induces the global $\alpha$-H\"older continuity of the function $I_*$. It remains to establish the H\"older continuity of the processes $J_{1,*}$ and $J_{2,*}$. 
 For this, we need the following moment estimate for $V_*$. The proof is similar to that of Lemma~\ref{V2p} and is hence omitted. 
 
 \begin{lemma} \label{V*sup}
 	For any $p>0$,  there exists a constant $C>0$ such that for all $T\geq 0$ and $i=1,2$, 
 	\beqnn
 	\sup_{t\in[0,T]}\mathbf E \Big[\big|J_{i,*}(t)\big|^{2p} \Big ]\leq C\cdot (1+T)^{2p\alpha}
 	\quad \mbox{and}\quad 
 	\sup_{t\in[0,T]}\mathbf E \Big[\big|V_*(t)\big|^{2p} \Big ]\leq C\cdot (1+T)^{2p\alpha} . 
 	\eeqnn 
 \end{lemma}
 
  The corresponding result for the process $J_{2,*}$ can be established using the same arguments as in \cite[Theorem 2.9, Lemma 5.4]{Xu2024b} with Lemma \ref{V*sup} replacing Lemma 5.3 therein. 
 
 \begin{lemma} \label{holderJ2}
 For each $p\geq 1$, there exists a constant $C>0$ such that for all $T\geq 0$ and $t_1,t_2\in [0,T]$,
  \begin{align*}
  \mathbf{E}\Big [\big|J_{2,*}(t_2)- J_{2,*}(t_1)\big|^{2p} \Big]\leq C\cdot (1+T)^{p\alpha}\cdot \big| t_2-t_1 \big|^{p\alpha}.
  \end{align*}
 Hence $J_{2,*}$ is H\"older continuous of any order strictly less than $\alpha/2$.  
 \end{lemma}
 
 \begin{lemma} \label{holderJ1}
 For each $p\geq 1$, there exists a constant $C>0$ such that for all $T\geq 0$ and $t_1,t_2\in [0,T]$,
 \beqlb \label{eqn.503}
  \mathbf{E}\Big[ \big|J_{1,*}(t_2) - J_{1,*}(t_1)\big|^{2p} \Big]\leq C\cdot (1+T)^{p\alpha}\cdot \big| t_2-t_1 \big|^{p(2\alpha-1)}.
 \eeqlb
 Hence $J_{2,*}$ is H\"older continuous of any order strictly less than $\alpha-\frac{1}{2}$.  
 \end{lemma}
 \proof The H\"older continuity of $J_{1,*}$ follows directly from (\ref{eqn.503}) and the Kolmogorov continuity theorem. 
 To prove (\ref{eqn.503}) we use the same arguments as in the proof of Lemma~\ref{V2p} to find a constant $C>0$ that depends only on $p$ such that for any $0\leq t_1<t_2<\infty$, 
 \beqnn
  \mathbf{E}\Big[\big|J_{1,*}(t_2)-J_{1,*}(t_1)\big|^{2p}\Big] 
  \ar\leq\ar C\cdot \mathbf{E}\bigg[\Big|\int_{0}^{t_2} \big | f^{\alpha,\gamma}(t_2-s)-f^{\alpha,\gamma}(t_1-s)\big|^2 V_*(s)ds\Big|^{p}\bigg] \cr 
  \ar\leq\ar C\cdot \mathbf{E}\bigg[\Big|\int_{t_1}^{t_2} \big | f^{\alpha,\gamma}(t_2-s)\big|^2 V_*(s)ds\Big|^{p}\bigg]\cr 
  \ar\ar  + \ C\cdot \mathbf{E}\bigg[\Big|\int_{0}^{t_1} \big| f^{\alpha,\gamma}(t_2-s)-f^{\alpha,\gamma}(t_1-s)\big|^2 V_*(s)ds\Big|^{p}\bigg]. 
 \eeqnn
 By using H\"older's inequality\footnote{For two function $f,g$ and $p\geq 1$, we have 
 	\beqnn
 	\Big|\int f(s)g(s)ds\Big|^p\leq \Big|\int |f(s)|^{1-1/p}\cdot |f(s)|^{1/p}g(s)ds\Big|^p \leq  \Big|\int |f(s)|ds\Big|^{p-1}\cdot \int |f(s)| |g(s)|^p ds .
 	\eeqnn} and then Fubini's theorem along with Lemma~\ref{V*sup} to the last two expectations, we have uniformly in $0\leq t_1<t_2\leq T$, 
 \beqlb\label{eqn.501}
 \mathbf{E}\bigg[\Big|\int_{t_1}^{t_2} \big | f^{\alpha,\gamma}(t_2-s)\big|^2 V_*(s)ds\Big|^{p}\bigg]
 \ar\leq\ar \Big|\int_{t_1}^{t_2} \big | f^{\alpha,\gamma}(t_2-s)\big|^2 ds\Big|^{p-1} \cdot \int_{t_1}^{t_2} \big | f^{\alpha,\gamma}(t_2-s)\big|^2 \mathbf{E} \Big[\big|V_*(s)\big|^p\Big]ds \cr 
 \ar\leq\ar  C\cdot (1+T)^{p\alpha} \cdot \Big|\int_{0}^{t_2-t_1} \big | f^{\alpha,\gamma}(s)\big|^2 ds\Big|^{p}
 \eeqlb
 and similarly,
 \beqlb\label{eqn.502}
 \lefteqn{\mathbf{E}\bigg[\Big|\int_{0}^{t_1} \big | f^{\alpha,\gamma}(t_2-s)-f^{\alpha,\gamma}(t_1-s)\big|^2 V_*(s)ds\Big|^{p}\bigg]}\ar\ar\cr 
  \ar\leq\ar C \cdot (1+T)^{p\alpha} \cdot \Big|\int_{0}^{t_1} \big |  f^{\alpha,\gamma}(t_2-t_1+s)-f^{\alpha,\gamma}(s)\big|^2  ds\Big|^{p}. 
 \eeqlb
 By \eqref{upper bound of mittag leffler density functions}, we have $ f^{\alpha,\gamma}(s) \leq C\cdot s^{\alpha-1}$ uniformly in $s>0$ and $\int_{0}^{t_2-t_1} \big | f^{\alpha,\gamma}(s)\big|^2 ds \leq C\cdot |t_2-t_1|^{2\alpha-1}$ uniformly in $t_2\geq t_2\geq 0$. 
 Plugging this back into the last term in (\ref{eqn.501}), 
 \beqnn
 \mathbf{E}\bigg[\Big|\int_{t_1}^{t_2} \big | f^{\alpha,\gamma}(t_2-s)\big|^2 V_*(s)ds\Big|^{p}\bigg] 
 \leq  C\cdot (1+T)^{p\alpha} \cdot \big| t_2-t_1 \big|^{p(2\alpha-1)} .
 \eeqnn
 Moreover, we also have $|(f^{\alpha,\gamma})'(r)|\leq C\cdot r^{\alpha-2}$ uniformly in $r>0$ and hence 
 \beqnn
 \big |  f^{\alpha,\gamma}(t_2-t_1+s)-f^{\alpha,\gamma}(s)\big| \ar\leq\ar 
 \big( \sup_{r\in[s, t_2-t_1+s]} f^{\alpha,\gamma}(s)\big) \wedge \Big(|t_2-t_1|\cdot \sup_{r\in[s, t_2-t_1+s]} \big|(f^{\alpha,\gamma})'(r)\big|\Big) \cr
 \ar\leq\ar C\cdot \Big(s^{\alpha-1}\wedge \big(|t_2-t_1|\cdot s^{\alpha-2}\big)\Big)
 \eeqnn
 for some constant $C>0$ independent of $t_1,t_2,s$, which immediately yields that 
 \beqnn
 \lefteqn{\mathbf{E}\bigg[\Big|\int_{0}^{t_1} \big | f^{\alpha,\gamma}(t_2-s)-f^{\alpha,\gamma}(t_1-s)\big|^2 V_*(s)ds\Big|^{p}\bigg] }\ar\ar\cr
 \ar\leq\ar C \cdot (1+T)^{p\alpha} \cdot \Big|\int_{0}^{t_2-t_1} \big |  f^{\alpha,\gamma}(t_2-t_1+s)-f^{\alpha,\gamma}(s)\big|^2  ds\Big|^{p} \cr
 \ar\ar +  C \cdot (1+T)^{p\alpha} \cdot \Big|\int_{t_2-t_1}^\infty \big |  f^{\alpha,\gamma}(t_2-t_1+s)-f^{\alpha,\gamma}(s)\big|^2  ds\Big|^{p}\cr
 \ar\leq\ar C \cdot (1+T)^{p\alpha} \cdot \Big|\int_{0}^{t_2-t_1} s^{2\alpha-2}  ds\Big|^{p} \cr 
 \ar\ar  +  C \cdot (1+T)^{p\alpha} \cdot \Big||t_2-t_1|^2\cdot\int_{t_2-t_1}^\infty  s^{2\alpha-4}  ds\Big|^{p}\cr 
 \ar\leq\ar C \cdot (1+T)^{p\alpha} \cdot \big|t_2-t_1\big|^{p(2\alpha-1)}.
 \eeqnn
 The desired upper bound (\ref{eqn.503}) follows by putting all estimates above together. 
 \qed 
 
 Armed with the preceding lemmas, we are now ready to prove Theorem~\ref{MainThm.02}. 
 The proof will need the Garsia-Rodemich-Rumsey inequality; see Lemma~1.1 in \cite{GarsiaRodemichRumseyRosenblatt1970} with $\psi(u)=|u|^p$ and $p(u)=|u|^{q+1/p}$ for $q>1/p$.
 It states that for a continuous function $f$ on $\mathbb{R}_+$, there exists a constant $C_{p,q}>0$ such that for any $x_2>x_1\geq 0$, 
 \beqlb\label{eqn.405}
 \big|f(x_2)-f(x_1) \big| \leq C_{p,q} \cdot \big| x_2-x_1\big|^{pq-1} \int_{x_1}^{x_2} ds \int_{x_1}^{x_2} \frac{|f(s)-f(r)|}{|s-r|^{pq+1}}dr.
 \eeqlb
 
 \textbf{\textit{Proof of Theorem~\ref{MainThm.02}.}} 
 By Lemma~\ref{V*sup} and \ref{holderJ1}, the process $V_*$ is H\"older continuous of order strictly less than $\alpha\wedge (\alpha/2)\wedge (\alpha-1/2)= \alpha-1/2$ and claim (1) holds.
 
 For claim (2), for any $\kappa\in (0,\alpha-1/2)$, the non-randomness of $I_*$ tells that $\|I_*\|_{C^\kappa_T} <\infty$ uniformly in $T>0$. 
 Similarly as in the proof of Theorem~2.9 in \cite{Xu2024b}, we have uniformly in $T\geq 0$,
 \beqnn
 \mathbf{E}\Big[\big\|J_{2,*}\big\|_{C^\kappa_T}^p\Big] \leq  C\cdot(1+T)^{p(\alpha-\kappa)}. 
 \eeqnn
 We now consider $\|J_{1,*}\|_{C^\kappa_T} $. 
 By (\ref{eqn.405}) with $p>(\alpha-1/2-\kappa)^{-1}$ and $q=1/p+\kappa$, 
 \beqnn
 \|J_{1,*}\|_{C^\kappa_T}^p 
 \ar =\ar  \sup_{0\leq t_1,t_2\leq T} \frac{|J_{1,*}(t_2)-J_{1,*}(t_1)|^p}{|t_2-t_1|^{p\kappa}} 
 \leq C\cdot \int_{0}^T ds \int_0^T \frac{|J_{1,*}(s)-J_{1,*}(r)|^p}{|s-r|^{p\kappa+2}}dr.
 \eeqnn
 Taking expectations on both sides of the preceding inequality and then using Fubini's theorem along with Lemma~\ref{holderJ1}, 
 \beqnn
 \mathbf{E}\Big[  \|J_{1,*}\|_{C^\kappa_T}^p \Big] \leq  C\cdot (1+T)^{p\alpha/2} \int_0^T ds \int_0^T |s-r|^{p(\alpha-1/2)-p\kappa-2}dr 
 \leq C\cdot  (1+T)^{p(3\alpha/2-1/2-\kappa)}.
 \eeqnn
 Putting all these estimates together, notice that $\alpha-\kappa>3\alpha/2-1/2-\kappa$ we have uniformly in $T\geq 0$, 
 \beqnn
  \mathbf{E}\Big[  \|V_{*}\|_{C^\kappa_T}^p \Big] \leq C\cdot \Big(\|I_*\|_{C^\kappa_T}^p + \mathbf{E}\Big[  \|J_{1,*}\|_{C^\kappa_T}^p \Big]+\mathbf{E}\Big[  \|J_{2,*}\|_{C^\kappa_T}^p \Big] \Big) \leq C\cdot (1+ T)^{p(\alpha-\kappa)}
 \eeqnn
 and hence claim (2) holds. 
 
 Finally, we prove claim (3). For each $\kappa\in (0,\alpha-1/2)$, 
 \beqnn
 \sup_{t\in[0,T]}\big| V_*(t) \big|^p
 \ar\leq \ar C\cdot \sup_{t\in[0,T]}\big| V_*(t) -V_*(0)\big|^p + C\cdot \big|V_*(0)\big|^p
 \leq C\cdot  \|V_{*}\|_{C^\kappa_T}^p \cdot T^{p\kappa} +  C\cdot \big|V_*(0)\big|^p.
 \eeqnn
 Taking expectations on both sides of these inequality and then using claim (2), 
 \beqnn
 \mathbf{E}\bigg[  \sup_{t\in[0,T]}\big| V_*(t) \big|^p \bigg] \leq  C\cdot \mathbf{E}\Big[ \|V_{*}\|_{C^\kappa_T}^p\Big] \cdot T^{p\kappa} +  C\cdot \big|V_*(0)\big|^p
 \leq C (1+T)^{p\alpha},
 \eeqnn
 for some constant $C>0$ independent of $T$ and claim (3) holds.
 \qed

%% file: NolinearVolterraRiccatiNew.tex
 \section{Laplace functionals and weak uniqueness}
 \label{Sec.VolterraRiccati}
 \setcounter{equation}{0}
 
 So far we have proved the $C$-tightness of the sequence of rescaled processes $\{V^{(n)}\}_{n\geq 1}$ and identified the stochastic equation \eqref{eq1} that each weak accumulation point $V_*$ must satisfy.  To establish the weak uniqueness of accumulation points, we need to determine the characteristic functional of $V_*$. In this section we prove that for any given $\lambda\geq 0$ and any $g\in L^\infty(\mathbb{R}_+;\mathbb{R}_+)$, 
 \begin{align}\label{charF}
 	\textbf{E}\Big[\exp\big\{-\lambda\cdot V_*(T)-g*V_*(T)\big\}\Big]=
 	\exp \Big\{-V_*(0)\cdot L_K*\psi_g^\lambda(T)-\frac{a}{b}*\psi_g^\lambda(T) \Big\},\quad T> 0,
 \end{align}
 where the function $\psi_g^\lambda$ satisfies the nonlinear Volterra equation \eqref{PSI}.

  \subsection{Well-posedness of generalized Volterra-Riccati equation}

 In view of the second equality in \eqref{LK} the Volterra equation \eqref{PSI} can be brought into the more convenient form 
 \begin{align}\label{psi2}
 	\psi_g^\lambda=\lambda\cdot  K+g*K
 	-
 	\frac 1 2 \cdot \Big(\frac{\zeta_*^\mathtt{m}\sqrt{\lambda_*^\mathtt m }}{b}\Big)^2\cdot |\psi_g^\lambda|^2*K-  \big(\mathcal{V}\circ \psi_g^\lambda\big)*K-\psi_g^\lambda*K.
 \end{align}
 
 \begin{remark}
Convoluting both sides of the above equation with the function $L_K$ and then using the first equality in \eqref{LK} we obtain that
 \beqnn
     \psi_g^\lambda*L_K=\lambda +g*1 -  \frac 1 2 \cdot \Big(\frac{\zeta_*^\mathtt{m}\sqrt{\lambda_*^\mathtt m }}{b}\Big)^2\cdot |\psi_g^\lambda|^2*1-  \big(\mathcal{V}\circ \psi_g^\lambda\big)*1-\psi_g^\lambda*1.
 \eeqnn
 Differentiating both sides of this equation and then recalling the definitions of Riemann-Liouville fractional derivative operator and integral operator, we see that
the equation \eqref{PSI} and hence the equation \eqref{psi2} can be rewritten as  
 \beqnn
 \frac{1}{\gamma} \cdot D^\alpha\psi_g^\lambda = g  -  \frac 1 2 \cdot \Big(\frac{\zeta_*^\mathtt{m}\sqrt{\lambda_*^\mathtt m }}{b}\Big)^2\cdot |\psi_g^\lambda|^2 -  \big(\mathcal{V}\circ \psi_g^\lambda\big) -\psi_g^\lambda 
 \quad \mbox{with}\quad
 I^{1-\alpha} \psi_g^\lambda(0)= \gamma\cdot \lambda, 
 \eeqnn 
 where $D^\alpha$ and $D^\alpha$ denote the \textsl{Riemann-Liouville fractional derivative operator} and \textsl{integral operator} that acts on a function $f:\mathbb R\rightarrow \mathbb R$ according to
 \beqnn
 D^\alpha f(x):= \frac{1}{\Gamma(1-\alpha)}\frac{d}{dx}\int_0^x \frac{f(t)}{(x-t)^\alpha}dt
 \quad {and}\quad 
 I^\alpha f(x) :=   \frac{1}{\Gamma(\alpha)} \int_0^x \frac{f(t)}{(x-t)^{1-\alpha}}dt
 \eeqnn
 \end{remark}

 Any solution to the nonlinear Volterra equation \eqref{PSI} behaves like the function $\lambda \cdot f^{\alpha,\gamma}$ in a vicinity of the origin. Since 
  \beqnn
 	f^{\alpha,\gamma}(t)\sim C\cdot t^{\alpha-1} \quad \mbox{as} \quad t\to 0, 
 \eeqnn	
we expect to find solutions $\psi^\lambda_g$ to the equation \eqref{psi2} on the interval $[0,T]$ in the space 
 \beqnn
  \mathcal{A}_T^\lambda := \Big\{ f \in L^1\big((0,T],\mathbb{R} \big) : \| f \|_{L^\infty_{T,\alpha}} <\infty \mbox{ if }\lambda >0 \mbox{ and }  \| f \|_{L^\infty_{T}} <\infty \mbox{ if }\lambda =0 \Big\}, 
 \eeqnn
 where  these two norms are defined by
 \beqnn
  \| f \|_{L^\infty_{T,\alpha}} := \sup_{t \in (0,T]}  t^{1-\alpha}|f(t)|
  \quad \mbox{and}\quad 
  \| f \|_{L^\infty_{T}} := \sup_{t \in (0,T]}  |f(t)|. 
 \eeqnn

 \begin{definition}
 \begin{itemize}
	\item[(i)]  A pair $(T,\psi_g^\lambda)\in (0,\infty)\times \mathcal{A}_T^\lambda$ is called a {\rm local solution} to the equation \eqref{PSI} if the function $\psi_g^\lambda$ satisfies the equation \eqref{PSI} on $(0,T]$. 
	\item[(ii)] A pair $(T_g^\lambda,\psi_g^\lambda)$ is called a {\rm noncontinuable solution} to the equation \eqref{PSI} if $(T,\psi_g^\lambda)\in (0,T_g^\lambda)\times \mathcal{A}_T^\lambda$ is a local solution and if 
  \beqnn
  \lim_{T \uparrow T^\lambda_g}\big\|\psi_{g}^\lambda\big\|_{L^\infty_{T,\alpha}}=\infty\quad\mbox{whenever}\quad T_g^\lambda<\infty.
  \eeqnn 
 The function $\psi_{\lambda,g}$ is called a {\rm global solution} if $T_g^\lambda=\infty$. 
 \end{itemize}
 \end{definition}

 It follows from Lemma \ref{psi nonnegative} below that noncontinuable solutions are global solutions. To prove the existence of noncontinuable solutions we shall first establish the existence and uniqueness of local solutions for which we introduce, for any $T>0$ and $M>0$ the sets 
 \beqlb\label{eqn.100}
 	\mathcal{A}_{T,M}:=\big\{ f\in L^1\big((0,T],\mathbb{R} \big):\ 	\| f \|_{L^\infty_{T,\alpha}}\leq M \big\}  \quad \mbox{and} \quad 
	\mathcal{B}_{T,M}:=\big\{ f\in L^\infty([0,T];\mathbb{R}): \|f\|_{L_T^\infty}\leq M \big\}. 
 \eeqlb

 
 In what follows we shall repeatedly use the following auxiliary results. The first provides important estimates for the operator $\mathcal{V}$ acting on these sets; for a proof we refer to Propositions 6.2-6.5 in \cite{Xu2024b}.  The second is a direct consequence of Theorem 2.1 in \cite{DentonVatsala2010} and Theorem~2.7 in \cite{VasundharaDeviMcRaeDrici2012}.

 \begin{proposition}
 	$(i)$ There exists a constant $C>0$, such that for any $M,T>0$, $f\in\mathcal{A}_{T,M}$ and $t\in(0,T]$,
 	\begin{gather}
 		\big|(\mathcal{V}\circ f)*f^{\alpha,\gamma}(t)\big|\leq CM^2\cdot e^{\frac{M}{\alpha}t^\alpha}\cdot t^{2\alpha-1} \label{6.1}.
 	\end{gather}
 	Moreover, for any $\rho\in(1,(1-\alpha)^{-1})$, there exists a constant $C>0$ such that for any $M,T>0$ and $f_1,f_2\in\mathcal{A}_{T,M}$,
 	\begin{align}
 		\big\|(\mathcal{V}\circ f_1-\mathcal{V}\circ f_2)*f^{\alpha,\gamma}\big\|_{L_T^\rho}\leq CM\cdot T^\alpha\cdot e^{\frac{M}{\alpha}T^\alpha}\cdot \big\|f_1-f_2\big\|_{L_T^\rho}.
 		\label{6.2}
 	\end{align}
 	$(ii)$ There exists a constant $C>0$ such that for any $M,T>0$, $f\in\mathcal{B}_{T,M}$ and $t\in(0,T]$,
 	\begin{align}
 		\big|(\mathcal{V}\circ f)*f^{\alpha,\gamma}(t)\big|\leq C\cdot M^2 e^{Mt}\cdot t.
 		\label{6.3}
 	\end{align}
 	Furthermore, there exists a constant $C>0$ such that for any $T,M>0$ and $f_1,f_2\in\mathcal{B}_{T,M}$,
 	\begin{align}
 		\big\|(\mathcal{V}\circ f_1-\mathcal{V}\circ f_2)*f^{\alpha,\gamma}\big\|_{L_T^1}\leq C\cdot MT e^{M T}\cdot \big\|f_1-f_2\big\|_{L_T^1}.
 		\label{6.4}
 	\end{align}
 \end{proposition}

%
%
%
 
 \begin{proposition}\label{Prop.Comparison}
  For $T>0$, let $w,v$ be continuous functions on $(0,T]$ that satisfy  
 \beqnn
 w(t)\geq \hat{w}_0(t)+  \big( |w|^2-w \big)*K(t) 
 \quad \mbox{and}\quad 
 v(t)\leq \hat{v}_0(t)+  \big( |v|^2-v \big)*K(t),\quad t\in(0,T],
 \eeqnn 
where $K$ denotes the function introduced in \eqref{func:K}. Then $w(t)>v(t)$ for all $t\in(0,T]$ if one of the above inequalities is strict and one of the following two conditions holds:  
 \begin{enumerate}
 	\item[(i)] $\hat{w}_0(t)=w_0 \cdot K(t)$ and $\hat{v}_0(t)=v_0 \cdot K(t)$ for two constants $w_0,v_0 \in \mathbb{R} $ such that $w_0>v_0$. 
 	
 	\item[(ii)] $\hat{w}_0(t)\equiv w_0  $ and $\hat{v}_0(t)\equiv v_0 $ for two constants $w_0,v_0 \in \mathbb{R} $ such that $w_0\geq v_0$. 
 \end{enumerate}
 \end{proposition}
 
 The next proposition establishes upper and lower bounds on the functions $ \mathcal{V}\circ f$ for $f \in \mathcal{A}_T$.
 
 \begin{proposition}
 	For any $T>0$ and $f \in \mathcal{A}_T^\lambda$, there exists a constant $C>0$ such that for all $t\in[0,T]$, 
 	\beqlb\label{eqn.601}
  \mathcal{V}\circ f(t) \geq 0
  \quad \mbox{and}\quad 
  	0\leq (\mathcal{V}\circ f)*K(t)\leq
 	C \cdot |f|^2*K(t).
 	\eeqlb
 \end{proposition}
 \proof Using the inequality $e^{-z}-1+z\geq 0$ for all $z\in \mathbb{R}$, we see that $\mathcal{V}\circ f(t)\geq 0$ and hence that $(\mathcal{V}\circ f)*K(t)\geq 0$ for all $t\in[0,T]$. 
 
 To prove the upper bound we use the inequality $e^{-z}-1+z\leq |z|^2e^{|z|}$ for any $z\in\mathbb R$. Since $f \in L^1\big((0,T],\mathbb{R} \big)$,  
 \beqlb\label{Vpsi1} 
  \mathcal{V}\circ f(t)
  \ar\leq\ar \exp\bigg\{\Big|\int_0^t\frac{\zeta_*^\mathtt{l}}{b}\cdot f(s)ds\Big| \bigg\}\cdot\int_0^\infty \Big|\int_{(t-y)^+}^t\frac{\zeta_*^\mathtt{l}}{b}\cdot f(s)ds\Big|^2\nu(dy)\cr
  \ar\leq\ar C \cdot \int_0^\infty \Big|\int_{(t-y)^+}^t f(s)ds\Big|^2\nu(dy),
 \eeqlb
 uniformly in for any $t\in (0,T]$. Moreover, by H\"{o}lder's inequality,
 \begin{align*}
  \Big|\int_{(t-y)^+}^tf(s)ds\Big|^2\leq \int_{(t-y)^+}^t|f(s)|^2ds \cdot (t\wedge y), \quad t \in [0,T].
 \end{align*}
 Taking this back into \eqref{Vpsi1} and recalling the function $L_K$ defined in \eqref{func:K} yields that
 \beqnn
 \mathcal{V}\circ f(t) 
 \ar\leq\ar C\cdot \int_0^\infty \int_{(t-y)^+}^t|f(s)|^2ds \cdot (t\wedge y)\frac{dy}{y^{\alpha+2}}\cr
 \ar=\ar C\cdot \int_0^t \int_{t-y}^t|f(s)|^2ds  \frac{dy}{y^{\alpha+1}}  + C\cdot t\cdot \int_t^\infty \int_{0}^t|f(s)|^2ds  \frac{dy}{y^{\alpha+2}}\cr
 \ar\leq\ar C\cdot \int_0^t \int_{t-y}^t|f(s)|^2ds  \frac{dy}{y^{\alpha+1}}  +     C\cdot L_K(t)\int_{0}^t|f(s)|^2ds . 
 \eeqnn
 Applying Fubini's theorem to the first double integral on the right side of the last equality shows that 
 \beqnn
 \int_0^t \int_{t-y}^t|f(s)|^2ds  \frac{dy}{y^{\alpha+1}}
 \ar=\ar \int_0^t|f(s)|^2ds \int_{t-s}^t  \frac{dy}{y^{\alpha+1}}
 \leq  C\cdot L_K*|f|^2(t).
 \eeqnn
 As a result, there exists a constant $C>0$ such that for all $t\in  [0,T]$,
 \beqnn
  \mathcal{V}\circ f(t) \leq C\cdot L_K*|f|^2(t) +  C\cdot L_K(t)\int_{0}^t|f(s)|^2ds,
 \eeqnn
 which, along with the identity $K*L_K= L_K*K\equiv 1$ (see \eqref{LK}) shows that 
 \beqnn
 \big(\mathcal{V}\circ f\big)*K(t) \leq C\cdot K* L_K*|f|^2(t) +  C\cdot K* L_K(t) \cdot \int_{0}^t|f(s)|^2ds 
 \leq C\cdot \int_{0}^t|f(s)|^2ds .
 \eeqnn
 Moreover, by the monotonicity of $K$ we have 
 \beqnn
 \int_{0}^t|f(s)|^2ds \leq \int_{0}^t \frac{K(t-s)}{K(t)}|f(s)|^2ds \leq \frac{K*|f|^2(t)}{K(T)}, 
 \eeqnn
 uniformly in $t\in[0,T]$ and then the second desired upper bound in \eqref{eqn.601} holds.
 \qed 
 
 Next, we show that if a noncontinuable solution to our Volterra equations exists, then a continuous global solution exists.  
 
 \begin{lemma}\label{psi nonnegative}
 If  the equation \eqref{PSI} admits a noncontinuable solution $(T^\lambda_g, \psi^\lambda_g)$, then the following hold.
 \begin{enumerate}
  \item[(1)] For any $T < T^\lambda_g$, there exists a constant $C>0$ such that 
 \beqnn
	 0\leq \psi_g^\lambda(t) \leq C\cdot\big(1+ t^{\alpha-1}\big) \quad \mbox{on} \quad (0,T]. 
 \eeqnn
 
 \item[(2)] There exists a continuous global solution to the equation \eqref{PSI}.
 	
 \end{enumerate} 
 \end{lemma}
 \proof We prove these two claims separately. 
 \begin{itemize}
 
\item[(1)]  
Applying the first inequality in \eqref{eqn.601} to \eqref{PSI} and then using the estimate \eqref{upper bound of mittag leffler density functions} yields that 
 \beqnn
 \psi^\lambda_{g}  (t) \leq \lambda\cdot f^{\alpha,\gamma}(t)+g*f^{\alpha,\gamma}(t) 
 \leq C \cdot\lambda\cdot t^{\alpha-1} + \|g\|_{L^\infty} \leq C(1 +t^{\alpha-1})
 \eeqnn
 on any set $(0,T]$ where the function is defined.  
This shows the upper bound. To prove the lower bound, let $\widetilde\psi_g^\lambda:=-\psi_g^\lambda$.
 Multiplying both sides of the equation \eqref{psi2} by $-1$ gives that
 \beqnn
 \widetilde\psi_g^\lambda
 \ar=\ar -\lambda\cdot  K-g*K
 +
 \frac 1 2 \cdot \Big(\frac{\zeta_*^\mathtt{m}\sqrt{\lambda_*^\mathtt m }}{b}\Big)^2\cdot | \widetilde \psi_g^\lambda|^2*K+ \big(\mathcal{V}\circ \psi_g^\lambda\big)*K-\widetilde \psi_g^\lambda*K.
 \eeqnn
 
 We now distinguish two case.

\begin{itemize} 
	\item If $\lambda>0$,  we have $\widetilde \psi^\lambda_g \in \mathcal{A}_{T,M}$ for some $M>0$.  Since $K(t)\sim  t^{\alpha-1}$ it follows from the two estimates \eqref{upper bound of mittag leffler density functions} and \eqref{6.1} that for small $t>0$,
	\[
	g * K(t) \sim t^{\alpha}, \quad |\widetilde \psi^\lambda_g|^2 * K(t) \sim t^{3\alpha-2},
	\quad  |\widetilde \psi_g^\lambda*K(t)| \sim t^{2\alpha-1}, \quad 
	| \big(\mathcal{V}\circ \psi_g^\lambda\big)*K (t)| \leq C \cdot t^{2\alpha - 1}.	
	\]  
 These immediately yield that
 \[
	 t^{1-\alpha} \widetilde\psi_g^\lambda(t) \to -\frac{\gamma\cdot \lambda}{\Gamma(\alpha)} <0 \quad \mbox{as} \quad t \to 0+.
 \]
Moreover, follows from the second inequality in \eqref{eqn.601} that there exists a constant $C>0$ such that for any $t\in(0,T]$,
 \begin{align*}
 \widetilde\psi_g^\lambda(t)
& ~ \leq  -\lambda\cdot  K(t) -g*K(t) + C \cdot \Big( \big|\widetilde\psi_g^\lambda\big|^2- \widetilde\psi_g^\lambda  \Big)*K(t) \\
 & ~ < -\frac{\lambda}{2}\cdot  K(t)  + C \cdot \Big( \big|\widetilde\psi_g^\lambda\big|^2- \widetilde\psi_g^\lambda  \Big)*K(t),\quad t\in[0,T].
 \end{align*}
Applying Proposition~\ref{Prop.Comparison}(i) to $v = \widetilde \psi^\lambda_g$ and $w=0$ then yields that $\widetilde\psi_g^\lambda(t) \leq 0$. 

 \item If $\lambda=0$, we have $\widetilde\psi_g^\lambda \in \mathcal{B}_{T,M}$ for some $M>0$. 
 Similarly as in the preceding case, by using the estimates \eqref{upper bound of mittag leffler density functions} and \eqref{6.3} we also have for small $t>0$,
 \[
 g * K(t) \sim t^{\alpha}, \quad |\tilde \psi^\lambda_g|^2 * K(t) \sim t^{\alpha},
 \quad  |\widetilde \psi_g^\lambda*K(t)| \sim t^{\alpha}, \quad 
 | \big(\mathcal{V}\circ \psi_g^\lambda\big)*K (t)| \leq C .	
 \]  
 Hence we have 
 \[
 t^{1-\alpha} \widetilde\psi_g^\lambda(t) \to 0 \quad \mbox{as} \quad t \to 0+.
 \]
  and then
 \beqnn
 \widetilde\psi_g^\lambda(t)
 \ar<\ar   C \cdot \Big( \big|\widetilde\psi_g^\lambda\big|^2- \widetilde\psi_g^\lambda  \Big)*K(t),\quad t\in[0,T].
 \eeqnn
 In this case, it thus follows from Proposition~\ref{Prop.Comparison}(ii) that $\widetilde\psi_g^\lambda(t)\leq 0$.
%
\end{itemize}
 
 	\item[(2)] In view of the estimate established in (1) any noncontinuable solution is a global solution and that the solution is locally square integrable. The continuity of the solution hence follows from the fact that for any two functions $f\in L^p(\mathbb{R}_+;\mathbb{R})$ and $h\in L^q(\mathbb{R}_+;\mathbb{R})$ with $p,q \in[1,\infty]$ and $1/p+1/q=0$, the convolution $f*h$ is continuous on $\mathbb{R}_+$. In particular, the right hand side of \eqref{PSI} is continuous. 
  \end{itemize}
  \qed
 
 With the preceding preparations in hand, we are now ready to prove our existence and uniqueness of solutions result for the Volterra equation \eqref{PSI}. 
We recall the two sets $\mathcal{A}_{T,M}$ and $\mathcal{B}_{T,M}$ defined in \eqref{eqn.100}, which are closed, bounded and convex subsets of $L^\rho([0,T];\mathbb{R})$ for all $\rho\in(1,(1-\alpha)^{-1})$ and of $L^1([0,T];\mathbb{R})$, respectively.

 \medskip
 \textbf{\textit{Proof of Theorem~\ref{thmcf}(1).}} 
 It suffices to prove the existence  of a noncontinuable solution and the uniqueness of global solutions. 
 Here we just deal with the case of $\lambda >0$. 
 The case of $\lambda=0$ can be proved in the same way by replacing $\mathcal{A}_{T,M}$  with $\mathcal{B}_{T,M}$.
 We proceed in three steps, starting with the existence of local solutions. 
 
 \medskip
 
 {\bf Step 1: Existence of local solutions}.
 We will prove that for some $M_0>0$ and $\delta>0$ that will be specified later, there exists a local solution $\big(\delta,\psi_g^\lambda\big)$ to \eqref{PSI} with  $\psi_g^\lambda \in\mathcal A_{\delta,M_0}$. 
 For this, we introduce the mapping $\mathcal{R}_0$ that acts on functions $f\in\mathcal{A}_{\delta,M_0}$ according to
 \begin{align*}
 	\mathcal{R}_0\circ f(t):= 
 	\lambda\cdot f^{\alpha,\gamma}(t)+g*f^{\alpha,\gamma}(t)- \frac{1}{2} \cdot \Big(\frac{\zeta_*^\mathtt{m}\sqrt{\lambda_*^\mathtt m }}{b}\Big)^2\cdot |f|^2*f^{\alpha,\gamma}(t)-  \big(\mathcal{V}\circ f\big)*f^{\alpha,\gamma}(t),\quad t> 0.
 \end{align*}
 Since the set $\mathcal{A}_{T,M}$ is closed, bounded and convex subsets of the spaces $L^\rho([0,T];\mathbb{R})$ for all $\rho\in(1,(1-\alpha)^{-1})$, the existence of a local solution follows from Banach's fixed point theorem if we can prove that $\mathcal{R}_0$ is a contraction w.r.t.~a suitable $L^\rho$-norm.    
 
 To this end, we first provide a pointwise estimate for the function $\mathcal{R}_0\circ f$. 
 The first two summands can be bounded using the estimate \eqref{upper bound of mittag leffler density functions}. In view of this estimate it holds uniformly in $t>0$ that
 \beqnn
 \lambda\cdot f^{\alpha,\gamma}(t)+g*f^{\alpha,\gamma} (t)
 \leq C\cdot \big(1+ t^{\alpha-1}  \big).
 \eeqnn 
 Moreover, there exists a constant $C>0$ that is independent of $M_0$ and $\delta$ such that for any $f\in\mathcal{A}_{\delta,M_0}$ and $t\in(0,\delta]$,    
 \begin{align*}
 	|f|^2*f^{\alpha,\gamma}(t) \leq  C \cdot M_0^2 \int_0^t (t-s)^{2\alpha-2}s^{\alpha-1}ds ~ \leq ~ C
 	M_0^2\cdot t^{3\alpha-2}\leq
 	CM_0^2\delta^{2\alpha-1}\cdot t^{\alpha-1}.
 \end{align*}
 Together with the inequality \eqref{6.1} this shows that there exists a constant $C_0>0$ that is independent of $M_0$ and $\delta$ such that uniformly in $f\in\mathcal{A}_{\delta,M_0}$ and $t\in(0,\delta]$, 
 \beqlb
 |\mathcal{R}_0\circ f(t)|
 \leq C_0\Big(1+M_0^2\delta^{2\alpha-1}+M_0^2e^{M_0\delta^\alpha/\alpha}\cdot \delta^{\alpha} \Big)\cdot  t^{\alpha-1}.
 \label{R0f}
 \eeqlb
 
 Furthermore, it follows from the inequality \eqref{6.2} that for $f_1,f_2\in\mathcal{A}_{\delta,M_0}$, there exists a constant $C>0$ that does not depend on $\delta$ and $M_0$, such that
 \begin{equation}\label{R0}
 	\begin{split}
 	 \big\|\mathcal{R}_0\circ f_1-\mathcal{R}_0\circ f_2\big\|_{L^\rho_\delta} & ~ \leq C\cdot\Big(
 		\big\|(f_1^2-f_2^2)*f^{\alpha,\gamma}\big\|_{L^\rho_\delta} 
 		+\big\|(\mathcal{V}\circ f_1-\mathcal{V}\circ  f_2)*f^{\alpha,\gamma} \big\|_{L^\rho_\delta}\Big) \\
		& \leq ~ C\cdot\Big( \big\|(f_1^2-f_2^2)*f^{\alpha,\gamma}\big\|_{L^\rho_\delta}  + \delta^\alpha\cdot e^{\frac{M_0}{\alpha}\delta^\alpha} \big\|f_1-f_2 \big\|_{L^\rho_\delta}  \Big) . 
 	\end{split}
 \end{equation}
 To further bound the above sum we now choose three constants $1\leq p,q,r\leq \infty$ that satisfy 
 \beqnn
	p \leq r\leq \infty, \quad \frac{1}{p}+\frac{1}{q}=\frac{1}{\rho}+1 \quad \mbox{and} \quad 
	\frac{1}{r}+\frac{1}{\rho}=\frac{1}{p}.
 \eeqnn 
 By first applying Young's convolution inequality\footnote{Young's convolution inequality implies that $\|f*g\|_{L^\rho_\delta}\leq \|f\|_{L^p_\delta}\cdot\|g\|_{L^q_\delta}$, where $\frac{1}{p}+\frac{1}{q}=\frac{1}{\rho}+1$ and $1\leq p,q\leq \infty$.} and then H\"{o}lder's inequality, we obtain that  
 \beqnn
 \big\|(f_1^2-f_2^2)*f^{\alpha,\gamma}\big\|_{L^\rho_\delta}
 \leq \big\|f^{\alpha,\gamma}\big\|_{L^q_\delta} \cdot \big\|f_1^2-f_2^2\big\|_{L^p_\delta}
 \leq \big\|f^{\alpha,\gamma}\big\|_{L^q_\delta}\cdot
 \big\|f_1+f_2\big\|_{L^r_\delta}\cdot \big\|f_1-f_2\big\|_{L^{\rho}_\delta}.
 \eeqnn  
 In particular, since $f_1,f_2\in\mathcal{A}_{\delta,M_0}$,  we may choose $r=q=2$ to get that
 \begin{align*}
  \big\|f^{\alpha,\gamma}\big\|_{L^q_\delta}\cdot
 \big\|f_1+f_2\big\|_{L^r_\delta}\leq
 C\cdot M_0\cdot  \int_0^\delta s^{2(\alpha-1)}ds  \leq C\cdot M_0\cdot \delta^{2\alpha-1}.
 \end{align*}
 Thus,
 \begin{align*}
 \big\|(f_1^2-f_2^2)*f^{\alpha,\gamma}\big\|_{L^\rho_\delta}
 \leq C\cdot M_0\cdot \delta^{2\alpha-1}\cdot \big\|f_1-f_2\big\|_{L^{\rho}_\delta}.
 \end{align*}
In particular, there exists a constant $C>0$ that is independent of $M_0$ and $\delta$ such that 
 \begin{align}
 \big\|\mathcal{R}_0\circ f_1-\mathcal{R}_0\circ f_2\big\|_{L^\rho_\delta} \leq  C \cdot  M_0 \cdot\Big( e^{M_0\delta^\alpha/\alpha}\cdot \delta^{\alpha}
 	+ \delta^{2\alpha-1} \Big) \cdot \big\|f_1-f_2\big\|_{L^{\rho}_\delta}.
 	\label{R0contractive}
 \end{align}
 
 We can thus choose $M_0$ large, and $\delta$ small enough in \eqref{R0f} and \eqref{R0contractive} such that for all $\ f, f_1,f_2\in\mathcal{A}_{\delta,M_0}$ and $t\in(0,\delta]$,
 \beqnn
 |\mathcal{R}_0\circ f(t)|\leq M_0t^{\alpha-1}\ \  \text{and}\ \  \big\|\mathcal{R}_0\circ f_1-\mathcal{R}_0\circ f_2\big\|_{L^\rho_\delta} <\big\|f_1-f_2\big\|_{L^{\rho}_\delta}.
 \eeqnn
 
 This shows that the operator $\mathcal{R}_0$ maps the set $\mathcal{A}_{\delta,M_0}$ into itself and that it is a contraction in the $L^\rho$-sense. 
 In particular, for any two fixed points $\psi_{\delta,i} \in \mathcal{A}_{\delta,M_0}$, $i=1,2$ it holds that
 \[
 	\int_{0}^\delta | \psi_{\delta,1}(t) - \psi_{\delta,2}(t) | dt = 0.  
 \]
 As a result, $\mathcal{R}_0$ admits an a.e.~uniquely defined fixed point $\psi_\delta\in \mathcal{A}_{\delta,M_0}$. In particular, there exists a local solution $(\delta,\psi_\delta)$ of the equation \eqref{PSI} with the solution being unique up to a set of Lebeque-measure zero.

 	\medskip
 
 \noindent {\bf Step 2: Existence of global solutions.} We are now going to show that any local solution can be extended to a local solution on a larger interval. 
%
%

To this end, let $(t_0,\psi_0)$ be a local solution of the equation \eqref{PSI} with $\psi_0 \in \mathcal{A}_{t_0}$.   
%
%
 To extend this solution to a larger interval we introduce the bounded function $H_1: \mathbb R_+ \to \mathbb R$ defined by
 \beqnn
 	H_1(t) :=\lambda f^{\alpha,\gamma}(t_0+t)+
 	g* f^{\alpha,\gamma}(t_0+t) -\int_0^{t_0}
 	\Big( \frac 1 2\Big(\frac{\zeta^*_\mathtt m \sqrt{\lambda_\mathtt m^*}}{b}\Big)^2|\psi_0(s)|^2+, 
 	\mathcal{V}\circ \psi_0(s) \Big) f^{\alpha,\gamma}(t_0+t-s)ds
 \eeqnn
and fix two constants $\delta_1\in(0,1)$ and $M_1>0$ to be specified later. Furthermore, we introduce the operator $\mathcal{R}_1$ that acts on the set $\mathcal B_{\delta_1,M_1}$ according to 
 \begin{align*}
 	\mathcal{R}_1\circ f
 	(t):=H_1(t)-\int_0^{t}
 	\Big[  \frac{1}{2} \Big(\frac{\zeta^*_\mathtt m \sqrt{\lambda_\mathtt m^*}}{b}\Big)^2|f(s)|^2+ \mathcal{V}\circ f(s) \Big] 
    f^{\alpha,\gamma}(t-s)ds, \quad t\in [0,\delta_1]. 
 \end{align*}
%

It follows from the inequalities (\ref{6.3}) and \eqref{6.4} that for any $t\in(0,\delta_1]$ and all functions $f,f_1,f_2\in \mathcal{B}_{\delta_1,M_1}$, there exists a constant $C_1<\infty$, such that
 \begin{align*}
 |\mathcal{R}_1\circ f(t)|
 &\leq C_1 \cdot \Big (1+ M_1^2 \cdot \delta_1^{\alpha}+M_1^2e^{M_1 \delta_1}\cdot \delta_1\Big), \quad t\in [0,\delta_1]
 \end{align*}
and
 \begin{align*}
 \big\|\mathcal{R}_1\circ f_1-\mathcal{R}_1\circ f_2\big\|_{L^1_{\delta_1}}\leq
 	C_1\cdot\big(M_1\delta_1^\alpha+ e^{M_1 t_1}\cdot M_1\delta_1\big)\cdot\big\|f_1-f_2\big\|_{L^1_{\delta_1}}.
 \end{align*}
 Choosing first $M_1$ large and then $\delta_1$ small enough, we see that 
 \beqnn
 \big\|\mathcal{R}_1\circ f\big\|_{L^\infty_{\delta_1}} \leq M_1 
 \quad \mbox{and}\quad 
 \big\|\mathcal{R}_1\circ f_1-\mathcal{R}_1\circ f_2\big\|_{L^1_{\delta_1}}<\big\|f_1-f_2\big\|_{L^1_{\delta_1}}.
 \eeqnn
 
 Following the same arguments as in Step 1, we conclude that the operator $\mathcal{R}_1$ admits an a.e.~uniquely defined fixed point $\psi_1$ in the set $\mathcal{B}_{\delta_1,M_1}$, viewed as a closed, convex, bounded subset of $L^1([0,T];\mathbb{R})$.   Let us now set 
 \beqnn
 \psi_{g}^\lambda(t):=\psi_0(t)\cdot \mathbf{1}_{(0,t_0]}(t) + \psi_1(t{\blue -t_0})\cdot \mathbf{1}_{(t_0,{\blue t_0+}\delta_1]}(t), \quad t\in(0,t_0+\delta_1] . 
 \eeqnn 
 
 The function $\psi^\lambda_g$ is a local solution to the equation \eqref{PSI} on the interval $(0,t_0+\delta_1]$. In fact, $\psi^\lambda_g(t) = \psi_0(t)$ on $(0,t_0]$ while for $t\in (t_0,t_0+\delta_1]$,  
 \begin{align*}
 \psi^\lambda_g (t) & ~= \psi_1 (t-t_0) = \mathcal R_1\circ \psi_1(t-t_0) \\
 & ~= \lambda f^{\alpha,\gamma}(t)+ g* f^{\alpha,\gamma}(t) -\int_0^{t_0} \Big( \frac 1 2\Big(\frac{\zeta^*_\mathtt m \sqrt{\lambda_\mathtt m^*}}{b}\Big)^2|\psi_0(s)|^2+ 
 \mathcal{V}\circ \psi_0(s) \Big) f^{\alpha,\gamma}(t-s)ds \\
 &\qquad- \int_0^{t-t_0}
 \Big[  \frac{1}{2} \Big(\frac{\zeta^*_\mathtt m \sqrt{\lambda_\mathtt m^*}}{b}\Big)^2|\psi_1(s-t_0)|^2+ \mathcal{V}\circ \psi_1(s-t_0) \Big]*f^{\alpha,\gamma}(t-t_0-s)ds\\
 &~= \lambda f^{\alpha,\gamma}(t)+ g* f^{\alpha,\gamma}(t) -\int_0^t \Big( \frac 1 2\Big(\frac{\zeta^*_\mathtt m \sqrt{\lambda_\mathtt m^*}}{b}\Big)^2|\psi^\lambda_g(s)|^2+ 
 \mathcal{V}\circ \psi^\lambda_g(s) \Big) f^{\alpha,\gamma}(t-s)ds.
 \end{align*}
 
 Let us now denote by $\mathcal{I}$ the largest interval on which local solutions can be defined. In view of the above, $\mathcal I = (0,T^\lambda_g)$ is an open interval, and there exists a function $\psi_{g}^\lambda$ on $\mathcal I$ such that for all $T<T_g^\infty$ the function $\psi_g^\lambda$ is local solution on $(0,T]$.  Since any noncontinuable solution is a global solution, due to Lemma \ref{psi nonnegative}, either $T_g^\lambda = +\infty$ or $T_g^\lambda < \infty$ and $\|\psi_g^\lambda\|_{L^\infty_{T_g^\lambda,\alpha}}<\infty$. In the latter case $\psi^\lambda_g$ is well-defined on $(0,T^\lambda_g ]$ and so the function $\psi^\lambda_g$ can be extended to a larger interval, which contradicts the definition of $T^\lambda_g$. 

 \medskip 
 
 {\bf Step 3: Uniqueness of global solution.} This follows from the uniqueness of local solutions as any global solution is a local solution when restricted to finite time intervals.  \qed

  \subsection{Proof of Theorem~\ref{thmcf}(2)}

  In this section we determine the Laplace functional of the solutions $V_*$ to the equation \eqref{eq1} and link them with the unique global solution $\psi_{g}^\lambda$ to the nonlinear Volterra equation (\ref{PSI}).  For convenience, we rewrite this equation as 
  \beqlb\label{eqn.700}
  \psi^\lambda_{g} (t)=\lambda\cdot f^{\alpha,\gamma}(t)+(g-\varphi)*f^{\alpha,\gamma}(t) 
  \quad \mbox{with}\quad 
  \varphi(t):=  \frac{1}{2} \cdot \Big(\frac{\zeta_*^\mathtt{m}\sqrt{\lambda^\mathtt{m}_*}}{b}\Big)^2|\psi^\lambda_{g} (t)|^2 + \mathcal{V}\circ\psi^\lambda_{g} (t) .
  \eeqlb

  
  Our proof used a series of auxiliary results whose proofs will be given below. The key is to construct an exponential martingale associated with the pair $(V_*,\psi^\lambda_{g}) $. To this end, we introduce the semimartingale  
  \beqlb\label{eqn.705}
  Z_T(t):= \mathbf{E}\Big[ \lambda \cdot V_*(T) + (g -\varphi)*V_*(T)  \,\big|\, \mathscr{F}_t \Big] + \int_0^t\varphi(T-s)V_*(s)ds, \quad t \in [0,T]. 
  \eeqlb
  
  \begin{lemma}\label{Lemma.701}
  	In terms of the function $L_K$ defined in \eqref{func:K} the semimartingale $Z_T$ admits the following  alternative representation:
  	\beqlb\label{eqn.701}
  	Z_T(t)\ar=\ar 
  	V_*(0)\cdot L_K*\psi^\lambda_{g}(T)+\frac{a}{b} \int_0^T \psi^\lambda_{g}(s)ds +	\int_0^t\varphi(T-s)V_*(s)ds\cr 
  	\ar\ar  + \int_0^t \frac{\zeta_*^\mathtt{m}\sqrt{\lambda^\mathtt{m}_*}}{b} \cdot\psi^\lambda_{g}(T-s)\sqrt{V_*(s)}dB(s) \cr 
  	\ar\ar +\int_0^t\int_0^{\infty}\int_0^{V_*(s)}\Big(\int_{(T-s-y)^+}^{T-s}
  	\frac{\zeta_*^\mathtt{l}}{b}\cdot \psi^\lambda_{g}(r)dr\Big)\widetilde{N}(ds,dy,dz), \quad t \in [0,T]. 	
  	\eeqlb 
  	In particular, 
  	\beqlb\label{eqn.704}
  	Z_T(T)=\lambda \cdot V_*(T) + g*V_*(T)
  	\quad \mbox{and}\quad 
  	Z_T(0)=  V_*(0)\cdot L_K*\psi^\lambda_{g}(T)+\frac{a}{b} \int_0^T \psi^\lambda_{g}(s)ds .
  	\eeqlb
  \end{lemma}
  
 Applying It\^o's formula (see Theorem~5.1 in \cite[p.66]{IkedaWatanabe1989}) to $e^{-Z_T(t)}$ and using the representation (\ref{eqn.701}) and the second equality in (\ref{eqn.700}) we obtain that 
  \beqlb \label{eqn.703}
  e^{-Z_T(t)} \ar=\ar e^{-Z_T(0)} - \int_0^t e^{-Z_T(s)} \frac{\zeta_*^\mathtt{m}\sqrt{\lambda^\mathtt{m}_*}}{b} \cdot\psi^\lambda_{g}(T-s)\sqrt{V_*(s)}dB(s)\cr 
  \ar\ar + \int_0^t \int_0^{\infty}\int_0^{V_*(s)} e^{-Z_T(s)} \bigg( \exp\Big\{ -\int_{(T-s-y)^+}^{T-s}
  \frac{\zeta_*^\mathtt{l}}{b}\cdot \psi^\lambda_{g}(r)dr \Big\}-1 \bigg)  \widetilde{N}(ds,dy,dz).
  \eeqlb
  
  \begin{lemma}\label{Lemma.702}
  	The process $\{ e^{-Z_T(t)}: t\in[0,T] \}$ is a true $(\mathscr{F}_t)$-martingale. 
  \end{lemma}
  
Plugging the equation (\ref{eqn.704}) into (\ref{eqn.703}) and then taking expectations on both sides of the equation yields the desired representation (\ref{CF}). 
 It hence remains to prove the preceding two lemmas.
  
  \medskip
  \textbf{\textit{Proof of Lemma~\ref{Lemma.701}.}} 
 Let us denote by $Y_T(t)$ the conditional expectation in (\ref{eqn.705}), i.e.,
  \beqlb\label{eqn.711}
  Y_T(t)\ar:=\ar\mathbf{E}\Big[ \lambda \cdot V_*(T) + (g - \varphi)*V_*(T)  \,\big|\, \mathscr{F}_t \Big] \cr 
  \ar=\ar \lambda \cdot \mathbf{E}\big[V_*(T)\,\big|\, \mathscr{F}_t\big] 
  + \int_0^T(g - \varphi)(T-r)\mathbf{E}\big[V_*(r)\,\big|\, \mathscr{F}_t\big]dr.
  \eeqlb
  For $0\leq r,t\leq T$, taking conditional expectations on both side of (\ref{eq1}) yields that 
  \beqnn
  \mathbf{E}\big[V_*(r)\,\big|\, \mathscr{F}_t\big] 
  \ar=\ar  V_*(0) \cdot \big(1-F^{\alpha,\gamma}(r)\big) + \frac{a}{b}\cdot F^{\alpha,\gamma}(r)
  + \int_0^{r\wedge t} f^{\alpha,\gamma}(t-s)\cdot\frac{\zeta_*^\mathtt m\sqrt{\lambda^\mathtt{m}_*}}{b}\sqrt{V_*(s)}dB(s) \cr 
  \ar\ar +\int_0^{r\wedge t}\int_0^{\infty}\int_0^{V_*(s)}\Big(\int_{(t-s-y)^+}^{t-s} \frac{\zeta_*^\mathtt{l}}{b}\cdot f^{\alpha,\gamma}(r)dr \Big) \widetilde{N}(ds,dy,dz).
  \eeqnn
  Plugging this back into the right-hand side of the second equality in (\ref{eqn.711}) shows that
  \beqnn
  Y_T(t)  
  \ar=\ar Y_0 +Y_1(t) +Y_2(t),
  \eeqnn
  where 
  \beqlb
  Y_0 \ar:=\ar V_*(0) \cdot \lambda\cdot \big(1-F^{\alpha,\gamma}(T)\big) + V_*(0) \cdot\int_0^T(g - \varphi)(T-r)   \big(1-F^{\alpha,\gamma}(r)\big) dr \cr
  \ar\ar +  \frac{a}{b}\cdot \lambda \cdot F^{\alpha,\gamma}(T) +  \frac{a}{b}\cdot \int_0^T(g - \varphi)(T-r)   F^{\alpha,\gamma}(r) dr, \label{eqn.Y1} \\
  Y_1(t)\ar:=\ar \int_0^{t} \lambda f^{\alpha,\gamma}(T-s)\cdot\frac{\zeta_*^\mathtt m\sqrt{\lambda^\mathtt{m}_*}}{b}\sqrt{V_*(s)}dB(s)  \cr
  \ar\ar 
  +\int_0^T(g - \varphi)(T-r) \int_0^{r\wedge t} f^{\alpha,\gamma}(t-s)\cdot\frac{\zeta_*^\mathtt m\sqrt{\lambda^\mathtt{m}_*}}{b}\sqrt{V_*(s)}dB(s)dr, \label{eqn.Y2} \\
  Y_2(t)\ar:=\ar \int_0^{t}\int_0^{\infty}\int_0^{V_*(s)} \Big(\lambda\int_{(T-s-y)^+}^{T-s} \frac{\zeta_*^\mathtt{l}}{b}\cdot f^{\alpha,\gamma}(r)dr \Big) \widetilde{N}(ds,dy,dz)\cr
  \ar\ar + \int_0^T(g - \varphi)(T-r)\int_0^{r\wedge t}\int_0^{\infty}\int_0^{V_*(s)}\Big(\int_{(t-s-y)^+}^{t-s} \frac{\zeta_*^\mathtt{l}}{b}\cdot f^{\alpha,\gamma}(r)dr \Big) \widetilde{N}(ds,dy,dz)dr.  \label{eqn.Y3} 
  \eeqlb
  
  The representation (\ref{eqn.701}) can now be obtained by bringing the above quantities into a more convenient form. We start with $Y_0$. 
  
  Convoluting both sides of the second equation in (\ref{LK}) with function $L_K$ and then using the first inequality shows that 
  \beqnn
  1=L_K* f^{\alpha,\gamma}(t)+ F^{\alpha,\gamma}(t) 
  \quad \mbox{and hence} \quad 
  1-F^{\alpha,\gamma}(t)=L_K* f^{\alpha,\gamma}(t),\quad t\geq 0. 
  \eeqnn
  Plugging this into the right side of equation (\ref{eqn.Y1}) and then using the first equality in (\ref{eqn.700}), we have that
  \beqnn
  Y_0 \ar=\ar V_*(0) \cdot L_K* \Big(f^{\alpha,\gamma}(T) +  (g - \varphi)*   f^{\alpha,\gamma}(T) \Big)\cr
  \ar\ar
  +  \frac{a}{b}\cdot \int_0^T \Big( \lambda \cdot f^{\alpha,\gamma}(s) +   (g - \varphi)*f^{\alpha,\gamma}(s)\Big) dr \cr
  \ar\ar\cr
  \ar =\ar V_*(0) \cdot L_K* \psi^\lambda_{g} (T). 
  \eeqnn
  
  Now now turn to the process $Y_1$. Applying the stochastic Fubini theorem to the second term on the right side of equation (\ref{eqn.Y2}), this term equals
  \beqnn
  \int_0^t (g - \varphi)* f^{\alpha,\gamma}(T-s) \cdot\frac{\zeta_*^\mathtt m\sqrt{\lambda^\mathtt{m}_*}}{b}\sqrt{V_*(s)}dB(s) .
  \eeqnn
  Taking this back into (\ref{eqn.Y2}), merging it with the first term and then using  the first equality in (\ref{eqn.700}) shows that 
  \beqnn
  Y_1(t)
  \ar=\ar \int_0^{t} \big( \lambda f^{\alpha,\gamma}(T-s)+ (g - \varphi)* f^{\alpha,\gamma}(T-s) \big)\cdot\frac{\zeta_*^\mathtt m\sqrt{\lambda^\mathtt{m}_*}}{b}\sqrt{V_*(s)}dB(s)\cr
  \ar=\ar \int_0^{t} \psi^\lambda_{g}(T-s)\cdot\frac{\zeta_*^\mathtt m\sqrt{\lambda^\mathtt{m}_*}}{b}\sqrt{V_*(s)}dB(s) .
  \eeqnn
  
  We finally consider the process $Y_2$. Similarly as in \cite[Proposition~6.7]{Xu2024b}, applying the stochastic Fubini theorem (see Theorem~D2 in \cite{Xu2024b}) to the second term on the right side of (\ref{eqn.Y3}) and then using Fubini's theorem to the integrand, this term equals 
  \beqnn
  \lefteqn{\int_0^t\int_0^{\infty}\int_0^{V_*(s)}\Big(\int_0^{T-s}(g - \varphi)(T-s-r) \int_{(r-y)^+}^{r} \frac{\zeta_*^\mathtt{l}}{b}\cdot f^{\alpha,\gamma}(x)dx dr \Big) \widetilde{N}(ds,dy,dz)}\qquad \ar\ar\cr
  \ar=\ar \int_0^t\int_0^{\infty}\int_0^{V_*(s)}\Big(\frac{\zeta_*^\mathtt{l}}{b}\cdot\int_{(T-s-y)^+}^{T-s}(g - \varphi)*  f^{\alpha,\gamma}(r) dr \Big) \widetilde{N}(ds,dy,dz) . 
  \eeqnn
  Taking this back into (\ref{eqn.Y3}), merging this term with the first term and then using the first equality in (\ref{eqn.700}), we arrive at 
  \beqnn
  Y_2(t)\ar=\ar  \int_0^t\int_0^{\infty}\int_0^{V_*(s)}\Big(\frac{\zeta_*^\mathtt{l}}{b}\cdot\int_{(T-s-y)^+}^{T-s} \big(\lambda f^{\alpha,\gamma}(r) + (g - \varphi)*  f^{\alpha,\gamma}(r)\big) dr   \Big) \widetilde{N}(ds,dy,dz)\cr
  \ar=\ar \int_0^t\int_0^{\infty}\int_0^{V_*(s)}\Big(\frac{\zeta_*^\mathtt{l}}{b}\cdot\int_{(T-s-y)^+}^{T-s}  \psi^\lambda_{g}(r)  dr   \Big) \widetilde{N}(ds,dy,dz).  
  \eeqnn
  \qed

  \medskip
  \textbf{\textit{Proof of Lemma~\ref{Lemma.702}.}} 
  We proceed in three steps. 
  In the first step we show that the martingale property follows from the boundedness of the running maximum of the weak solution $V_*$ under a suitable equivalent measure. 
  In the second step we establish a representation of the process $V_*$ under the equivalent measure that allows us to establish the desired boundedness in a third step.

  \medskip
  {\it Step 1: Measure change.} We consider a martingale $\{U_T(t):t\in[0,T]\}$ with $U_T(t):= M_1(t)+ M_2(t)$ with  
  \beqnn
  M_1(t) \ar:=\ar -\int_0^t \frac{\zeta_*^\mathtt{m}\sqrt{\lambda^\mathtt{m}_*}}{b} \cdot\psi^\lambda_{g}(T-s)\sqrt{V_*(s)}dB(s),\cr 
  M_2(t) \ar:=\ar \int_0^t \int_0^{\infty}\int_0^{V_*(s)}   \bigg( \exp\Big\{ -\int_{(T-s-y)^+}^{T-s}
  \frac{\zeta_*^\mathtt{l}}{b}\cdot \psi^\lambda_{g}(r)dr \Big\}-1 \bigg)  \widetilde{N}(ds,dy,dz).
  \eeqnn
  Let $\mathcal{E}_{U_T}$ be the Dol\'ean-Dade exponential of $U_T$ defined by the unique solution to the SDE
  \beqnn
  \mathcal{E}_{U_T}(t) = 1 + \int_0^t \mathcal{E}_{U_T}(s)dU_T(s),\quad t\in[0,T].
  \eeqnn
  By It\^o's formula (see Theorem~5.1 in \cite[p.66]{IkedaWatanabe1989} or Theorem~37 \cite[p.84]{Protter2005}) and using the representation \eqref{eqn.700}, 
  \beqnn
  \mathcal{E}_{U_T}(t)
  \ar=\ar \exp\bigg\{- \int_0^t \varphi (T-s)V_*(s)ds   - \int_0^t \frac{\zeta_*^\mathtt{m}\sqrt{\lambda^\mathtt{m}_*}}{b} \cdot\psi^\lambda_{g}(T-s)\sqrt{V_*(s)}dB(s) \cr  
  \ar\ar - \int_0^t \int_0^{\infty}\int_0^{V_*(s)}  
  \Big( \int_{(T-s-y)^+}^{T-s}
  \frac{\zeta_*^\mathtt{l}}{b}\cdot \psi^\lambda_{g}(r)dr   \Big)  \widetilde{N}(ds,dy,dz) \bigg\},
  \eeqnn
  from which we see that 
  \beqnn
  e^{-Z_T(t)}= e^{-Z_T(0)}\cdot  \mathcal{E}_{U_T}(t), \quad t\in[0,T]. 
  \eeqnn
  
  We hence need to show that the local martingale $\mathcal{E}_{U_T}$ is a true martingale. Since $\mathcal{E}_{U_T}$ is a supermartingale it suffices to show that  
  \beqnn
  \mathbf{E}\big[\mathcal{E}_{U_T}(T_0)\big]=1,\quad\text{for all }T_0 \in [0,T].
  \eeqnn
  To this end, we introduce the localizing sequence of stopping times
  \beqnn
  \tau_n := \{t\geq 0: V_*(t)\geq n\}\wedge T_0, \quad \text{for}\quad T_0\in [0,T],
  \eeqnn
  so that the stopped process $ \mathcal{E}^n_{U_T}(t):=\mathcal{E}_{U_T}(\tau_n\wedge t) $ for $t\in[0,T]$ is a martingale for each $n\in\mathbb{N}$. Thus,
  \beqnn
  1=\mathbf{E}\big[\mathcal{E}^n_{U_T}(T_0)\big]
  \ar=\ar
  \mathbf{E}\big[\mathcal{E}^n_{U_T}(T_0);\tau_n\geq T_0\big]+\mathbf{E}\big[\mathcal{E}^n_{U_T}(T_0);\tau_n<T_0\big]\\
  \ar=\ar \mathbf{E}\big[\mathcal{E}_{U_T}(T_0);\tau_n\geq T_0\big]+\mathbf{E}\big[\mathcal{E}^n_{U_T}(T_0);\tau_n<T_0\big].
  \eeqnn
  Since  $\tau_n\wedge T_0\to T_0$ for $n\to\infty$, it follows from the monotone convergence theorem that
  \[
  \lim_{n\to\infty}\mathbf{E}\big[\mathcal{E}_{U_T}(T_0);\tau_n\geq T_0\big]=\mathbf{E}\big[\mathcal{E}_{U_T}(T_0)\big].
  \] 
  Therefore, it suffices to show that as $n\to\infty$,
  \beqlb \label{eqn.801}
  \mathbf{E}\big[\mathcal{E}^n_{U_T}(T_0);\tau_n<T_0\big]\to 0 .
  \eeqlb
  
  To prove this, we introduce the following probability measure $\mathbf{Q}^n$ associated with  the martingale $ \mathcal{E}^n_{U_T}$:
  \beqnn
  \frac{d \mathbf{Q}^n}{d \mathbf{P}}  =  \mathcal{E}^n_{U_T}(T_0) 
  \quad \mbox{and hence}\quad
  \frac{d \mathbf{Q}^n}{d \mathbf{P}} \Big|_{\mathscr{F}_t} =  \mathcal{E}^n_{U_T}(t),\quad t\in[0,T_0]. 
  \eeqnn
  We write $\mathbf{E}^{\mathbf{Q}^n}$ for the expectation under the law $\mathbf{Q}^n$. By Chebyshev's inequality,  
  \beqnn
  \mathbf{E}\big[\mathcal{E}^n_{U_T}(T_0);\tau_n<T_0\big] =  \mathbf{Q}^n \big\{  \tau_n<T_0 \big\}
  = \mathbf{Q}^n \bigg\{  \sup_{t\in[0,T_0]}V_*(t)\geq n \bigg\}\leq \frac{1}{n} \cdot \mathbf{E}^{\mathbf{Q}^n }\bigg[ \sup_{t\in[0,T_0]}V_*(t) \bigg].
  \eeqnn
  To obtain (\ref{eqn.801}) it thus suffices to show that 
  \begin{equation} \label{boundQ2}
  	\sup_{n\geq 1} \mathbf{E}^{\mathbf{Q}^n }\bigg[ \sup_{t\in[0,T_0]}V_*(t) \bigg] < \infty.
  \end{equation}
  
  \medskip
  
  {\it Step 2: Representation under $\mathbf{Q}^n$.}
  To prove the desired boundedness result, we first establish an equivalent representation of the stochastic equation (\ref{eq1}) under the law $\mathbf{Q}^n$. 
  
  By Girsanov's theorem (see Theorem~3.11 and 3.17 in \cite[p.168, 170]{JacodShiryaev2003})  we have that under the law $\mathbf{Q}^n$, the Brownian motion
  $B$ is continuous martingale with quadratic variation $\langle B\rangle_t=   t\wedge \tau_n $ and the Poisson random measure $N(ds,dy,dz)$ turns to be a random point measure on $(0,\infty)^3$ with intensity 
  \beqnn
  \mathbf{1}_{\{  s\leq \tau_n\}} \cdot \exp\bigg\{ -\int_{(T-s-y)^+}^{T-s}
  \frac{\zeta_*^\mathtt{l}}{b}\cdot \psi^\lambda_{g}(r)dr \bigg\} ds \cdot \nu_*(dy)\cdot dz.
  \eeqnn
  For each $t_1\in[0,T]$, we now define the following auxiliary process 
  \begin{align}    \label{eq1011}
  	\begin{split}
  		V_{*,t_1}(t) :=  & ~ V_*(0) \cdot \big(1-F^{\alpha,\gamma}(t)\big) + \frac{a}{b}\cdot F^{\alpha,\gamma}(t)
  		+ \int_0^t f^{\alpha,\gamma}(t_1-s)\cdot\frac{\zeta_*^\mathtt m\sqrt{\lambda^\mathtt{m}_*}}{b}\sqrt{V_{*}(s)}dB(s)\\  
  		&+\int_0^t\int_0^{\infty}\int_0^{V_{*}(s)}\Big(\int_{(t_1-s-y)^+}^{t_1-s} \frac{\zeta_*^\mathtt{l}}{b}\cdot f^{\alpha,\gamma}(r)dr \Big) \widetilde{N}(ds,dy,dz),
  	\end{split}
  \end{align}
  which is a $(\mathscr{F}_t)$-semimartingale under $\mathbf{P}$ and 
 \beqnn
  	V_{*,t_1}(t_1)\overset{\rm a.s.}=V_{*}(t_1).
 \eeqnn 
 Applying  Theorem~3.24 in \cite[p.172]{JacodShiryaev2003} to the process $V_{*,t_1}$ we see that 
under the law $\mathbf{Q}^n$ this process is a strong solution to the stochastic equation 
  \beqnn
  \lefteqn{V_{*,t_1}(t)
  	= V_*(0) \cdot \big(1-F^{\alpha,\gamma}(t)\big) + \frac{a}{b}\cdot F^{\alpha,\gamma}(t)   - \frac{|\zeta_*^\mathtt m|^2\lambda^\mathtt{m}_*}{b^2}\int_0^t \big|f^{\alpha,\gamma}(t_1-s)\big|^2\cdot \mathbf{1}_{\{  s\leq \tau_n\}} \cdot V_{*}(s)ds } \ar\ar\cr
  \ar\ar  +\int_0^t \int_0^\infty \int_{(t_1-s-y)^+}^{t_1-s} \frac{\zeta_*^\mathtt{l}}{b}\cdot f^{\alpha,\gamma}(r)dr   \Big( \exp\Big\{ -\int_{(T-s-y)^+}^{T-s} \frac{\zeta_*^\mathtt{l}}{b}\cdot \psi^\lambda_{g}(r)dr \Big\}-1 \Big) \nu_*(dy) \cdot \mathbf{1}_{\{  s\leq \tau_n\}} \cdot V_{*}(s)ds\cr
  \ar\ar + \int_0^t f^{\alpha,\gamma}(t_1-s)\cdot\frac{\zeta_*^\mathtt m\sqrt{\lambda^\mathtt{m}_*}}{b}\sqrt{V_{*}(s)}dB(s)+\int_0^t\int_0^{\infty}\int_0^{V_{*}(s)}\Big(\int_{(t_1-s-y)^+}^{t_1-s} \frac{\zeta_*^\mathtt{l}}{b}\cdot f^{\alpha,\gamma}(r)dr \Big) \widetilde{N}(ds,dy,dz).
  \eeqnn
  Setting $t=t_1$ and using that $t_1\in [0,T]$ is arbitrary, the stochastic equation (\ref{eq1})  under $\mathbf{P}$ is equal (in law) to the following stochastic equation under $\mathbf{Q}^n$: 
  \beqlb \label{eqn.802}
  V_{*}(t)
  = V_*(0) \cdot \big(1-F^{\alpha,\gamma}(t)\big) + \frac{a}{b}\cdot F^{\alpha,\gamma}(t)  + A^n_1(t)   + A^n_2(t) + M^n_1(t) + M^n_2(t) , 
  \eeqlb
  with 
  \beqnn
  A^n_1(t) \ar:=\ar  - \frac{|\zeta_*^\mathtt m|^2\lambda^\mathtt{m}_*}{b^2}\int_0^t \big|f^{\alpha,\gamma}(t-s)\big|^2\cdot \mathbf{1}_{\{  s\leq \tau_n\}} \cdot V_{*}(s)ds,  \cr
  A^n_2(t) \ar:=\ar \int_0^t  \mathbf{1}_{\{  s\leq \tau_n\}} \cdot V_{*}(s)ds \int_0^\infty  \int_{(t-s-y)^+}^{t-s} \frac{\zeta_*^\mathtt{l}}{b}\cdot f^{\alpha,\gamma}(r)dr\cr
  \ar\ar \qquad \times    \Big( \exp\Big\{ -\int_{(T-s-y)^+}^{T-s} \frac{\zeta_*^\mathtt{l}}{b}\cdot \psi^\lambda_{g}(r)dr \Big\}-1 \Big) \nu_*(dy), \cr
  M^n_1(t)\ar:=\ar  \int_0^t f^{\alpha,\gamma}(t-s)\cdot\frac{\zeta_*^\mathtt m\sqrt{\lambda^\mathtt{m}_*}}{b}\sqrt{V_{*}(s)}dB(s) , \cr
  M^n_2(t)\ar:=\ar \int_0^t\int_0^{\infty}\int_0^{V_{*}(s)}\Big(\int_{(t-s-y)^+}^{t-s} \frac{\zeta_*^\mathtt{l}}{b}\cdot f^{\alpha,\gamma}(r)dr \Big) \widetilde{N}(ds,dy,dz).
  \eeqnn

  \medskip
  
  {\it Step 3: Moment estimates under $\mathbf{Q}^n$.} 
  We now prove the upper bound in \eqref{boundQ2}. 
  For some constant $\theta>0$ to be specified later, we multiply both sides of the equation (\ref{eqn.802}) by $e^{-\theta t}$ and then consider their respective suprema over $[0,T]$ to obtain that
  \beqlb\label{eqn.7001}
  \sup_{t\in[0,T]} e^{-\theta t} V_*(t) 
  \ar\leq\ar \Big(V_*(0) + \frac{a}{b}\Big) + \sup_{t\in[0,T]} \big| e^{-\theta t} A^n_1(t)\big| + \sup_{t\in[0,T]} \big| e^{-\theta t} A^n_2(t)\big|\cr
  \ar\ar  +\sup_{t\in[0,T]} \big|   M^n_1(t)\big|+ \sup_{t\in[0,T]} \big| M^n_2(t)\big|.
  \eeqlb
  
  A simple calculation along with \eqref{upper bound of mittag leffler density functions} shows that one can find a large enough constant $\theta>0$ that does not depend  of $n$ such that almost surely
  \beqlb\label{eqn.7002}
  \sup_{t\in[0,T]} \big| e^{-\theta t} A^n_1(t)\big| 
  \ar\leq\ar \sup_{t\in[0,T]} e^{-\theta t}V_*(t)\cdot \frac{|\zeta_*^\mathtt m|^2\lambda^\mathtt{m}_*}{b^2} \cdot \int_0^T e^{-\theta s} \big|f^{\alpha,\gamma}(s)\big|^2ds 
  \leq \frac{1}{4}\cdot \sup_{t\in[0,T]} e^{-\theta t} V_*(t).
  \eeqlb 
  Similarly,since $|e^{-x}-1|\leq x$ uniformly in $x\geq 0$, we also have that 
  \beqnn
  \sup_{t\in[0,T]} \big| e^{-\theta t} A^n_2(t)\big| 
  \ar\leq\ar \sup_{t\in[0,T]} e^{-\theta t}V_*(t)\cdot \sup_{t\in[0,T]} \int_0^t  e^{-\theta (t-s)} ds \int_0^\infty  \int_{(t-s-y)^+}^{t-s} \frac{\zeta_*^\mathtt{l}}{b}\cdot f^{\alpha,\gamma}(r)dr\cr
  \ar\ar \qquad \times    \Big| \exp\Big\{ -\int_{(T-s-y)^+}^{T-s} \frac{\zeta_*^\mathtt{l}}{b}\cdot \psi^\lambda_{g}(r)dr \Big\}-1 \Big| \nu_*(dy)\cr
  \ar\leq\ar \sup_{t\in[0,T]} e^{-\theta t}V_*(t)\cdot  \sup_{t\in[0,T]} \int_0^t  e^{-\theta (t-s)}  ds \int_0^\infty  \int_{(t-s-y)^+}^{t-s} \frac{\zeta_*^\mathtt{l}}{b}\cdot f^{\alpha,\gamma}(r)dr\cr
  \ar\ar \qquad \times    \int_{(T-s-y)^+}^{T-s} \frac{\zeta_*^\mathtt{l}}{b}\cdot \psi^\lambda_{g}(r)dr \nu_*(dy)\cr
  \ar\leq\ar \sup_{t\in[0,T]} e^{-\theta t}V_*(t)\cdot \bigg( \sup_{t\in[0,T]} \int_0^t   e^{-\theta (t-s)} ds \int_0^\infty  \Big(\int_{(t-s-y)^+}^{t-s} \frac{\zeta_*^\mathtt{l}}{b}\cdot f^{\alpha,\gamma}(r)dr\Big)^2\nu_*(dy)\cr
  \ar\ar \qquad +  \sup_{t\in[0,T]} \int_0^t  e^{-\theta (t-s)}  ds \int_0^\infty  \Big(\int_{(T-s-y)^+}^{T-s} \frac{\zeta_*^\mathtt{l}}{b}\cdot \psi^\lambda_{g}(r)dr\Big)^2 \nu_*(dy)\bigg)\cr
  \ar=\ar \sup_{t\in[0,T]} e^{-\theta t}V_*(t)\cdot \bigg(  \int_0^T   e^{-\theta s} ds \int_0^\infty  \Big(\int_{(s-y)^+}^{s} \frac{\zeta_*^\mathtt{l}}{b}\cdot f^{\alpha,\gamma}(r)dr\Big)^2\nu_*(dy)\cr
  \ar\ar \qquad +   \sup_{t\in[0,T]} \int_0^t  e^{-\theta (t-s)}  ds  \int_0^\infty \Big(\int_{(T-s-y)^+}^{T-s} \frac{\zeta_*^\mathtt{l}}{b}\cdot \psi^\lambda_{g}(r)dr\Big)^2 \nu_*(dy)\bigg)\cr
  \ar=\ar \sup_{t\in[0,T]} e^{-\theta t}V_*(t)\cdot \bigg(  \int_0^T   e^{-\theta s}  ds \int_0^\infty  \Big(\int_{(s-y)^+}^{s} \frac{\zeta_*^\mathtt{l}}{b}\cdot f^{\alpha,\gamma}(r)dr\Big)^2\nu_*(dy)\cr
  \ar\ar \qquad +   \sup_{t\in[0,T]} \int_0^t  e^{-\theta (t-s)}  ds  \int_0^\infty \Big(\int_{(T-s-y)^+}^{T-s} \frac{\zeta_*^\mathtt{l}}{b}\cdot \psi^\lambda_{g}(r)dr\Big)^2 \nu_*(dy)\bigg).
  \eeqnn
  By \eqref{eqn.210}, we have as $\theta \to\infty$,
  \beqnn
  \int_0^T   e^{-\theta s}  ds \int_0^\infty  \Big(\int_{(s-y)^+}^{s} \frac{\zeta_*^\mathtt{l}}{b}\cdot f^{\alpha,\gamma}(r)dr\Big)^2\nu_*(dy)
  \to 0.
  \eeqnn
  Moreover, by \eqref{upper bound of mittag leffler density functions} we also have as $\theta\to\infty$,
  \beqnn
  \lefteqn{\sup_{t\in[0,T]} \int_0^t  e^{-\theta (t-s)}  ds   \Big(\int_{(T-s-y)^+}^{T-s} \frac{\zeta_*^\mathtt{l}}{b}\cdot \psi^\lambda_{g}(r)dr\Big)^2 \nu_*(dy) }\ar\ar\cr
  \ar\leq\ar  C\cdot \sup_{t\in[0,T]} \int_0^t  e^{-\theta (t-s)}  ds   \Big(
  \big(|(T-s-y)^+|^{\alpha-1}\cdot y\big) \wedge (T-s)^{\alpha}\wedge y^\alpha \Big)^2 \nu_*(dy)\cr
  \ar\to\ar 0
  \eeqnn
  As a result, there exists a constant $\theta>0$ that is independent of $n$ such that 
  \beqlb\label{eqn.7003}
  \sup_{t\in[0,T]} \big| e^{-\theta t} A^n_2(t)\big| 
  \ar\leq\ar \frac{1}{4}\cdot \sup_{t\in[0,T]} e^{-\theta t}V_*(t).
  \eeqlb
  
  Taking \eqref{eqn.7002} and \eqref{eqn.7003} back into the right side of \eqref{eqn.7001} and then using the fact $F^{\alpha,\gamma}(t)\leq 1$ we see that 
  \beqnn
  \sup_{t\in[0,T]} e^{-\theta t} V_*(t) 
  \ar\leq\ar 2\cdot \Big( \Big(V_*(0) + \frac{a}{b}\Big) +  \sup_{t\in[0,T]} \big|   M^n_1(t)\big|+ \sup_{t\in[0,T]} \big| M^n_2(t)\big| \Big)
  \eeqnn  
  and hence that
  \beqnn
  \mathbf{E}^{\mathbf{Q}^n} \bigg[\sup_{t\in[0,T]} V_*(t) \bigg]
  \ar\leq\ar 2  e^{\theta T}\cdot \Big( \Big(V_*(0) + \frac{a}{b}\Big) +  \mathbf{E}^{\mathbf{Q}^n} \bigg[\sup_{t\in[0,T]} \big|   M^n_1(t)\big| \bigg]+ \mathbf{E}^{\mathbf{Q}^n} \bigg[\sup_{t\in[0,T]} \big| M^n_2(t)\big|\bigg] \Big).
  \eeqnn
  
  It remains to prove the boundedness of the last two expectations on the right side of this inequality.
  Similarly as in the proofs of Lemma~\ref{V*sup} and \ref{V2p}, we also have for any $p\geq 0$,
  \beqnn
  \sup_{n\geq 1} \sup_{t\in[0,T]} \mathbf{E}^{\mathbf{Q}^n}\Big[ \big| V_{*}(t) \big|^p \Big] \leq C\cdot (1+T)^{p\alpha} . 
  \eeqnn
  Armed with this moment estimate, the same arguments as given in the proofs of Lemma~\ref{V*sup} and \ref{holderJ2} show that for any $p\geq 0$ and $0\leq t_1,t_2\leq T$,
  \beqnn
  \sup_{n\geq 1}\mathbf{E}^{\mathbf{Q}^n} \Big[\big| M^n_1(t_2)-M^n_1(t_1) \big|^{2p} \Big] 
  \ar\leq\ar C\cdot (1+T)^{p\alpha} |t_2-t_1|^{p(2\alpha-1)} ,\cr
  \sup_{n\geq 1}\mathbf{E}^{\mathbf{Q}^n} \Big[\big| M^n_2(t_2)-M^n_2(t_1) \big|^{2p} \Big] 
  \ar\leq\ar C\cdot (1+T)^{p\alpha} |t_2-t_1|^{p\alpha} . 
  \eeqnn
  Following the argument as in the proof of Theorem~\ref{MainThm.02} we have 
  \beqnn
  \sup_{n\geq 1}\mathbf{E}^{\mathbf{Q}^n}\bigg[\sup_{t\in[0,T]} \big|   M^n_1(t)\big| \bigg]+  \sup_{n\geq 1} \mathbf{E}^{\mathbf{Q}^n}\bigg[\sup_{t\in[0,T]} \big| M^n_2(t)\big|\bigg]<\infty . 
  \eeqnn
  This yields the desired inequality \eqref{boundQ2}. 

%
%
%